\documentclass[11pt, reqno]{amsart}

\usepackage{amssymb}
\usepackage{hyperref} 
\usepackage{mathrsfs,bbold}
\usepackage{amsmath}
\usepackage{amsthm}
\usepackage{enumerate}
\usepackage{dsfont}
\usepackage{color}
\usepackage[a4paper]{geometry}
\usepackage[utf8]{inputenc}
\usepackage{stmaryrd}
\usepackage{titlesec}
\usepackage{mathabx}
\usepackage{appendix}
\usepackage{todonotes, fancyhdr}

\usepackage[all]{xy}
\renewcommand\thesection       {\arabic{section}}
\renewcommand\thesubsection    {\thesection{\boldmath $.$}\arabic{subsection}}
\renewcommand\thesubsubsection    {\thesection{\boldmath $.$}\arabic{subsection}{\boldmath $.$}\arabic{subsubsection}}

\newcommand{\ssk}{\smallskip}

\renewcommand{\epsilon}{\varepsilon}

\newcommand\bbR{\mathbb{R}}

\newcommand{\mcA}{\mathcal{A}}
\newcommand{\mcB}{\mathcal{B}} 
\newcommand{\mcC}{\mathcal{C}} 
\newcommand{\mcD}{\mathcal{D}}
\newcommand{\mcE}{\mathcal{E}}
\newcommand{\mcG}{\mathcal{G}}
\newcommand{\mcH}{\mathcal{H}}
\newcommand{\mcI}{\mathcal{I}}

\newcommand{\mcM}{\mathcal{M}}

\newcommand{\mcP}{\mathcal{P}}
\newcommand{\mcQ}{\mathcal{Q}}
\newcommand\mcS{\mathcal{S}}

\newcommand\bfT{\textsf{\textbf{T}}}
\newcommand\bfR{\textsf{\textbf{R}}}

\newcommand{\bfZ}{\textsf{\textbf{Z}}}

\newcommand{\scrC}{\ensuremath{\mathscr{C}}}

\newenvironment{Dem}[1][\unskip]{%
    \begin{list}{\hspace{0.5cm}{\sf \textbf{Proof #1 --}}}{%
        \setlength{\topsep}{0pt}%
        \setlength{\leftmargin}{0pt}%
        \setlength{\rightmargin}{0pt}%
        \setlength{\listparindent}{0pt}%
        \setlength{\itemindent}{0pt}%
        \setlength{\parsep}{0pt}%
        \addtolength{\leftmargin}{20pt}%
        \addtolength{\rightmargin}{0pt}%
    } \item }{\hfill $\rhd$\end{list}\smallskip}

\titleformat{\section}[block]
{\filcenter\normalfont\sffamily\bfseries\Large}
{{\hspace{-0.7cm}}\thesection \hspace{0.2em} --\vspace{0.3cm}}{0.5em}{}

\titleformat{\subsection}[block]
{\filcenter\normalfont\sffamily\bfseries\large}  						  
{\hspace{-0.7cm}\thesubsection \hspace{0.5em} \vspace{0.3cm}}{.5em}{}  
\titlespacing{\subsection}{-0pc}{1.5ex plus .1ex minus .2ex}{0pc}

\titleformat{\subsubsection}[block]
{\filcenter\normalfont\sffamily\bfseries}					  
{\hspace{-0.7cm}\thesubsubsection \hspace{0.5em} \vspace{0.3cm}}{.5em}{}  
\titlespacing{\subsection}{-0pc}{1.5ex plus .1ex minus .2ex}{0pc}

\newtheorem{theorem}{Theorem}[section]

\newtheorem{definition}[theorem]{Definition}

\numberwithin{theorem}{section}
\numberwithin{equation}{section}

\newtheoremstyle{mystyle}
{3pt}               
{3pt}               
{\it }                      
{}                      
{\sffamily\bfseries}             
{}                      
{0.5em}                 
{#1 #2{\boldmath $.$}}

\theoremstyle{mystyle}

\newtheorem{thm}{Theorem}

\newtheorem{cor}[thm]{\hspace{-0.15cm}  {Corollary} }
\newtheorem{lem}[thm]{\hspace{-0.2cm}  {Lemma} }
\newtheorem{prop}[thm]{\hspace{-0.2cm} {Proposition}}
\newtheorem{defn}[thm]{ \hspace{-0.3cm} {Definition}}
\newtheorem*{defn*} {Definition}
\newtheorem*{prop*} {Proposition}
\newtheorem*{lem*} {Lemma}
\newtheorem*{cor*} {Corollary}
\newtheorem{rem}[thm]{\hspace{-0.15cm} {Remark}}

\newtheoremstyle{mystyle2}
{3pt}               
{3pt}               
{\it }                      
{}                      
{\sffamily\bfseries}             
{}                      
{0.5em}                 
{\llap{#2 }#1{\boldmath$.$}}

\theoremstyle{mystyle2}

\newtheorem*{definition*}{Definition}
\numberwithin{equation}{section} 

\newcommand\restr[2]{{
  \left.\kern-\nulldelimiterspace 
  #1 
  \vphantom{\big|} 
  \right|_{#2} 
  }}

\begin{document}
\date{\today}

\vspace*{3ex minus 1ex}
\begin{center}
{\Huge\sffamily{High order paracontrolled calculus  \vspace{0.5cm}}}
\end{center}

\begin{center}
{\sf I. BAILLEUL\footnote{I.Bailleul thanks the U.B.O. for their hospitality, part of this work was written there.} and F. BERNICOT\footnote{F. Bernicot's research is supported by the ERC Project FAnFArE no. 637510.   \\

{\sf AMS Classification:} 60H15; 35R60; 35R01 \vspace{0.1cm}

{\sf Keywords:} Stochastic singular PDEs; paracontrolled calculus; generalised Parabolic Anderson Model equation; generalised KPZ equation}}
\end{center}

\vspace{1cm}

\begin{center}
\begin{minipage}{0.8\textwidth}
\renewcommand\baselinestretch{0.7} \scriptsize \textbf{\textsf{\noindent Abstract.}} 
We develop in this work a general version of paracontrolled calculus that allows to treat analytically within this paradigm a whole class of singular partial differential equations with the same efficiency as regularity structures. This work deals with the analytic side of the story and offers a toolkit for the study of such equations, under the form of a number of continuity results for some operators, while emphasizing the simple and systematic mechanics of computations within paracontrolled calculus, via the introduction of two model operations ${\sf E}$ and $\sf F$. We illustrate the efficiency of this elementary approach on the example of the generalised parabolic Anderson model equation
$$
(\partial_t + L) u = f(u)\zeta,
$$
on a $3$-dimensional closed manifold, and the generalised KPZ equation
$$
(\partial_t + L) u = f(u)\zeta + g(u)(\partial u)^2,
$$
driven by a $(1+1)$-dimensional space/time white noise.
\end{minipage}
\end{center}

\vspace{1cm}

{\sf 
\begin{center}
\begin{minipage}[t]{11cm}
\baselineskip =0.35cm
{\scriptsize 

\center{\textbf{Contents}}

\vspace{0.1cm}

\textbf{ 1.~Paracontrolled calculus\dotfill 
\pageref{SectionParaControlledCalculus}}

\textbf{ 2.~High order paracontrolled expansion\dotfill 
\pageref{SectionTaylor}}

\textbf{ 3.~A toolkit for paracontrolled calculus
\dotfill \pageref{SectionToolKit}}

\textbf{ 4.~Nonlinear singular PDEs
\dotfill \pageref{SectionNonlinearPDEs}}

\textbf{ A.~Parabolic setting
\dotfill \pageref{SectionAppendix}}

\textbf{ B.~Proof of the high order paracontrolled expansion
\dotfill \pageref{SectionAppendixTaylor}}

\textbf{ C.~Continuity results
\dotfill \pageref{SectionAppendixContinuity}}}
\end{minipage}
\end{center}
}   \vspace{1cm}

\section[\hspace{0.7cm} Paracontrolled calculus]{Paracontrolled calculus}
\label{SectionParaControlledCalculus}

Thirty years after T. Lyons' seminal work on controlled differential equations \cite{Lyons98}, it is now well-understood that the construction of a robust approximation theory for continuous time stochastic systems, such as stochastic differential equations or stochastic partial differential equations, requires a twist in the notion of noise that allows to treat the resolution of such equations in a two step process.   \vspace{0.15cm}

\begin{enumerate}
   \item[{\sf (a)}] Enhance the noise $\zeta$ into an enriched object $\widehat{\zeta}$ that lives in some space of analytic objects -- this is a purely \textit{probabilistic} step.   \vspace{0.1cm}
   
   \item[{\sf (b)}] Given \textit{any} such object $\widehat{\zeta}$ in this space, one can introduce a Banach space $\mcS$ such that the equation makes sense for the unknown in the image of $\mcS$ by a simple $\widehat{\zeta}$-dependent map. The equation can be formulated in $\mcS$, and solved uniquely by a \textit{deterministic} analytic argument, such as the contraction principle, which gives the continuity of the solution as a function of $\widehat{\zeta}$.   \vspace{0.15cm}
\end{enumerate}

These two steps are very different in nature and require totally different tools. The present work deals with the deterministic side of the story, point {\sf (b)}, for the study of \textit{singular partial differential equations} (PDEs). The term \textit{singular} refers here to the fact that the 'noise' in the equation is not regular enough for all the  expressions in the equation to make sense analytically, given the expected regularity of the solution in terms of the regularity of the 'noise'. Recall that one can generically not make sense of the product of a distribution with a continuous function.

\bigskip

\subsection{Overview}
\label{SubsectionOverview}

Hairer's theory of regularity structures \cite{HairerRegularity} provides undoubtedly the most complete picture for the study of a whole class of singular stochastic PDEs from the above point of view -- the class of the so-called singular subcritical parabolic stochastic PDEs. It comes with a very rich algebraic structure and an entirely new setting that are required to give flesh to the guiding principle that a solution should be described by the datum at each point in space-time of its high order 'jet' in a basis given by the elements of the enhanced noise. Regularity structures are introduced as a tool for describing these jets. At the same time that Hairer built his theory, Gubinelli-Imkeller-Perkowski implemented in \cite{GIP} this idea of giving a local/global description of a possible solution in a different way, using the language of paraproducts and avoiding the introduction of any new setting, but providing only a first order description of the objects under study. This is what we shall call from now on the \textit{first order paracontrolled calculus}. While this kind of approach may seem far from being as powerful as Hairer's machinery, the first order paracontrolled approach to singular stochastic PDEs has been successful in recovering and extending a number of results that can be proved within the setting of regularity structures, on the parabolic Anderson model and Burgers equations \cite{GIP, BB15, BBF15, ChoukCannizzaro}, the KPZ equation \cite{KPZReloaded}, the scalar $\Phi^4_3$ equation \cite{CatellierChouk}, the stochastic Navier-Stokes equation \cite{ZZ1, ZZ2, ZZ3}, or the study of the continuous Anderson Hamiltonian \cite{ChoukAllez}, to name but a few.

\medskip

We develop in this work a high order version of paracontrolled calculus that allows to treat analytically within this paradigm some parabolic singular partial differential equations that are beyond the scope of the original formulation of the theory. We refer to our setting as \textit{paracontrolled calculus}. By a 'noise' in an equation we shall simply mean a function/distribution-valued parameter $\zeta$ -- realisations of a white noise are typical examples. Within our setting, and given as input a noise $\zeta$ and some initial condition, the resolution process of a typical parabolic equation 
\begin{equation}
\label{EqCanonicalSPDE}
\mathscr{L} u := (\partial_t + L) u = f(u,\partial u, \zeta),
\end{equation}
involves the following elementary steps. Write $\mathscr{L}^{-1} := (\partial_t + L)^{-1}$ for the resolution operator, and keep in mind that we have in hands two space-time paraproducts ${\sf P}$ and $\widetilde{\sf P}$, related by the intertwining relation 
$$
\mathscr{L}^{-1}\circ{\sf P} = \widetilde{\sf P}\circ\mathscr{L}^{-1};
$$
all the objects are properly introduced below.   \vspace{0.15cm}

\begin{enumerate}
   \item[\textsf{\textbf{1{\boldmath $.$}}}] \textsf{\textbf{Paracontrolled ansatz{\boldmath $.$}}} \textit{The irregularity of the noise $\zeta$, and the form of the equation, dictate the choice of a solution space made up of functions/distributions of the form
   \begin{equation}
   \label{EqAnsatzU}
   u = \sum_{i=1}^{k_0} \widetilde{\sf P}_{u_i}Z_i + u^\sharp,
   \end{equation}
   for reference functions/distributions $Z_i$ that depend formally only on $\zeta$, to be determined later. The 'derivatives' $u_i$ of $u$ also need to satisfy a similar a structure equation ; their derivatives as well, and so on. See Definition \ref{DefnPCSystem} below. One sees the above description \eqref{EqAnsatzU} of $u$ as a paracontrolled expansion at order $k_0$ for it. Denote by $\widehat{u}$ the datum of $u$ and all its derivatives, and by $\widehat{u}^\sharp$ the datum of all their derivatives.   }   \vspace{0.15cm}
   
   \item[\textsf{\textbf{2{\boldmath $.$}}}] \textsf{\textbf{Right hand side{\boldmath $.$}}} \textit{The use of a high order paracontrolled expansion formula, and a number of continuity results for some operators, allow together to rewrite the right hand side $f(u, \partial u, \zeta)$ of equation \eqref{EqCanonicalSPDE} in the canonical form
   \begin{equation}
   \label{EqRHS}
   f\big(u,\partial u, \zeta\big) = \sum_{j=1}^{k_0} {\sf P}_{v_j}Y_j + (\flat)    
   \end{equation}
where $(\flat)$ is a nice remainder and the distributions $Y_j$ depend only on $\zeta$ and the $Z_i$.}   \vspace{0.15cm}
   
   \item[\textsf{\textbf{3{\boldmath $.$}}}]  \textsf{\textbf{Fixed point{\boldmath $.$}}} \textit{Denote by $P$ the resolution of the free heat equation 
   $$
   Pu_0:=(\tau,x) \mapsto \big(e^{-\tau L}u_0\big)(x).
   $$ 
   Then the fixed point relation
   \begin{equation*}
   \begin{split}
   u &= Pu_0 + \mathscr{L}^{-1}\big(f(u, \partial u, \zeta)\big)   \\
      &= Pu_0 + \sum_{j=1}^{k_0} \mathscr{L}^{-1}\Big({\sf P}_{v_j}Y_j\Big) + \mathscr{L}^{-1}(\flat)   \\
      &= Pu_0 + \sum_{j=1}^{k_0}\widetilde{\sf P}_{v_j}Z_j + \mathscr{L}^{-1}(\flat),
   \end{split}
   \end{equation*} 
imposes some consistency relations on the choice of the $Z_i = \mathscr{L}^{-1}(Y_i)$ that define them uniquely as functions of $\zeta$, and induces a fixed point relation for $\widehat{u}$, or rather $\widehat{u}^\sharp$.}   \vspace{0.15cm}
\end{enumerate}   
   
See Proposition \ref{prop:derive} for a justification of the name 'derivative' for the $u_i$ in identity \eqref{EqAnsatzU}. The expressions inside the $Y_i$'s that have no proper analytical sense need to be given a priori as components of the enhanced distribution $\widehat{\zeta}$. As expected, the enhanced noise $\widehat{\zeta}$ contains what is needed to make sense of the corresponding ill-defined products in the regularity structures setting. We shall not touch in this work on renormalisation matters, and we invite the reader to read the latest developments of Hairer \& co-authors on the subject; see \cite{BrunedChandraChevyrevHairer, BrunedHairerZambotti, ChandraHairer}. In any case, we shall assume here that the enhancement $\widehat{\zeta}$ of $\zeta$ is given. 

\medskip

We single out here the notion of paracontrolled system involved in point \textbf{\textsf{1}} above, in the setting that will be sufficient to deal with the generalised (PAM) and (KPZ) equations studied in Section \ref{SectionNonlinearPDEs}. It corresponds to having $k_0=3$. Fix $2/5<\alpha<1/2$. For each $1\leq i\leq 3$, we denote below by $Z_i$ a finite collection $\big(Z_i^{(n_i)}\big)$ of spacetime functions of parabolic regularity $i\alpha$. Given a collection $\big(u_i^{(n_i)}\big)$ of bounded functions indexed by the same set $(n_i)$ as $Z_i$, we write $\widetilde{\sf P}_{u_i}Z_i$ for $\sum_{n_i} \widetilde{\sf P}_{u_i^{(n_i)}}Z_i^{(n_i)}$.

\medskip

\begin{defn}
\label{DefnPCSystem}
A third order \textbf{\textsf{paracontrolled system}} is a family 
$$
\widehat{u} := \big(u, u_i, u_{ij},u_{ijk}\big)_{1\leq i,i+j,i+j+k\leq 3}
$$ 
of collections of bounded parabolic functions defined by the datum of remainders 
$$
u^\sharp\in\mcC^{4\alpha}, \; u_i^\sharp\in\mcC^{(4-i)\alpha}, \; u_{ij}^\sharp\in\mcC^{(4-i-j)\alpha},\; u_{ijk}^\sharp\in\mcC^\alpha, \qquad 1\leq i, i+j,i+j+k\leq 3
$$
via the identity
\begin{equation*}
\begin{split}
&u = \sum_{i=1..3}\widetilde{\sf P}_{u_i}Z_i + u^\sharp,   \\
&u_i = \sum_{i+j=1..3} \widetilde{\sf P}_{u_{ij}}Z_j + u_i^\sharp,   \\
&u_{ij} = \sum_{i+j+k=1..3} \widetilde{\sf P}_{u_{ijk}}Z_k + u_{ij}^\sharp,   \\
&u_{ijk} = u_{ijk}^\sharp.
\end{split}
\end{equation*}
\end{defn}

\medskip

The above system is a shorthand notation for 
$$
u_i^{(n_i)} = \sum_{i+j=1..3} \widetilde{\sf P}_{u^{(n_i)}_{ij}}Z_j + \big(u_i^{(n_i)}\big)^\sharp,
$$
and so on. There is only one function of the form $u_{ijk}$ above, namely $u_{111}$. 

\medskip

One sees on the above synopsis that we shall in particular obtain solutions of the (gPAM) and (gKPZ) equations under the form
$$
u = {\sf P}_{f(u)}Z_1 + (2\alpha).
$$
Corollary 1.11 of Hairer and Pardoux' work \cite{HairerPardoux} follows then from Proposition \ref{prop:derive}, proved in Appendix \ref{SectionAppendixContinuity}. We denote here by $[0,T]$ a time interval on which the solution of the equation is defined.

\medskip

\begin{cor}
\label{CorEstimate}   {\it 
Let $f$ be $C_b^5$. For $0<t<T$, there exists a positive constant $C$ such that one has the estimate
$$
\big| u(e') - u(e) - f\big(u(e)\big)\big(Z_1(e') - Z_1(e)\big)\big| \leq C \left(\sqrt{|\tau-\sigma|}+|y-x|\right)^{2\alpha},
$$
uniformly in $e' = (\sigma,y)$ and $e=(\tau,x)$ with $|\tau-\sigma| \leq \frac{T}{2}$.   }
\end{cor}

\bigskip

\subsection{The mechanics of computations within paracontrolled calculus}
\label{SubsectionMechanicsComputations}

The basics of the paracontrolled analysis of singular PDEs are easily grasped by a parallel with It\^o calculus. Denote by $a,b,c$ three generic continuous martingales. The following computational rules appear as fundamental in stochastic calculus.   \vspace{0.1cm}

\begin{itemize}
   \item[\textcolor{gray}{$\bullet$}] \textit{The basic It\^o formula}
   $$
   d(ab) = a\,db + b\,da + d\langle a, b\rangle.   
   $$
   
   \item[\textcolor{gray}{$\bullet$}]\textit{ It\^o formula}
   $$
   d\big(f(a)\big) = f'(a)\, da + \frac{1}{2}f''(a) \,d\langle a, a\rangle.   
   $$

   \item[\textcolor{gray}{$\bullet$}] \textit{Bracket rule for stochastic integrals}
  	$$
  	d\left\langle \int adb , c\right\rangle - a \,d\langle b,c\rangle = 0.  	
  	$$ 
\end{itemize}

\ssk

The building blocks of the first order paracontrolled calculus devised by Gubinelli, Imkeller and Perkowski in \cite{GIP} are the exact counterparts of the above three points, with the paraproduct operator in the role of the (derivative of the) stochastic integral and the diagonal operator in the role of the (derivative of the) bracket. For $a,b,c$ functions or distributions with some precise regularity properties, we have the following facts.   \vspace{0.1cm}

\begin{itemize}
   \item[\textcolor{gray}{$\bullet$}] \textit{Paraproduct decomposition}
   $$
   ab = {\sf P}_ab + {\sf P}_ba + {\sf \Pi}(a,b)   
   $$
   where ${\sf P}$ is the paraproduct and $\Pi$ the resonant part
  \item[\textcolor{gray}{$\bullet$}] \textit{Bony's paralinearisation}
   $$
   f(a) = {\sf P}_{f'(a)}a + (\textrm{remainder})
   $$   

   \item[\textcolor{gray}{$\bullet$}] \textit{Fundamental corrector estimate}. The operator
   \begin{equation}
   \label{EqDefnCorrector}
   {\sf C}(a,b,c) := {\sf \Pi}\big({\sf P}_ab,c\big) - a{\sf \Pi}(b,c)   
   \end{equation}
	is continuous for certain ranges of regularity exponents for its arguments.
\end{itemize}

\ssk

The three step resolution process of Section \ref{SubsectionOverview} for the study of a singular PDE requires that we refine these tools. The task of writing the right hand side $f(u,\partial u, \zeta)$ of a singular PDE under the form \eqref{EqRHS} can be divided into two sub-tasks. First, making sense of the ill-defined products that appear in $f(u,\partial u, \zeta)$, up to the datum of a  number of formal multilinear functions of the noise $\zeta$ -- the enhancement $\widehat{\zeta}$. Second, writing $f(u,\partial u, \zeta)$ under the form \eqref{EqRHS}, needed to run the fixed point argument -- step 3 in Section \ref{SubsectionOverview}. Ill-defined quantities appear under the form of (possibly multi-)linear maps ${\sf E}(\cdot)$, defined on H\"older spaces, and formally taking values in a H\"older space $\mcC^\gamma$ of negative regularity exponent $\gamma<0$.   \vspace{0.15cm}

\begin{itemize}
   \item In the line of identity \eqref{EqDefnCorrector}, to deal with the first task, we systematically use the paracontrolled structure to write the following kind of decomposition. Let ${\sf E}$ stand for a map that formally sends $\mcC^\beta$ into $\mcC^\gamma$. For $a\in\mcC^\alpha$ and $b\in\mcC^\beta$, with $\alpha,\beta>0$, one has
\begin{equation}
\label{EqFundExpansion1}
\begin{split}
{\sf E}({\sf P}_ab) &= a {\sf E}(b) + {\sf E}^+(a,b)   \\
				&= {\sf P}_a{\sf E}(b) + {\sf P}_{{\sf E}(b)}a + {\sf \Pi}\big(a,{\sf E}(b)\big) + {\sf E}^+(a,b),
\end{split}
\end{equation}
for some linear operator ${\sf E}^+(\cdot,b)$ that formally takes values in the H\"older space $\mcC^{\gamma+\alpha}$. The map ${\sf \Pi}\big(\cdot,{\sf E}(b)\big)$ also formally takes values in $\mcC^{\gamma+\alpha}$, while the two paraproduct terms make perfect sense if ${\sf E}(b)$ does. The regularity exponents of the two possibly ill-defined terms ${\sf \Pi}\big(a,{\sf E}(b)\big)$ and ${\sf E}^+(a,b)$ have been increased by $\alpha$, compared to the regularity exponent of ${\sf E}({\sf P}_ab)$. One can iterate the kind of expansion given by \eqref{EqFundExpansion1} on ${\sf \Pi}\big(a,{\sf E}(b)\big)$ and ${\sf E}^+(a,b)$ if $a$ has a paracontrolled structure. Iterating the expansion \eqref{EqFundExpansion1} as many times as necessary, if possible, eventually leads to a decomposition of the initial quantity into a sum of well-defined terms, up to the a priori datum of a number of terms, like the above ${\sf E}(b)$.   \vspace{0.1cm}
   
   \item The second task involves giving a special form to the evaluation on some paraproduct terms of well-defined continuous linear functions ${\sf F}$. We systematically use for that purpose the following kind of decomposition. If ${\sf F}$ takes values in a H\"older space $\mcC^\gamma$, \emph{whatever} $\gamma\in\bbR$, and for $a\in\mcC^\alpha$ and $b\in\mcC^\beta$, one has
\begin{equation}
\label{EqFundExpansion2}
{\sf F}({\sf P}_ab) = {\sf P}_a{\sf F}(b) + {\sf F}^+(a,b),
\end{equation}
for some continuous bilinear operator ${\sf F}^+$ that takes values in $\mcC^{\gamma+\alpha}$.
\end{itemize}   \vspace{0.15cm}

As an example of computation, let $a_1, a_2, b$ be $\alpha$-H\"older functions, with $0<\alpha<1$, and let ${\sf E}$ send the space of $\alpha$-H\"older functions into the space of $\gamma$-H\"older functions, for some negative regularity exponent $\gamma$ with $-2\alpha<\gamma<-\alpha$. One has from formula \eqref{EqFundExpansion1} 
$$
{\sf E}({\sf P}_ab) = {\sf P}_a{\sf E}(b) + {\sf P}_{{\sf E}(b)}a + {\sf \Pi}\big(a,{\sf E}(b)\big) + {\sf E}^+(a,b),
$$
where the first term is $\gamma$-H\"older, the second term is $(\gamma+\alpha)$-H\"older, an elementary property of the paraproduct operator, and the third and fourth terms have formal regularity $\gamma+\alpha$. If $a={\sf P}_{a_1}a_2$, then we have the ${\sf F}$-type decomposition for the $(\gamma+\alpha)$-H\"older valued function $a\mapsto{\sf P}_{{\sf E}(b)}a$,
$$
{\sf P}_{{\sf E}(b)}a = {\sf P}_{a_1}\big({\sf P}_{{\sf E}(b)}a_2\big) + {\sf \Pi}^+_{{\sf E}(b)}(a_1,a_2),
$$
with the first term of regularity $(\gamma+\alpha)$ and the second term of regularity $(\gamma+2\alpha)$. We also have the two ${\sf E}$-type decompositions
\begin{equation*}
\begin{split}
{\sf \Pi}\big(a,{\sf E}(b)\big) &= a_1{\sf \Pi}\big(b_1,{\sf E}(b)\big) + {\sf \Pi}\big(\cdot,{\sf E}(b)\big)^+(a_1,b_1)   \\
								   &= {\sf P}_{a_1}{\sf \Pi}\big(b_1,{\sf E}(b)\big) + {\sf P}_{{\sf \Pi}\big(b_1,{\sf E}(b)\big)}a_1 + {\sf \Pi}\big(a_1,{\sf \Pi}\big(b_1,{\sf E}(b)\big)\Big) + {\sf \Pi}\big(\cdot,{\sf E}(b)\big)^+(a_1,b_1),
\end{split}
\end{equation*}
and 
\begin{equation*}
\begin{split}
{\sf E}^+(a,b) &= a_1{\sf E}^+(a_2,b) + \{{\sf E}^+(\cdot,b)\}^+(a_1,a_2)   \\
			   &= {\sf P}_{a_1}{\sf E}^+(a_2,b) + {\sf P}_{{\sf E}^+(a_2,b)}a_1 + {\sf \Pi}\big(a_1,{\sf E}^+(a_2,b)\big) + \{{\sf E}^+(\cdot,b)\}^+(a_1,a_2),
\end{split}
\end{equation*}
In both expressions, the first term has regularity $(\gamma+\alpha)$, the second term has regularity $(\gamma+2\alpha)$, while the third and fourth terms have formal regularity $(\gamma+2\alpha)$. If one is only interested in having a description of ${\sf E}({\sf P}_ab)$ up to terms $(\sharp)$ of true or formal positive regularity $(\gamma+2\alpha)$, one then has, with $a={\sf P}_{a_1}a_2$,
$$
{\sf E}({\sf P}_{{\sf P}_{a_1}a_2}b) = {\sf P}_{{\sf P}_{a_1}a_2}{\sf E}(b) + \Big\{{\sf P}_{a_1}\big({\sf P}_{{\sf E}(b)}a_2\big) + {\sf P}_{a_1}{\sf \Pi}\big(a_2,{\sf E}(b)\big) + {\sf P}_{a_1}{\sf E}^+(a_2,b)\Big\} + (\sharp).
$$

None of the computations of Section \ref{SectionNonlinearPDEs}, dealing with the paracontrolled analysis of concrete examples of singular PDEs, is more complicated than what we have just done. Convenient notations will be used to work with iterated operators; they are introduced in Section \ref{SectionToolKit}. (Note here that such expansions can be done on systems of singular PDEs as well.) On a technical level, three ingredients are used to run the three step scheme of Section \ref{SubsectionOverview}.   \vspace{0.15cm}

\begin{enumerate}
   \item[\textcolor{gray}{$\bullet$}] The pair $\big({\sf P},\widetilde{\sf P}\,\big)$ of \textbf{\textsf{intertwined paraproducts}} introduced in \cite{BBF15}. It is used to define a continuous map $\Phi$ from the solution space $\mcS$ to itself.   \vspace{0.15cm}
   
   \item[\textcolor{gray}{$\bullet$}] A \textbf{\textsf{high order paracontrolled expansion formula}} generalizing Bony's paralinearization formula is used to give a paracontrolled expansion of a non-linear function of any $\alpha$-H\"older function $u$, for $0<\alpha<1$. See section \ref{SectionTaylor}.   \vspace{0.15cm}
   
   \item[\textcolor{gray}{$\bullet$}] \textbf{\textsf{ Continuity results.}} We introduce in Section \ref{SectionToolKit} a number of operators and prove their continuity. These are the operators that appear in applying the decompositions \eqref{EqFundExpansion1} and \eqref{EqFundExpansion2} in the analysis of a generic right hand side $f\big(u,\partial u, \zeta\big)$ for equation \eqref{EqCanonicalSPDE}.
\end{enumerate}

\bigskip

\subsection{Setting and results}
\label{SubsectionSettingResults}

We adopt in this work essentially the same geometric and functional setting as in our previous work \cite{BBF15}, slightly restricted so as not to bother here the reader with the use of weighted functional spaces. All this work could be formulated in the more general geometric/functional setting of \cite{BBF15}; we refrain from doing this as it may blur the simple ideas that we want to promote in this work. Let then $(M,d,\mu)$ stand for a compact smooth Riemannian manifold equipped with a doubling measure $\mu$ that may differ from the canonical Riemmanian volume measure. Let $V_1,\dots, V_{\ell_0}$ stand for some smooth vector fields on $M$, identified with first order differential operators. Given a tuple $I = (i_1,\dots,i_k)$ in $\{1,\dots,\ell_0\}^k$, we shall set $|I| := k$ and 
$$
V_I := V_{i_k}\cdots V_{i_1}.
$$
Set 
$$
L := -\sum_{\ell=1}^{\ell_0} V_\ell^2
$$
and assume that $L$ \textit{is elliptic}, so that the $V_i$ span smoothly at every point of $M$ the whole tangent space. So there exist smooth functions $(\gamma_i)_{1\leq i\leq \ell_0}$ such that for every function $f\in C^1(M,{\mathbb R})$ and $x\in M$ we have
$$ 
\nabla f(x) = \sum_{\ell=1}^{\ell_0} \gamma_\ell(x) V_\ell(f)(x) V_\ell(x).
$$
The operator $L$ is then a sectorial operator in $L^2(M)$, it is injective on the quotient space of $L^2(M)$ by the space of constant functions, it has a bounded $H^\infty$-calculus on $L^2(M)$, and $-L$ generates a holomorphic semigroup $(e^{-tL})_{t>0}$ on $L^2(M)$ -- see \cite{DM}. The above class of operators includes obviously the Laplacian on the flat torus. Note that under the above smoothness and ellipticity conditions, the semigroup $e^{-tL}$ has regularity estimates at any order, by which we mean that for every tuple $I$, the operators $\Big(t^\frac{| I |}{2}V_I\Big)e^{-tL}$ and $e^{-tL}\Big(t^\frac{| I |}{2}V_I\Big)$ have kernels $K_t(x,y)$ satisfying the Gaussian estimate
\begin{equation*}
\Big| K_t(x,y) \Big| \lesssim  \frac{1}{\mu\big(B(x,\sqrt{t})\big)}\, e^{-c\,\frac{d(x,y)^2}{t}}
\end{equation*}
and the following regularity estimate, for some constants which may depend on $|I|$. For $d(x,z)\leq \sqrt{t}$
\begin{equation*}
\Big| K_t(x,y) - K_t(x',y)\Big| \lesssim \frac{d(x,x')}{\sqrt{t}}  \frac{1}{\mu\big(B(x,\sqrt{t})\big)} \, e^{-c\,\frac{d(x,y)^2}{t}}.
\end{equation*} 
Note again that we could equally well develop paracontrolled calculus in the more general setting adopted in our previous work \cite{BBF15}; we refrain from doing that here as it could obscure the simplicity of the ideas put forward here.

\medskip

Given a finite time horizon $T$, we define the parabolic space $\mcM$ as 
$$
\mcM := [0,T]\times M,
$$
and equip it with the parabolic metric 
$$ 
\rho\big((\tau,x), (\sigma,y)\big) = \sqrt{|\tau-\sigma|} + d(x,y)
$$
and the parabolic measure $\nu= dt \otimes \mu$. 
Then $(\mcM,\rho,\nu)$ is a doubling space of homogeneous type. Note that for  $(\tau,x) \in \mcM$ and small positive radius $r$, the parabolic ball $B_\mcM\big((\tau,x),r\big)$ has volume 
$$ 
\nu \Big( B_\mcM \big( (\tau,x),r \big) \Big) \approx r^2\,\mu\big( B_M(x,r) \big).  \vspace{0.15cm}
$$
We shall denote by $e=(\tau,x)$ a generic element of the parabolic space $\mcM$.

\medskip 

We have chosen to work in the scale of H\"older spaces; this makes life easier, although we could equally develop paracontrolled calculus in the functional setting of Sobolev spaces, in the line of what we did in our previous work \cite{BB15}. For a real number $s$, we will denote by $C^s=C^s(M)$, the H\"older space on $M$ of order $s$, defined in terms of Besov spaces; and $\mcC^s=\mcC^s(\mcM)$, the parabolic H\"older space. We refer the reader to Appendix \ref{SectionAppendix} for more details on these spaces. Following our previous work \cite{BBF15}, one can define from $L$ only parabolic paraproduct and resonant operators that have good continuity properties in the scale of parabolic H\"older spaces -- see Appendix \ref{SubsectionParaproducts}. The high order paracontrolled expansion formula and the continuity results stated in Section \ref{SectionTaylor} and Section \ref{SectionToolKit} respectively, and fully proved in Appendix \ref{SectionAppendixTaylor} and \ref{SectionAppendixContinuity}, make use of these operators and provide the spine of paracontrolled calculus. They are the main contributions of this work. 

\bigskip

We illustrate our approach of the study of singular PDEs by proving well-posedness results for the $3$-dimensional generalised parabolic Anderson model equation (gPAM)
\begin{equation}
\label{EqgPAM}
\mathscr{L} u := (\partial_t + L) u = f(u)\zeta,
\end{equation}
and the generalised KPZ equation
\begin{equation}
\label{EqgKPZ}
\mathscr{L} u = f(u)\zeta + g(u)(\partial u)^2,
\end{equation}
on the one-dimensional torus, with a space/time white noise. The generalised parabolic Anderson model equation is a natural nonlinear generalisation of its linear counterpart, for which $f(u)=u$. The latter equation describes the evolution of a Brownian particle in a white noise environment. The generalised (KPZ) equation appears in the study of the random motion of a string on a manifold \cite{HairerRubber}; its study in the setting of regularity structures is the object of Hairer and co-authors' works \cite{HairerRubber, BrunedGabrielHairerZambotti}. The renormalisation of the 70ish terms that appear in the models for this equation motivated the development of systematic renormalisation procedures, such as done in the recent works of Bruned-Hairer-Zambotti \cite{BrunedHairerZambotti} and Chandra-Hairer \cite{ChandraHairer}. In the present work, we assume that a proper enhancement $\widehat{\zeta}$ of the noise $\zeta$ is given. Defining $\widehat{\zeta}$ in a stochastic setting is a very different question that is not addressed here.

\medskip

We have organised this work as follows. Section \ref{SectionTaylor} is dedicated to the statement and proof of a high order paracontrolled expansion formula generalising Bony's paralinearisation formula. The core results of the mechanics of computations within paracontrolled calculus are obtained in Section \ref{SectionToolKit}, under the form of continuity results for a number of operators. We test paracontrolled calculus on the examples of the $3$-dimensional generalised parabolic Anderson model equation \eqref{EqgPAM} (Theorems \ref{thm:PAM1} and \ref{thm:PAM2}), and the generalised KPZ equation \eqref{EqgKPZ} (Theorem \ref{thm:KPZ}), in Section \ref{SectionNonlinearPDEs}. Appendix \ref{SectionAppendix} contains all the relevant details about the parabolic setting, approximation operators, H\"older spaces and paraproducts. Appendices \ref{SectionAppendixTaylor} and \ref{SectionAppendixContinuity} contain the proofs of a number of statements.

\bigskip

\noindent {\textbf{\textsf{Aknowledgements.}} \textit{The authors warmly thank the reviewers for their thorough readings of previous versions of the present work.}}

\bigskip

\section[\hspace{0.7cm} High order paracontrolled expansion]{High order paracontrolled expansion}
\label{SectionTaylor}

We explain in this section a simple procedure for getting an arbitrary high order expansion of a nonlinear map of a given H\"older function $u$ defined on the parabolic space $\mcM$, in terms of its parabolic regularity properties. It provides, in the setting of H\"older spaces, a refinement of Bony's paralinearisation theorem in the form of a viable alternative to the paradifferential calculus of Chemin \cite{CheminDuke}; see also \cite{CheminBull}, Theorem 2.5, p.18, for a more readable account of \cite{CheminDuke} in the case of a second order expansion. 

\medskip

In its simplest form, the classical paraproduct operator $(f,g)\mapsto\Pi^0_fg$ on the $d$-dimensional torus is defined via Fourier analysis by modulation of the high frequencies of a given 'reference' function/distribution $g$ by the low frequencies of another function/distribution $f$. For a function $f$ on the torus, we denote by $f=\sum \Delta_i(f)$ its usual Littlewood-Paley representation, where $\Delta_i(f)$ is the dyadic bloc with Fourier coefficients essentially only at the frequency scale $2^i$. Consider the Littlewood-Paley decompositions of two functions 
$$
f = \sum_i \Delta_i(f),\quad g = \sum_j \Delta_j(g),
$$
as sums of smooth functions with localized frequencies; the paraproduct of $g$ by $f$ is defined as
\begin{equation}
\label{EqBonyParaProduct}
\Pi^0_f\,g = \sum_{i<j-1} \Delta_i(f) \Delta_j(g),
\end{equation}
and the resonant part as
$$
\Pi^0(f,g) = \sum_{|i-j|\leq 1} \Delta_i(f) \Delta_j(g),
$$
so we have the product decomposition
$$ 
fg = \Pi^0_g f + \Pi^0_f g + \Pi^0(f,g).
$$
In the parabolic setting of Section \ref{SubsectionSettingResults}, one can define some paraproduct and resonant operators associated with the operator $L$ and its semigroup, that have the same regularity properties in the scale of parabolic H\"older spaces as the operator $\Pi^0$ in the scale of spatial H\"older spaces. We denote by ${\sf P}$ this paraproduct, introduced in \cite{BBF15}, and whose definition is recalled in Appendix \ref{SubsectionParaproducts}. It depends implicitly on an integer-valued parameter $b$ that is chosen once and for all, and whose precise value is irrelevant for our purposes. It is not crucial at that stage to go into the details of the definition of ${\sf P}$.

\ssk

The mechanics of the proof of our general Taylor expansion formula is fairly simple and better understood in the light of the proof of Bony's paralinearisation theorem given by Gubinelli, Imkeller and Perkowski in \cite{GIP}, which we recall first. 

\medskip

\begin{thm}[{\bf Bony's Paralinearisation}]  {\it
Let $f : \bfR\mapsto\bfR$ be a $C^2_b$ function and $u$ be a real-valued $\alpha$-H\"older function on the $d$-dimensional torus, with $0<\alpha<1$. Then
$$
f(u) = \Pi^0_{f'(u)} u + f(u)^\sharp
$$
for some remainder $f(u)^\sharp$ of spatial H\"older regularity $2\alpha$.   }
\end{thm}

\medskip

\begin{Dem}
This is a copy and paste from \cite{GIP}. Denote by $K_i$ the kernels of the Fourier projectors $\Delta_i$ corresponding to the Littlewood-Paley decomposition operator, and write $K_{\leq k}$ for $\sum_{i\leq k} K_i$, with associated operator $S_{k}:=\sum_{i\leq k} \Delta_i$. Note that by their definition we have, for any $i\geq 1$,
\begin{equation}
\label{EqNullMean}
\int_{\bfR^d} K_i(y)\,dy = 0;
\end{equation}
or more properly $\int_{\bfR^d} K_i(x,y)\,dy = 0$, for any $x\in\bfR^d$. The trick is then simply to write
$$
f(u) - \Pi^0_{f'(u)}(u) = \sum \Delta_i\big(f(u)\big) - S_{i-2}\big(f'(u)\big)\Delta_i(u) =: \sum \epsilon_i
$$
with 
$$
\epsilon_i(x) = \int K_i(x,y)K_{\leq i-2}(x,z)\Big\{f\big(u(y)\big) - f'\big(u(z)\big)u(y)\Big\}\,dzdy,
$$
and to take profit from the fact that $K_i$ has null mean for $i\geq 1$, as put forward in identity \eqref{EqNullMean}, to see that one also has, for $i\geq 1$,
$$
\epsilon_i(x) = \int K_i(x,y)K_{\leq i-2}(x,z)\Big\{f\big(u(y)\big) - f\big(u(z)\big)  - f'\big(u(z)\big)\big(u(y)-u(z)\big)\Big\}\,dzdy.
$$
One thus has
\begin{equation*}
\big|\epsilon_i(x)\big| \lesssim \|f''\|_{\infty} \int \big| K_i(x,y)K_{\leq i-2}(x,z)\big|\, \big| u(y)-u(z) \big|^2\,dzdy \lesssim 2^{-2i\alpha} \, \|u\|_{\mcC^\alpha}^2.
\end{equation*}
Since the functions $\epsilon_i$ are frequentially supported on the annulus of frequencies at the scale $2^{i}$, and the $\epsilon_i$ are uniformly bounded by $2^{-2i\alpha}$, the serie $\sum_i \epsilon_i$ defines a $2\alpha$-H\"older function.
\end{Dem}

\medskip

One can play exactly the same game and prove a general expansion result in a parabolic setting, with our paraproduct ${\sf P}$ in the role of the comparison operator. 

\medskip

\begin{thm}[\textsf{\textbf{Higher order Taylor expansion}}]   \label{ThmTaylorExpansion}   {\it
Let $f : \bfR\mapsto\bfR$ be a $C^4_b$, and let $u$ be a real-valued $\alpha$-H\"older function on the parabolic space $\mcM$, with $0<\alpha<1$. Then
\begin{equation}
\label{EqTaylorExpansion}
\begin{split}
f(u) &= {\sf P}_{f'(u)}(u) + \frac{1}{2}\,\left\{{\sf P}_{f^{(2)}(u)}(u^2) - 2{\sf P}_{f^{(2)}(u)u}(u) \right\}   \\
&+ \frac{1}{3!}\Big\{ {\sf P}_{f^{(3)}(u)}(u^3) - 3{\sf P}_{f^{(3)}(u)u}(u^2) + 3{\sf P}_{f^{(3)}(u)u^2}(u) \Big\} + f(u)^\sharp
\end{split}
\end{equation}
for some remainder $f(u)^\sharp$ of parabolic H\"older regularity $4\alpha$. If $f\in C^5_b$ then the remainder term $f(u)^\sharp$ is a locally Lipschitz function of $u$, in the sense that 
\begin{equation}
\label{EqLipBehaviourRemainderTaylor}
\big\|f(u)^\sharp - f(v)^\sharp\big\|_{\mcC^{4\alpha}} \lesssim \big(1+\|u\|_{\mcC^\alpha}+\|v\|_{\mcC^\alpha}\big)^{4} \|u-v\|_{\mcC^\alpha}.
\end{equation}
 }
\end{thm}

\medskip

We give here a proof of this statement in the case where $u$ is a time-independent function on the $d$-dimension torus and we can use the elementary paraproduct $\Pi^0$ instead of ${\sf P}$. The full proof of Theorem \ref{ThmTaylorExpansion} is given in Appendix \ref{SectionAppendixTaylor},  Theorem \ref{ThmTaylorExpansion-bis}; we hope this way of proceeding will make the reasoning clear and technical-free. 

\medskip

\begin{Dem}
Let us prove the second order formula in the special case where $u : {\bfT}^d\rightarrow\bfR$, and we use the elementary paraproduct $\Pi^0$ in place of ${\sf P}$. The claim amounts in that case to proving that 
$$
(\star) := f(u) - \Pi^0_{f'(u)}(u) - \frac{1}{2}\,\left\{\Pi^0_{f^{(2)}(u)}(u^2) - 2\,\Pi^0_{f^{(2)}(u)u}(u) \right\} 
$$
is a $3\alpha$-H\"older function on the torus. As in the proof of Bony's paralinearisation result, write $(\star)$ under the form
{\small $$
\sum \Delta_i\big(f(u)\big) - S_{i-2}\big(f'(u)\big)\Delta_i(u) - \left\{\frac{1}{2}S_{i-2}\big(f^{(2)}(u)\big)\Delta_i(u^2) - S_{i-2}\big(f^{(2)}(u)u\big)\Delta_i(u)\right\}  =:  \sum \epsilon_i.
$$   }
For each $i\geq 1$, we have 
\begin{equation*}
\begin{split}
\epsilon_i(x) = \int &K_i(x,y)K_{\leq i-2}(x,z)   \\
							&\left\{\int_0^1 f^{(2)}\big(u(z)+t(u(y)-u(z))\big) \big(u(y)-u(z)\big)^2\,tdt\right.   \\
							&\quad\quad\quad\quad\quad\left. -\frac{1}{2}f^{(2)}\big(u(z)\big)\,u^2(y) + f^{(2)}\big(u(z)\big) \,u(z) u(y) \right\} \, dzdy, 
\end{split}
\end{equation*}
which we can rewrite as
\begin{equation*}
\begin{split}
\epsilon_i(x) = \int &K_i(x,y)K_{\leq i-2}(x,z)   \\
							&\int_0^1\int_0^1 f^{(3)}\big(u(z)+st(u(y)-u(z))\big)\, \big(u(y)-u(z)\big)^3\,ds\,t^2dt \, dzdy, 
\end{split}
\end{equation*}
using once again the fact that the kernels $K_i(x,\cdot)$ have null mean. One reads on this expression for $\epsilon_i$ that it is of order $2^{-3i\alpha}$, uniformly in $x$. See Appendix \ref{SectionAppendixTaylor} for a full proof of the statement in the parabolic setting.

\medskip

To prove the stability estimate \eqref{EqLipBehaviourRemainderTaylor}, write 
$$
f(u)^\sharp =: \sum_i \epsilon_{i,u}, \quad f(v)^\sharp =: \sum_i \epsilon_{i,v}
$$ 
with $\epsilon_{i,u}$ equal to the above $\epsilon_i$ and $\epsilon_{i,v}$ defined similarly with $v$ in place of $u$. A uniform estimate of $\epsilon_{i,u} - \epsilon_{i,v}$ provides an estimate on $f(u)^\sharp - f(v)^\sharp$. The expression for $\epsilon_{i,u} - \epsilon_{i,v}$ involves the difference
$$ 
f^{(3)}\big(u(z)+st(u(y)-u(z))\big)\, \big(u(y)-u(z)\big)^3 - f^{(3)}\big(v(z)+st(v(y)-v(z))\big)\, \big(v(y)-v(z)\big)^3.
$$
Use the fact that $f^{(3)}$ is Lipschitz continuous, that $u$ and $v$ are $\mcC^\alpha$, to see that each of the two terms above is controlled by $d(y,z)^{3 \alpha}$ and 
$$ 
\Big| f^{(3)}\Big(u(z)+st\big(u(y)-u(z)\big)\Big) - f^{(3)}\Big(v(z)+st\big(v(y)-v(z)\big)\Big) \Big| \lesssim \|f^{(4)}\|_\infty  \|u-v\|_{C^\alpha}
$$
and
$$ 
\big| \big(u(y)-u(z)\big)^3 - \big(v(y)-v(z)\big)^3\big| \lesssim \|u-v\|_{C^\alpha} \, d(y,z)^{3\alpha} \big(\|u\|_{C^\alpha} + \|v\|_{C^\alpha} \big)^2.
$$
Estimate \eqref{EqLipBehaviourRemainderTaylor} is obtained by combining these two estimates together.
\end{Dem}

\medskip

Observe that the expansion \eqref{EqTaylorExpansion} is exact, $f(u)^\sharp=0$, for a polynomial function $f$ of degree at most $3$. The above expansion formula for $f(u)$ is conveniently rewritten under the form
\begin{equation*}
\begin{split}
f(u) &= {\sf P}_{f'(u) - uf^{(2)}(u) + \frac{1}{2}\,u^2f^{(3)}(u)}(u) + \frac{1}{2}\,{\sf P}_{f^{(2)}(u) - uf^{(3)}(u)}(u^2) + \frac{1}{6}\,{\sf P}_{f^{(3)}(u)}(u^3)
+ f(u)^\sharp.
\end{split}
\end{equation*}

\medskip

Note here that the general paracontrolled expansion formula writes
$$
f(u) = \sum_{n=1}^k \frac{1}{n!} \sum_{j=0}^{n-1} (-1)^j {{n}\choose{j}} {\sf P}_{u^j f^{(n)}(u)}(u^{n-j}) + f(u)^\sharp,
$$
for a function $f$ of class $C^{k+1}$ with bounded $(k+1)^\textrm{th}$ derivative, and a remainder $f(u)^\sharp$ of parabolic H\"older regularity $(k+1)\alpha$. Note that each paraproduct $ {\sf P}_{u^j f^{(n)}(u)}(u^{n-j})$ is only of regularity $\alpha$, while the two brackets that appear in \eqref{EqTaylorExpansion} are respectively $2\alpha$ and $3\alpha$-H\"older. We shall see in Proposition \ref{PropTaylorPractical} how to write each bracket in \eqref{EqTaylorExpansion} as a sum of terms of regularity $2\alpha$ and $3\alpha$, respectively.

\bigskip

\section[\hspace{0.7cm} A toolkit for paracontrolled calculus]{A toolkit for paracontrolled calculus}
\label{SectionToolKit}

As said in Section \ref{SubsectionMechanicsComputations}, the mechanics of computations within paracontrolled calculus lies in the analysis of how some operators behave when estimated on paraproduct terms. As a rule of thumb, an ill-defined operator $\sf E$ satisfies an expansion of the form
$$
{\sf E}({\sf P}_ab) = a{\sf E}(b) + {\sf E}^+(a,b),
$$
and well-defined operators $\sf F$ satisfy an expansion of the form
$$
{\sf F}({\sf P}_ab) = {\sf P}_a{\sf F}(b) + {\sf F}^+(a,b).
$$
Both ${\sf E}^+$ and ${\sf F}^+$ are expected to take values in spaces of functions/distributions that are more regular than the typical elements of the spaces where ${\sf E}$ and ${\sf F}$ take values, respectively. Gubinelli, Imkeller and Perkowski's continuity result \eqref{EqDefnCorrector} is the archetype of such an expansion. To iterate this kind of expansion, we introduce in this section a number of operators and prove their continuity properties. Together with the Taylor formula of Section \ref{SectionTaylor}, the results of this section are our main contribution. 

The development of paracontrolled calculus beyond the first order calculus of \cite{GIP} requires the introduction of a modified paraproduct $\widetilde{\sf P}$, introduced in \cite{BBF15}, and the use of expansion formulas of the form
\begin{equation*}
\begin{split}
&{\sf E}(\widetilde{\sf P}_ab) = a{\sf E}(b) + {\sf E}^+(a,b),   \\
&{\sf F}(\widetilde{\sf P}_ab) = {\sf P}_a{\sf F}(b) + {\sf F}^+(a,b).
\end{split}
\end{equation*}
This technical point is needed to run the fixed point procedure described in step 3 of the three step process of the paracontrolled analysis of a given singular PDE; see Section \ref{SubsectionOverview}. Concretely, one works with $\widetilde{\sf P}$ as one works with ${\sf P}$. The modified paraproduct $\widetilde{\sf P}$ is given by the formula
$$
\widetilde{\sf P}_f g := \mathscr{L}^{-1}\Big({\sf P}_f(\mathscr{L} g)\Big),
$$
where $\mathscr{L}$ stands for the parabolic differential operator $(\partial_\tau + L)$ on the parabolic space $\mcM$. See Section 4.1 of \cite{BBF15} for a study of the continuity properties of $\widetilde{\sf P}$, and Appendix \ref{Subsectionapp2} for a digest. Recall from Section \ref{SubsectionMechanicsComputations} the parallel between the basic rules of stochastic calculus and the fundamentals of the first order paracontrolled calculus. This integral picture of paraproduct provides a useful guide for the intuition, where the time derivative $d$ plays the role of the operator $\mathscr{L}$. In that comparison, $\widetilde{\sf P}_fg$ corresponds to the formal quantity
$$
\int\left(\int f\,d^2g\right) \simeq \int f\,dg - \iint df\,dg,
$$
after an integration by parts. So the difference between $\sf P$ and $\widetilde{\sf P}$ is a kind of 'bracket' term, reminiscent of the It\^o-to-Stratonovich rule for stochastic integration.

\medskip

We provide in this section a number of continuity results for some operators involving the paraproduct and resonant operators, together with the modified paraproduct $\widetilde{\sf P}$.  \textsl{We state in this section our results in their general form, in the parabolic setting of Section \ref{SubsectionSettingResults}, and give proofs in the time-independent, space setting of the torus, of versions of each statement where we use $\Pi^0$ instead of $\widetilde{\sf P}$.} This should make it easier for the reader to go to the core of the machinery without fighting with some possibly overwhelming technicalities; full proofs are given in Appendix \ref{SectionAppendixContinuity}. We hope this way of proceding will convince the reader that the basic ideas involved here are elementary. \textsl{It is not necessary, for the purpose of solving a particular singular PDEs using the paracontrolled calculus method, to get into the details of the proofs of the different results given here. We invite the reader to have a look at the results only and then go directly to Section \ref{SectionNonlinearPDEs} to see them on stage.} 

\bigskip

\noindent \textbf{\textsf{A word of caution.}} We repeatedly use below the fact that ${\sf P}_{\bf 1}f = f$ for an arbitrary distribution. This is not true, strictly speaking, as one rather has ${\sf P}_{\bf 1}f = f + (\textrm{smooth})$, for an infinitly smooth additional term that is continuous and linear with respect to $f$, for $f$ in any H\"older space, with positive or negative H\"older exponent. Using the first identity rather than the second has no effect whatsoever on the analysis below, so we prefer not to burden the reader with these somewhat irrelevant additional terms and stick to the identity ${\sf P}_{\bf 1}f = f$.

\bigskip

\subsection[\hspace{-0.3cm} Corrector $\sf C$, commutator $\sf D$, and their iterates]{Corrector $\sf C$, commutator $\sf D$, and their iterates}
\label{SubSectionCorrectors}

We introduce in this section and four operators ${\sf C}, {\sf D}$ and ${\sf R}, {\sf S}$. The operators ${\sf C}$ and ${\sf D}$ are used to make sense of ill-defined products, while the operators ${\sf R}$ and ${\sf S}$ are used to write the right hand side $f(u,\partial u, \zeta)$ of Equation \eqref{EqCanonicalSPDE} in the form \eqref{EqRHS} needed  to run a fixed point argument. The definitions of these operators involves the paraproduct and resonant operators and the modified paraproduct $\widetilde{\sf P}$. We define similarly operators ${\sf C}^\circ, {\sf D}^\circ$ and ${\sf R}^\circ, {\sf S}^\circ$ using ${\sf P}$ instead of $\widetilde{\sf P}$.

\smallskip

To motivate the introduction of the different operators ${\sf C}, {\sf D}$ and their iterates, let set ourselves the task of making sense of the product $u\zeta$, of $u\in\mcC^\alpha$ and $\zeta\in\mcC^{\alpha-2}$, with $\frac{2}{5}<\alpha<\frac{1}{2}$. As a product of a $\beta$-H\"older function $(\beta)$  with $\zeta$ is well-defined if $\beta+\alpha-2>0$, and paraproducts are always well-defined in a H\"older setting, we concentrate here on the ill-defined terms that appear in the computations and write 
$$
(\beta)\zeta \sim 0
$$
to mean that the product $(\beta)\zeta$ is well-defined, and 
$$
u\zeta \sim {\sf \Pi}(u,\zeta)
$$
rather than 
$$
u\zeta = {\sf P}_u\zeta + {\sf P}_\zeta u + {\sf \Pi}(u,\zeta).
$$
This convention is only used in this paragraph. If the resonant term ${\sf \Pi}(u,\zeta)$ were defined it would have regularity $2\alpha-2$; we say that the term ${\sf \Pi}(u,\zeta)$ has formal regularity $2\alpha-2$. Assume 
$$
u = \sum_{i=1}^3 {\sf P}_{u_i}Z_i + (4\alpha),
$$
for some  functions $u_i\in\mcC^\alpha$, some functions $Z_i\in\mcC^{i\alpha}$ depending only on the noise $\zeta$, and a remainder $(4\alpha)\in\mcC^{4\alpha}$. We do not specify the structure of the remainder term $(4\alpha)$ as its product with $\zeta$ is well-defined. In the present work, all multilinear functions of the noise $\zeta$ only are assumed to be well-defined, even if they have negative formal regularity. Then 
$$
u\zeta \sim \sum_{i=1}^3 {\sf \Pi}({\sf P}_{u_i}Z_i,\zeta),
$$
and setting
$$
{\sf C}^\circ(a,b,c) := {\sf \Pi}({\sf P}_ab,c)-a{\sf \Pi}(b,c),
$$
for the above map $\Pi^+$, one has 
\begin{equation*}
\begin{split}
u\zeta &\sim \sum_{i=1}^3\Big(u_i\,{\sf \Pi}(Z_i,\zeta) + {\sf C}^\circ(u_i,Z_i,\zeta)\Big)   \\
		   &\sim \sum_{i=1}^3\Big({\sf \Pi}\big(u_i,{\sf \Pi}(Z_i,\zeta)\big) + {\sf C}^\circ(u_i,Z_i,\zeta)\Big).
\end{split}
\end{equation*}
Gubinelli, Imkeller and Perkowski proved in \cite{GIP} that ${\sf C}^\circ(a,b,c)$ is essentially well-defined whenever the sum of the H\"older exponents of $a,b$ and $c$ add up to a positive number, in which case it defines an element of regularity the sum of the regularity exponents. (There is a mild restriction on the range of the different regularity exponents.) The term ${\sf \Pi}(Z_i,\zeta)$ is of regularity $(i+1)\alpha-2$, so we have 
$$
u\zeta \sim \sum_{i=1}^2 \Big({\sf \Pi}\big(u_i,{\sf \Pi}(Z_i,\zeta)\big) + {\sf C}^\circ(u_i,Z_i,\zeta)\Big).
$$
The formal regularity $(i+2)\alpha-2$ of ${\sf C}^\circ(u_i,Z_i,\zeta)$ is negative. To proceed we assume that each $u_i$ is also given in paraproduct form
$$
u_i = \sum_{j^=1}^2 {\sf P}_{u_{ij}}Z_j + (3\alpha),
$$
for $u_{ij}\in\mcC^\alpha$ and a remainder $(3\alpha)$ of regularity $3\alpha$. Setting
$$
{\sf C}^\circ\big((a_1,a_2),b,c\big) := {\sf C}^\circ\big({\sf P}_{a_1}a_2,b,c\big) - a_1{\sf C}^\circ(a_2,b,c),
$$
for the operator $({\sf C}^\circ)^+$, we shall prove in Theorem \ref{ThmIteratedCorrector} that this map is well-defined whenever the sum of the H\"older exponents of $a_1,a_2,b$ and $c$ add up to a positive number. (Here again, there is a mild restriction on the range of the regularity exponents.) So we have
\begin{equation}   \label{EqCheatingRefinedCorrector}
\begin{split}
u\zeta &\sim\sum_{i=1}^2 \Big({\sf \Pi}\Big(u_{ij}, {\sf \Pi}\big(Z_j, {\sf \Pi}(Z_i,\zeta)\big)\Big) + {\sf C}^\circ\big((u_{ij},Z_j), Z_i, \zeta\big)\Big)   \\
		   &\sim {\sf \Pi}\Big(u_{11}, {\sf \Pi}\big(Z_1, {\sf \Pi}(Z_1, \zeta)\big)\Big) + {\sf C}^\circ\big((u_{11},Z_1),Z_1,\zeta\big).
\end{split}
\end{equation}
Assuming $u_{11}={\sf P}_{u_{111}}Z_1 + (2\alpha)$, for $u_{111}\in\mcC^\alpha$ and a $2\alpha$ remainder $(2\alpha)$, and iterating once more the $E$-expansion that we have just done twice shows that $u\zeta$ is be well-defined, under the assumption that all the multilinear functions of the noise $\zeta$ only that appear above are well-defined and have regularity their formal regularity. This kind of computation explains why we need the corrector operator ${\sf C}^\circ$ and its iterates. (The above mentioned restriction on the range of the regularity exponents in the continuity results for ${\sf C}^\circ$ imposes that the term ${\sf C}^\circ\big((3\alpha),Z_i,\zeta\big)$ in ${\sf C}^\circ(u_i,Z_i,\zeta)$ is treated differently from the others. A {\sl refined corrector} is introduced for that purpose in Section 3.1.1. Its use justifies identity \eqref{EqCheatingRefinedCorrector}.)

\smallskip

New things happen when we look at the product $u^2\zeta$, as 
$$
u^2\zeta \sim {\sf \Pi}\big(2{\sf P}_uu + {\sf \Pi}(u,u)\,,\,\zeta\big) \sim 2{\sf C}^\circ(u,u,\zeta) + {\sf \Pi}\big({\sf \Pi}(u,u),\zeta\big).
$$
The ${\sf C}^\circ$ term is dealt with as above. For the resonant term ${\sf \Pi}\big({\sf \Pi}(u,u),\zeta\big)$ we need first to use an ${\sf F}$-type expansion on ${\sf \Pi}(u,u)$ before using an ${\sf E}$-type expansion. We have
$$
{\sf \Pi}(u,u) = \sum_{i=1}^3 \Big({\sf P}_{u_i}{\sf \Pi}(Z_i,u) + {\sf D}^\circ(u_i,Z_i,u)\Big),
$$
with ${\sf D}^\circ(\cdot,\cdot,u) := \Pi^+(\cdot,\cdot,u)$, with the notation of Section \ref{SubsectionMechanicsComputations}. The resonant terms 
$$
{\sf \Pi}\big({\sf P}_{u_i}{\sf \Pi}(Z_i,u),\zeta\big)
$$ 
are dealt with as above, using the corrector operator ${\sf C}^\circ$. For the resonance of ${\sf D}^\circ(u_i,z_i,u)$ with $\zeta$, we use the paracontrolled structures of $u$ and $u_i$ and the fact that the operator ${\sf D}^\circ$ satisfies an ${\sf F}$-type expansion with respect to its first and third arguments -- this is proved in the next section. This gives
\begin{equation*}
\begin{split}
\sum_{i=1}^3{\sf D}^\circ(u_i,Z_i,u) &= {\sf D}^\circ(u_1,Z_1,u) + (4\alpha)  \\
												   &= {\sf P}_{u_{11}}\big({\sf P}_{u_1}{\sf D}^\circ(Z_1,Z_1,Z_1)\big) + (4\alpha)',
\end{split}
\end{equation*}
for remainders $(4\alpha), (4\alpha)'$ of regularity $4\alpha$, and 
$$
\sum_{i=1}^3{\sf \Pi}\big({\sf D}^\circ(u_i,Z_i,u),\zeta\big) \sim u_{11}{\sf \Pi}\big({\sf P}_{u_1}{\sf D}^\circ(Z_1,Z_1,Z_1),\zeta\big) \sim u_{11}u_1{\sf \Pi}\big({\sf D^\circ}(Z_1,Z_1,Z_1),\zeta\big) \sim 0,
$$
provided ${\sf \Pi}\big({\sf D}^\circ(Z_1,Z_1,Z_1),\zeta\big)$ is given a priori as an element of $\mcC^{4\alpha-2}$.

\bigskip

\textbf{\textsf{3{\boldmath $.$}1{\boldmath $.$}1  Corrector $\sf C$ and the outer centering operator $\mathscr{C}${\boldmath$.$}}} Define on the space $L^\infty$ of bounded measurable functions on the parabolic space $\mcM$ the \textsf{\textbf{corrector}} {\sf C} as the operator 
$$
{\sf C}(f,g,h) := {\sf \Pi}\big(\widetilde{\sf P}_f g, h\Big) - f\,{\sf \Pi}(g,h).
$$
The next theorem is the workhorse of the first order paracontrolled calculus, such as devised in \cite{GIP} by Gubinelli, Imkeller and Perkowski; we recall it here, together with its proof, as this is our starting point. Recall we denote by $C^\alpha$ the spacial H\"older spaces on the torus and by $\mcC^\alpha$ the parabolic H\"older spaces over the compact manifold $M$.

\medskip

\begin{thm}
\label{ThmCorrectors}   {\it 
Let $\alpha, \beta, \gamma$ be regularity exponents, with $\alpha\in (0,1), \beta\in (-3,3)$ and $\gamma\in (-\infty,3)$. Assume $\alpha+\beta<3$, and 
   $$
   0< \alpha+\beta+\gamma <1,   \quad\textrm{while }\quad \beta+\gamma < 0.
   $$
Then, the corrector ${\sf C}$ extends  continuously as a function from $\mcC^\alpha\times\mcC^\beta\times\mcC^\gamma$ to $\mcC^{\alpha+\beta+\gamma}$.   }
\end{thm}

\medskip

\begin{Dem}
As advertised above, we prove here this continuity result for a simplified version of the operator $\sf C$, and refer the reader to Proposition \ref{prop:C1} in Appendix \ref{SectionAppendixContinuity} for full proofs. Assume we are working in the time-independent setting of the $d$-dimensional torus, with the operator
$$
{\sf C}^0(f,g,h) := \Pi^0\Big(\Pi^0_f g, h\Big) - f\,\Pi^0(g,h).
$$
We prove the continuity of the corrector ${\sf C}^0$, as a function from  $C^\alpha\times C^\beta\times C^\gamma$ to $C^{\alpha+\beta+\gamma}$, under the above assumptions on $\alpha,\beta,\gamma$.

\medskip

$\bullet$ We first give a heuristic proof of the statement. The resonant operator is given by
\begin{equation}
\label{EqDiagonalTerm}
\Pi^0(a,b)\simeq \sum_{|i-j|\leq 1} \Delta_i(a)\,\Delta_j(b).
\end{equation}
Write 
$$
{\sf C}^0(f,g,h) = \sum_{|i-j|\leq 1}  \Delta_i\Big(\Pi^0_f g\Big)\Delta_jh - f\Delta_i(g)\,\Delta_j(h),
$$
and set
$$
\epsilon'_i := \Delta_i\Big(\Pi^0_f g\Big) - f\Delta_i(g),
$$
so we have
$$ 
{\sf C}^0(f,g,h) = \sum_{|i-j|\leq 1}\epsilon'_i \, \Delta_j(h).
$$
The fact that $\epsilon'_i$ has $L^\infty$-norm of order $2^{-i(\alpha+\beta)}$ can be guessed on the expression
\begin{equation*}
\begin{split}
\epsilon'_i(x) &= \int K_i(x,y) \Big\{\big(\Pi^0_fg\big)(y) - f(x)g(y)\Big\}\,dy   \\
			    	 &= \int K_i(x,y) \Big\{\Pi^0_{f-f(x){\bf 1}}g\Big\}(y)\,dy.
\end{split}
\end{equation*}
As $y$ is concentrated near $x$, at scale $2^{-i}$, and we are looking at the $i^\textrm{th}$ Littlewood-Paley block of ${\sf P}_{f-f(\cdot)}g$, we expect 
$$
\big|\varepsilon_i'(x)\big| \lesssim 2^{-i\beta}\,\Big\|\Pi^0_{f-f(x)}g\Big\|_{C^\beta} \lesssim 2^{-i\beta}\,\big\| f-f(x)\big\|_{L^\infty}\,\|g\|_{C^\beta},
$$
with a term $\big\|f-f(x)\big\|_{L^\infty}$ involving only the neighborhood of $x$ of size $2^{-i}$, that is with 
$$
\big\|f-f(x)\big\|_{L^\infty} \lesssim 2^{-i\alpha} \|f\|_{C^\alpha},
$$
since $f$ is $\alpha$-H\"older. Such an estimate would imply the continuity of the corrector ${\sf C}$ as a function from $C^\alpha\times C^\beta\times C^\gamma$ to $C^{\alpha+\beta+\gamma}$, if $\alpha+\beta+\gamma$ is positive, since $h$ is $\gamma$-H\"older. This heuristic argument, however, does not make it clear why we need $\beta+\gamma$ to be negative to get the result.

\ssk

$\bullet$ A mathematically correct version of the above sketch of proof is done by estimating the $L^\infty$-norm of the dyadic blocks of $\varepsilon_i'$. For $k\geq i+2$ then 
$$
\Delta_k\epsilon'_i = - \Delta_k\big(f\Delta_i(g)\big) \simeq - \Delta_k(f) \Delta_i(g);
$$
here and below the $\simeq$ means that the right hand side is equal to a finite sum of terms of the given form. This is here a direct consequence of frequency considerations with the frequency supports of the dyadic projections $\Delta_k$, for $k\geq i+2$. So we have
$$
\big\|\Delta_k\epsilon'_i\big\|_{L^\infty} \lesssim 2^{-k\alpha}\,2^{-i\beta}\,\|f\|_{C^\alpha} \|g\|_{C^\beta}.
$$
For $k\leq i-2$ then 
$$ 
\Delta_k\epsilon'_i = - \Delta_k(f\Delta_i(g)) \simeq - \Delta_k\big( \Delta_i(f)\,\Delta_i(g)\big)
$$
hence
$$ 
\big\|\Delta_k\epsilon'_i\big\|_{L^\infty} \lesssim 2^{-i(\alpha+\beta)}\,\|f\|_{C^\alpha} \|g\|_{C^\beta}.
$$
We adopt the classical notation $S_{k-2}f$ for the partial sum $\sum_{\ell\leq k-2} \Delta_\ell(f)$ of the Littlewood-Paley decomposition, so for $|i-k|\leq 2$ we have
$$ 
\Delta_k\epsilon'_i \simeq \Delta_k\Big( \Delta_i(g) S_{i-2}(f) - S_{k+2}(f) \Delta_i(g)\Big),
$$
hence
$$ 
\big\|\Delta_k\epsilon'_i\big\|_{L^\infty} \lesssim 2^{-i(\alpha+\beta)}\,\|f\|_{C^\alpha} \|g\|_{C^\beta}.
$$
As a consequence, we always have the following estimate
\begin{equation} 
\label{eq::}
\big\|\Delta_k\epsilon'_i\big\|_{L^\infty} \lesssim 2^{-i\beta}\,2^{-\max(i,k)\alpha}\,\|f\|_{C^\alpha} \|g\|_{C^\beta}. 
\end{equation}
We can then estimate ${\sf C}^0(f,g\,;h)$ in some H\"older space. For a non-negative integer $k$, we have
\begin{align*}
\Delta_k\Big({\sf C}^0\big(f,g,h\big)\Big) & = \sum_{|i-j|\leq 1} \Delta_k\Big(\epsilon'_i \,\Delta_j h\Big) \\
&\simeq \underset{{|i-j|\leq 1}}{\sum_{i< k-2}} \Delta_k(\epsilon'_i) \,\Delta_j(h) + \underset{{|i-j|\leq 1}}{\sum_{k< i-2}} \Delta_k\Big(\Delta_i(\epsilon'_i) \,\Delta_j(h)\Big)   \\
&\hspace{0.6cm}+ \underset{{|i-j|\leq 1}}{\sum_{|k-i|\leq 2}} \Delta_k\Big(S_i(\epsilon'_i) \,\Delta_j(h)\Big)
\end{align*}
which is then controlled, using estimate \eqref{eq::}, by
\begin{align*} 
\Big\| &\Delta_k\big({\sf C}^0\big(f,g,h\big)\big) \Big\|_{L^\infty}   \\
&\lesssim \left(\sum_{i< k-2} 2^{-i\beta}2^{-k\alpha}2^{-i\gamma} + \sum_{k< i-2} 2^{-i(\alpha+\beta+\gamma)} + \sum_{|k-i|\leq 2} 2^{-i(\alpha+\beta+\gamma)} \right) \|f\|_{C^\alpha} \|g\|_{C^\beta} \|h\|_{C^\gamma} \\
& \lesssim 2^{-k(\alpha+\beta+\gamma)} \|f\|_{C^\alpha} \|g\|_{C^\beta} \|h\|_{C^\gamma};
\end{align*}
we used the condition $\beta+\gamma<0$ to estimate the first sum, and the condition $\alpha+\beta+\gamma>0$, to estimate the second sum. The fact that the latter estimate holds uniformly in $k$ concludes the proof of the $(\alpha+\beta+\gamma)$-H\"older regularity of the corrector.
\end{Dem}

\medskip

We emphasize the importance of the above heuristic proof of continuity of the corrector {\sf C} by introducing a notation. 

\medskip

\begin{definition*}
Given an endomorphism $A$ of some function space, we denote by $\scrC f$, or $\scrC_x f$, the function 
$$
(\mathscr{C} f)(\cdot) := f(\cdot) - f(x), 
$$
recentered around its value at the 'running' variable $x$, so that 
$$
A(\scrC f)(x) = A\big(f-f(x)\big)(x).
$$
Strictly speaking, the operator $\scrC$ is an operator on the space of operators $A$. The choice of letter $\mathscr{C}$ for this operator is for 'centering', and we call $\mathscr{C}$ the \textbf{\textsf{outer centering operator}}.
\end{definition*}

\medskip

In those terms, we have 
\begin{equation} 
\label{eq:formulec}
{\sf C}\big(f,g,h\big) = {\sf \Pi}\big(\widetilde{\sf P}_{\mathscr{C}f}g, h\Big),
\end{equation}
and 
$$
{\sf \Pi}\big({\sf P}_{\mathscr{C}{\sf P}_{\mathscr{C} b} c} g \,,\, h\Big)(x) = {\sf \Pi}\big({\sf P}_{{\sf P}_{b-b(x)} c - ({\sf P}_{b-b(x)} c)(x)} g \,,\, h\Big)(x),
$$
for instance. The main property of this operator is the following. For a function $f\in C^\alpha(\bfT^d)$ with $\alpha$ positive, we have
\begin{align*}
S_k(\mathscr{C}f)(x) & = S_k\big(f-f(x)\big)(x) = S_k(f)(x) - f(x)   \\
& = \sum_{\ell \geq k+1} \Delta_{\ell}(f)(x). 
\end{align*}
Since $f$ has positive regularity, the $L^\infty$ size of the dyadic blocks $\Delta_\ell f$ are exponentially decreasing with $\ell$, so 
\begin{equation}
\label{eq:appox}
\|S_k(\mathscr{C}f) \|_\infty + \| \Delta_k f \|_\infty \lesssim 2^{-\alpha k}.
\end{equation}
A very similar property holds in the parabolic setting, which is used in the proofs of the continuity results of this section given in Appendix \ref{SectionAppendixContinuity}.

\bigskip

\textbf{\textsf{3{\boldmath $.$}1{\boldmath $.$}2  Dealing with remainders: refined corrector ${\sf C}_{(1)}${\boldmath$.$}}}  The outer centering operator allows to take profit from the H\"older continuity regularity property of a function. To take further profit from the $\gamma$-Lipschitz regularity property of a function, with $1<\gamma<2$, one introduces a \textbf{\textsf{refined corrector}}, defined as follows in the model setting of the flat torus; see Definition \ref{DefnRefinedCorrector} for the definition in the parabolic setting. This operator is used to take care of a number of remainder terms in the paracontrolled analysis of ill-defined terms. Set 
$$
{\sf C}^0_{(1)}\big(f,g,h\big) := \Pi^0\big({\sf P}_f g, h\big) - f\,\Pi^0(g,h) - f'\,\Pi^0_{(1)}(g,h),
$$
where
$$
\Pi^0_{(1)}(g,h)(x) := \sum_{\vert i-j\vert\leq 1} \Delta_i\big((\cdot-x)g\big)(x)\,(\Delta_j h)(x),
$$
for $x\in\bf T$. An elementary refinement of the proof of Theorem \ref{ThmCorrectors} gives the following statement. We refer the reader to Definition \ref{DefnRefinedCorrector} for the definition of the parabolic counterpart ${\sf C}_{(1)}$ of ${\sf C}^0_{(1)}$.

\medskip

\begin{thm}
\label{ThmCorrectors1}   {\it 
Let $\alpha, \beta, \gamma$ be regularity exponents, with $\alpha\in (1,2), \beta\in (-3,3)$ and $\gamma\in (-\infty,3)$. Assume $\alpha+\beta<3$, and 
   $$
   0< \alpha+\beta+\gamma <1,   \quad\textrm{while }\quad \beta+\gamma < 0.
   $$
Then, the corrector ${\sf C}_{(1)}$ extends  continuously as a function from $\mcC^\alpha\times\mcC^\beta\times\mcC^\gamma$ to $\mcC^{\alpha+\beta+\gamma}$.   }
\end{thm}

\medskip

The proof of the next statement is identical to the proof of the continuity result for ${\sf C}$ from Theorem \ref{ThmCorrectors}; it is left to the reader. Recall from Definition \ref{DefnRefinedCorrector} that ${\sf C}_{(1)}$ is defined using the operators ${\sf \Pi}_{(1)}^\ell$.

\medskip

\begin{thm}   \label{ThmPi1}
Let $\alpha, \beta, \gamma$ be regularity exponents, with $\alpha,\beta\in (0,1)$, and $\gamma\in (-\infty,3)$. Assume 
   $$
   \beta+\gamma+1 < 0 <\alpha+\beta+\gamma+1 < 1.
   $$
Then all the operators 
$$
(f,g,h) \mapsto {\sf \Pi}_{(1)}^\ell\big(\widetilde{\sf P}_fg,h\big) - f{\sf \Pi}_{(1)}^\ell(g,h)
$$
extend continuously as functions from $\mcC^\alpha\times\mcC^\beta\times\mcC^\gamma$ to $\mcC^{\alpha+\beta+\gamma+1}$, for $1\leq \ell\leq \ell_0$.
\end{thm}

\bigskip

\textbf{\textsf{3{\boldmath$.$}1{\boldmath$.$}3  Iterated correctors{\boldmath$.$}}} Given a tuple of bounded functions $(a,b,c,d)$, set
$$
\widetilde{\sf P}^{\circ 2}(a,b,c) := \widetilde{\sf P}_{\widetilde{\sf P}_a b} c
$$
and 
$$
\widetilde{\sf P}^{\circ 3}(a,b,c,d) := \widetilde{\sf P}_{\widetilde{\sf P}^{\circ 2}(a,b,c)}d,
$$
and give similar definitions of $\Pi^{\circ 2}(a,b,c)$ and $\Pi^{\circ 3}(a,b,c,d)$, using only $\sf P$ operators, and $(\Pi^0)^{\circ 2}(a,b,c)$ and $(\Pi^0)^{\circ 3}(a,b,c,d)$, using only $\Pi^0$ operators, respectively. Depending on whether or not such a paraproduct appears in the low frequency, in place of $f$, or high frequency, in place of $g$, in the formulas for the corrector ${\sf C}(f,g,h)$ or the commutator ${\sf D}(f,g,h)$, we shall talk about \textsf{\textbf{lower}} or \textsf{\textbf{upper iterated}} operators.

\medskip

$\bullet$ We define the $4$ and $5$-linear \textsf{\textbf{lower iterated correctors}} by the formulas
\begin{equation}
\label{eq:C2bis}
\begin{split}
{\sf C}\big((a_1,a_2),g,h\big) :=&\; {\sf C}\Big(\widetilde{\sf P}_{a_1}a_2,g,h\Big) - a_1\,  {\sf C}(a_2,g,h)   \\
												  =&\; {\sf \Pi}\big(\widetilde{\sf P}^{\circ 2}(a_1,a_2,g)\, ,\,h\Big) - \Big\{ \big(\widetilde{\sf P}_{a_1}a_2\big)\,{\sf \Pi}(g,h) + a_1\,{\sf \Pi}\big(\widetilde{\sf P}_{\mathscr{C} a_2} g, h\Big)\Big\}, 
\end{split}
\end{equation}
and 
\begin{equation}
\label{EqIteratedCorrector}
\begin{split}
{\sf C}\Big(\big((a_1,a_2),a_3\big),g,h\big) := {\sf C}\Big(\big(\widetilde{\sf P}_{a_1}a_2,a_3\big),g,h\Big) - a_1 \, {\sf C}\big((a_2,a_3),g,h\big),
\end{split}
\end{equation}
also equal to
{\small \begin{equation*}
{\sf \Pi}\big(\widetilde{\sf P}^{\circ 3}(a_1,a_2,a_3,g)\, ,\,h\Big) - \Big\{ \big(\widetilde{\sf P}^{\circ 2}(a_1,a_2,a_3)\big)\,{\sf \Pi}(g,h) + \big(\widetilde{\sf P}_{a_1}a_2\big)\,{\sf \Pi}\big(\widetilde{\sf P}_{\mathscr{C} a_3} g, h\Big) + a_1\,{\sf \Pi}\big(\widetilde{\sf P}_{\mathscr{C}\widetilde{\sf P}_{\mathscr{C} a_2} a_3} g\,,\, h\Big)\Big\}.
\end{equation*}   }
The conditions $(\cdots)<3$ that appear in Theorem \ref{ThmIteratedCorrector} below are purely technical; a choice of implicit constant $b$ in the definition of the paraproduct operator ${\sf P} = {\sf P}^{(b)}$ would change the bound $3$ for any other bound. In any concrete situation, one can assume that such a good choice of parameter $b$ has been done and forget about that condition.

\medskip

\begin{thm}
\label{ThmIteratedCorrector}   {\it 
Let $\alpha_1,\alpha_2,\mu$ be regularity exponents in $(0,1)$, and $\alpha_3\in (-3,3)$. Let $\nu\in (-\infty,3]$ be another regularity exponent. 
\begin{itemize}
   \item[\textcolor{gray}{$\bullet$}] Assume $(\alpha_1+\alpha_2+\mu) < 3$, and 
   \begin{equation*}
   \begin{split}
   &\alpha_2+\mu+\nu < 0, \quad \alpha_1+\mu+\nu<0  \\
   &(\alpha_1+\alpha_2+\mu+\nu) \in (0,1).
   \end{split}
   \end{equation*}
   Then the $4$-linear lower iterated corrector
   \begin{equation*}
   \begin{split}
   \mcC^{\alpha_1}\times\mcC^{\alpha_2}\times \mcC^{\mu}\times\mcC^{\nu} &\rightarrow \mcC^{\alpha_1+\alpha_2+\mu+\nu}   \\
   \big(a_1,a_2,g,h\big) &\mapsto {\sf C}\big((a_1,a_2),g,h\big)
   \end{split}
   \end{equation*} 
   is continuous.  \vspace{0.15cm}  
   
   \item[\textcolor{gray}{$\bullet$}]  Assume that $(\alpha_1+\alpha_2+\alpha_3+\mu)<3$, and 
   \begin{equation*}
   \begin{split}
   &(\alpha_3+\mu+\nu) < 0, \quad \alpha_2+\alpha_3+\mu+\nu<0  \\
   &(\alpha_1+\alpha_2+\beta+\alpha_3+\nu) \in (0,1).
   \end{split}
   \end{equation*}      
   The $5$-linear lower iterated corrector 
   \begin{equation*}
   \begin{split}   
   \mcC^{\alpha_1}\times\mcC^{\alpha_2}\times\mcC^{\alpha_3}\times \mcC^{\mu}\times\mcC^{\nu} &\rightarrow \mcC^{\alpha_1+\alpha_2+\alpha_3+\mu+\nu}   \\
   \big(a_1,a_2,a_3,g,h\big) &\mapsto {\sf C}\big(((a_1,a_2),a_3),g,h\big)
   \end{split}
   \end{equation*}    
   is continuous.     
\end{itemize}   }
\end{thm}

\medskip

\begin{Dem}
To get a clear idea of the mechanics at play, we only prove here the analogue statement of the time-independent setting of the flat torus. That means that we aim to prove that the formula
\begin{equation*}
\begin{split} 
&\Pi^0\Big((\Pi^0)^{\circ 3}(a_1,a_2,a_3,g)\, ,\,h\Big)   \\
&\qquad- \Big\{(\Pi^0)^{\circ 2}(a_1,a_2,a_3)\,\Pi^0(g,h) + \big(\Pi^0_{a_1}a_2\big)\,\Pi^0\Big(\Pi^0_{\mathscr{C} a_3} g, h\Big) + a_1\,\Pi^0\Big(\Pi^0_{\mathscr{C}\Pi^0_{\mathscr{C} a_2} a_3} g\,,\, h\Big)\Big\}
\end{split}
\end{equation*}
defines a continuous map from $\mcC^{\alpha_1}\times \mcC^{\alpha_2}\times\mcC^{\alpha_3}\times\mcC^{\mu}\times\mcC^{\nu}$ to $\mcC^{\alpha_1+\alpha_2+\alpha_3+\mu+\nu}$, under the above conditions on the regularity exponents.
We then let the reader to complete and adapt the proof in the full parabolic setting. \\
To see how the second term in the expansion arises, use formula \eqref{eq:formulec} for the corrector and write 
\begin{align*}
\Big\{\Pi^0\Big((\Pi^0)^{\circ 3}(a_1,a_2,a_3,g)\, ,\,h\Big) - (\Pi^0)^{\circ 2}(a_1,a_2,a_3)\,\Pi^0(g,h)\Big\}(x)  & = {\sf C}^0\Big((\Pi^0)^{\circ 2}(a_1,a_2,a_3), g\,;h\Big)(x)   \\
& = \Pi^0\Big(\Pi^0_{\mathscr{C}(\Pi^0)^{\circ 2}(a_1,a_2,a_3)} g, h\Big)(x).
\end{align*}
Note that since 
$$
\Pi^0_{a_1}a_2 = \Big(\Pi^0_{a_1}a_2\Big)(x) + \mathscr{C} \Pi^0_{a_1}a_2,
$$
we have the identity
$$
\mathscr{C}(\Pi^0)^{\circ 2}(a_1,a_2,a_3) = \big(\Pi^0_{a_1}a_2\big)(x)\,\mathscr{C} a_3 + \mathscr{C}\Pi^0_{\mathscr{C} \Pi^0_{a_1}a_2}a_3.
$$
It follows that $\Pi^0\Big((\Pi^0)^{\circ 3}(a_1,a_2,a_3,g)\, ,\,h\Big)$ is equal to
\begin{equation*}
(\Pi^0)^{\circ 2}(a_1,a_2,a_3)\,\Pi^0(g,h) + \big(\Pi^0_{a_1} a_2\big)\,\Pi^0\Big(\Pi^0_{\mathscr{C} a_3} g, h\Big) + \Pi^0\Big(\Pi^0_{\mathscr{C}\Pi^0_{\mathscr{C} \Pi^0_{a_1} a_2} a_3} g\,,\, h\Big).
\end{equation*}
Writing $a_1 = a_1(x) + \mathscr{C} a_1$, in the above expression for the remainder yields the formula
\begin{equation}\label{EqCWithScrC}
{\sf C}(a_1,a_2,a_3,g,h) = \Pi^0\Big(\Pi^0_{\mathscr{C}\Pi^0_{\mathscr{C} \Pi^0_{\mathscr{C} a_1} a_2} a_3} g\,,\, h\Big).
\end{equation}
The fact that it defines a $\big(\alpha_1+\alpha_2+\alpha_3+\mu+\nu\big)$-H\"older function if this exponent is positive can be seen as follows. For every $x$ we have 
\begin{align*}
\Pi^0\Big(\Pi^0_{\mathscr{C}\Pi^0_{\mathscr{C} \Pi^0_{\mathscr{C} a_1} a_2} a_3} g\,,\, h\Big)(x) & \simeq \sum_{k} \Delta_k\Big( \Pi^0_{\mathscr{C}\Pi^0_{\mathscr{C} \Pi^0_{\mathscr{C} a_1} a_2} a_3} g \Big)(x) \, \Delta_k (h)(x) \\
& \simeq \sum_{k} S_{k-2}\big(\mathscr{C}\Pi^0_{\mathscr{C} \Pi^0_{\mathscr{C} a_1} a_2 } a_3 \big)(x) \,\Delta_k(g)(x) \,\Delta_k(h)(x) \\
& \simeq \sum_{k} \Delta_{k-2}\big(\Pi^0_{\mathscr{C} \Pi^0_{\mathscr{C} a_1} a_2 } a_3\big)(x) \,\Delta_k(g)(x) \,\Delta_k(h)(x),
\end{align*}
where we used \eqref{eq:appox}. Iterating the reasoning, we get 
\begin{align}
\label{eq::ouf}
\Pi^0\Big(\Pi^0_{\mathscr{C}\Pi^0_{\mathscr{C} \Pi^0_{\mathscr{C} a_1} a_2 } a_3 } g\,,\, h\Big)(x) & \simeq \sum_{k} \Delta_{k-6}(a_1)(x) \,\Delta_{k-4}(a_2)(x)\, \Delta_{k-2}(a_3)(x) \,\Delta_k(g)(x) \,\Delta_k(h)(x) 
\end{align}
and so since $(\alpha_1+\alpha_2+\alpha_3+\mu+\nu)$ is non-negative, and setting
$$
m:=\|a_1\|_{C^{\alpha_1}} \|a_2\|_{C^{\alpha_2}} \|a_3\|_{C^{\beta}} \|g\|_{C^{\alpha_3}} \|h\|_{C^\nu},
$$
we conclude that
\begin{align*}
\left|\Pi^0\Big(\Pi^0_{\mathscr{C}\Pi^0_{\mathscr{C} \Pi^0_{\mathscr{C} a_1} a_2 } a_3 } g \,,\, h\Big)(x)\right| & \lesssim m \sum_{k} 2^{-k(\alpha_1+\alpha_2+\beta+\alpha_3+\nu)}   \\
& \lesssim m,
\end{align*}
uniformly in $x$, which yields that the main quantity defines a bounded function. Using \eqref{eq::ouf}, we can also obtain its H\"older character. For $x\neq y$, we have
\begin{align*}
& \left|\Pi^0\Big(\Pi^0_{\mathscr{C}\Pi^0_{\mathscr{C} \Pi^0_{\mathscr{C} a_1} a_2} a_3} g \,,\, h\Big)(x) - \Pi^0\Big(\Pi^0_{\mathscr{C}\Pi^0_{\mathscr{C} \Pi^0_{\mathscr{C} a_1} a_2} a_3} g \,,\, h\Big)(y)\right| \\
& \lesssim \sum_{k} \Big|\Delta_{k-6}(a_1)(x) \Delta_{k-4}(a_2)(x) \Delta_{k-2}(a_3)(x) \Delta_k(g)(x) \Delta_k(h)(x)   \\
&\hspace{1cm}-\Delta_{k-6}(a_1)(y) \Delta_{k-4}(a_2)(y) \Delta_{k-2}(a_3)(y) \Delta_k(g)(y) \Delta_k(h)(y) \Big| \\
& \lesssim m\,\left(\sum_{1\leq 2^k |x-y|}\ 2^{-k(\alpha_1+\alpha_2+\beta+\alpha_3+\nu)} + \sum_{1 \geq 2^k |x-y|} |x-y|\,2^{k-k(\alpha_1+\alpha_2+\beta+\alpha_3+\nu)}\right)
\\
& \lesssim m\,|x-y|^{\alpha_1+\alpha_2+\beta+\alpha_3+\nu};
\end{align*}
in the second sum, over $1 \geq 2^k |x-y|$, we have used the finite increment theorem together with the fact that differentiating one operator $\Delta_k$ is equivalent to multiplying it by $2^k$, together with the condition $(\alpha_1+\alpha_2+\alpha_3+\mu+\nu) \in (0,1)$.
\end{Dem}

\medskip

It is also necessary and possible to 'expand' the corrector $\sf C$ simultaneously on its first two arguments. Define
$$
{\sf C}\big((a_1,a_2),(b_1,b_2),h\big) := {\sf C}\big((a_1,a_2), {\sf P}_{b_1}b_2 ,h\big) - {\sf P}_{b_1} {\sf C}\big((a_1,a_2),b_2,h\big).
$$
We let the reader write this operator in terms of the outer centering operator $\mathscr{C}$, like in identity \eqref{EqCWithScrC}.

\medskip

\begin{prop}
Assume that $(\alpha_1+\alpha_2+\mu_1+\mu_2)<3$, and 
   \begin{equation*}
   \begin{split}
   &(\alpha_2+\mu_2+\nu) < 0, \quad \alpha_1+\mu_2+\nu<0  \\
   &(\alpha_1+\alpha_2+\mu_1+\mu_2+\nu) \in (0,1) \quad \textrm{and} \quad \mu_2\in(0,1).
   \end{split}
   \end{equation*}      
   Then the $5$-linear iterated corrector
   \begin{equation*}
   \begin{split}   
   \mcC^{\alpha_1}\times\mcC^{\alpha_2}\times\mcC^{\mu_1}\times \mcC^{\mu_2}\times\mcC^{\nu} &\rightarrow \mcC^{\alpha_1+\alpha_2+\mu_1+\mu_2+\nu}   \\
   \big(a_1,a_2,b_1,b_2,h\big) &\mapsto {\sf C}\big((a_1,a_2),(b_1,b_2),h\big)
   \end{split}
   \end{equation*}    
   is continuous.
\end{prop}

\medskip

\begin{Dem}
As in the proof of Theorem \ref{ThmIteratedCorrector}, observe that (in the flat torus setting) we have
\begin{align*}
{\sf C}\big((a_1,a_2), {\sf P}_{b_1}b_2 ,h\big) & \simeq \sum_{k} \Delta_{k-4}(a_1) \Delta_{k-2}(a_2) \Delta_k({\sf P}_{b_1} b_2) \Delta_k(h) \\
& \simeq \sum_k \Delta_{k-4}(a_1) \Delta_{k-2}(a_2) S_{k-2}(b_1)\Delta_k(b_2) \Delta_k(h).
\end{align*}
At the same time,
\begin{align*}
{\sf P}_{b_1} {\sf C}\big((a_1,a_2),b_2,h\big)& = \sum_{k} S_{k-2}(b_1) \Delta_k\Big( {\sf C}\big((a_1,a_2),b_2,h\big) \Big) \\
& = \sum_{k} S_{k-2}(b_1) \Delta_k\Big(\sum_{j}  \Delta_{j-4}(a_1) \Delta_{j-2}(a_2) \Delta_j(b_2) \Delta_j(h)\Big). 
\end{align*}
In the previous sum, the parameters $k,j$ have to be equivalent in order to have a non-vanishing contribution. With the normalization $\sum \Delta_k=\textrm{Id}$, we obtain
\begin{align*}
(\star)&:={\sf C}\big((a_1,a_2), {\sf P}_{b_1}b_2 ,h\big)(x) - {\sf P}_{b_1} {\sf C}\big((a_1,a_2),b_2,h\big)(x) \\
&\;\simeq \sum_{k,j} \Delta_k\Big( \Delta_{j-4}(a_1) \Delta_{j-2}(a_2) \big(S_{j-2}(b_1)-S_{k-2}(b_1)(x)\big)\Delta_j(b_2) \Delta_j(h)\Big)(x).
\end{align*}
From this decomposition, it is easy to check that for $\ell\geq 0$ and because $\alpha_1+\alpha_2+\mu_1+\mu_2+\nu \in(0,1)$, and $\mu_1\in(0,1)$, we have 
\begin{align*}
\|\Delta_\ell(\star)\|_{L^\infty} & \lesssim \sum_{k,j \gtrsim \ell} 2^{-j(\alpha_1+\alpha_2+\mu_2+\nu)} (2^{-j}+2^{-k})^{\mu_1} m \\
& \lesssim 2^{-\ell(\alpha_1+\alpha_2+\mu_1+\mu_2+\nu)} m,
\end{align*}
with
$$
m:=\|a_1\|_{C^{\alpha_1}} \|a_2\|_{C^{\alpha_2}} \|b_1\|_{C^{\mu_1}} \|b_2\|_{C^{\mu_2}} \|h\|_{C^\nu}.
$$
\end{Dem}

\medskip

The $5$-linear iterated corrector will never appear explicitly in our computations as it will provide a remainder term. The $4$ and $5$-linear \textsf{\textbf{upper iterated correctors}} are defined by the formulas
$$
{\sf C}\big(f,(a_1,a_2), h\big) := {\sf C}\Big(f,\widetilde{\sf P}_{a_1}a_2,h\Big) - a_1\,{\sf C}\big(f,a_2, h\big).
$$
and 
$$
{\sf C}\Big(f,\big(a_1,(a_2,a_3)\big), h\Big) := {\sf C}\Big(f,\big(a_1,\widetilde{\sf P}_{a_2}a_3\big), h\Big) - a_2\, {\sf C}\big(f,(a_1,a_3), h\big).
$$ 
 
\medskip

\begin{thm} \label{thm:upper} \label{ThmUpper}  {\it 
\begin{itemize}
   \item[\textsf{\textbf{(i)}}] Assume $\mu,\alpha_1\in (0,1)$, and the exponents $(\mu+\alpha_1+\nu)$ and $(\mu+\alpha_2+\nu)$ are negative, and
   $$
   (\mu+\alpha_1+\alpha_2+\nu) \in (0,1). 
   $$
   Then the $4$-linear upper iterated corrector 
   \begin{equation*}
   \begin{split}
   \mcC^\mu\times\mcC^{\alpha_1}\times\mcC^{\alpha_2}\times\mcC^{\nu} &\rightarrow \mcC^{\mu+\alpha_1+\alpha_2+\nu}   \\
   \big(f,a_1,a_2,h\big) &\mapsto {\sf C}\big(f,(a_1,a_2),h\big)  
   \end{split}
   \end{equation*}   
   is continuous.   \vspace{0.15cm}
   
   \item[\textsf{\textbf{(ii)}}] Assume $\mu,\alpha_1,\alpha_2\in (0,1)$, and the exponents $(\mu+\alpha_i+\nu)$ are all negative, for $1\leq i\leq 3$, and
   $$
   (\mu+\alpha_1+\alpha_2+\alpha_3+\nu) \in (0,1).
   $$
   Then the $5$-linear upper iterated corrector 
   \begin{equation*}
   \begin{split}
   \mcC^\mu\times\mcC^{\alpha_1}\times\mcC^{\alpha_2}\times\mcC^{\alpha_3}\times\mcC^{\nu} &\rightarrow \mcC^{\mu+\alpha_1+\alpha_2+\alpha_3+\nu}   \\
   \big(f,a_1,a_2,a_3,h\big) &\mapsto {\sf C}\big(f,(a_1,(a_2,a_3)),h\big)  
   \end{split}
   \end{equation*}   
    is continuous. 
\end{itemize}   }
\end{thm}

\medskip

\begin{Dem}
We only sketch the proof of the continuity result of the $4$-linear operator in the model case of the time-independent setting of the flat torus, and rely on formula \eqref{EqDiagonalTerm} for the diagonal operator $\Pi^0(\cdot,\cdot)$ for that purpose. See Proposition \ref{prop:C2} in Appendix \ref{Subsectionapp2} for a fully detailed proof in the parabolic setting. In the present setting, the quantity ${\sf C}^0(f\,;a_1,a_2\,;g)$ is then given by a sum of the form
$$
{\sf C}^0\big(f,(a_1,a_2),h\big) = \sum_{\vert i-j\vert\leq 1} \epsilon'_i \, \Delta_jh,
$$
with 
\begin{equation*} 
\epsilon'_i := \Big\{\Delta_i\big({\sf P}_f\big({\sf P}_{a_1}a_2\big)\big) - a_1\Delta_i\big({\sf P}_f a_2\big)\Big\} + f\Big\{a_1\Delta_i a_2 - \Delta_i\big({\sf P}_{a_1}a_2\big)\Big\}  
\end{equation*}
We read on the expression 
{\small \begin{equation*}
\begin{split}
\epsilon'_i(x) &= \int K_i(x,y)\Big\{{\sf P}_f\big({\sf P}_{a_1}a_2\big)(y) - a_1(x)\big({\sf P}_fa_2\big)(y) + (fa_1)(x)a_2(y) - f(x)\big({\sf P}_{a_1}a_2\big)(y)\Big\}\,dy   \\
				  	 &= \int K_i(x,y) \, {\sf P}_{f-f(x){\bf 1}}\Big({\sf P}_{a_1-a_1(x){\bf 1}}a_2\Big)(y)\,dy,
\end{split}
\end{equation*}   } 
that  
$$
\epsilon'_i = \Delta_i\Big({\sf P}_{\mathscr{C} f}\big({\sf P}_{\mathscr{C} a_1}a_2\big)\Big)
$$
has $L^\infty$-norm of order $2^{-i(\mu+\alpha_1+\alpha_2)}$, as a consequence of \eqref{eq:appox}. The proof is then not fully completed, since the block $\epsilon_i' \Delta_ih$ is not perfectly localized in frequency at scale $2^i$, so an extra decomposition is necessary. The same thing happens in the proof of Theorem \ref{ThmCorrectors}. We do not give the details here and refer the reader to the proof of Proposition \ref{prop:C2} in Appendix \ref{SectionAppendixContinuity}.
\end{Dem}

\bigskip

\textbf{\textsf{3{\boldmath$.$}1{\boldmath$.$}4  Commutator $\sf D${\boldmath$.$}}} Define on the space $L^\infty$ of bounded measurable functions on the parabolic space $\mcM$ the \textsf{\textbf{commutator}} {\sf D} as the operator 
$$
{\sf D}(f,g,h) := {\sf \Pi}\big(\widetilde{\sf P}_f g, h\Big) - {\sf P}_f\Big({\sf \Pi}(g,h)\Big),
$$

\smallskip

\begin{thm}
\begin{itemize}
   \item[\textsf{\textbf{(a)}}] For positive regularity exponents $\alpha,\beta$ and $\gamma$, the commutator ${\sf D}$ is continuous from $\mcC^\alpha\times\mcC^\beta\times\mcC^\gamma$ to $\mcC^{\alpha+\beta+\gamma}$.   \vspace{0.15cm}
   
   \item[\textsf{\textbf{(b)}}] The commutator ${\sf D}$ is bounded from $\mcC^{\alpha} \times \mcC^\beta \times \mcC^\gamma$ into $\mcC^{\alpha+\beta+\gamma}$ for $\alpha\in(-1,0)$ as soon as $\alpha+\beta+\gamma>0$.
\end{itemize}   
\end{thm}

\medskip

\begin{Dem}
Assume we are working in the time-independent setting of the $d$-dimensional torus, with the operator
$$
{\sf D}^0(f,g,h) := \Pi^0\Big(\Pi^0_f g, h\Big) - \Pi^0_f\Big(\Pi^0(g,h)\Big).
$$
\textsf{\textbf{(a)}} We refer the reader to  Proposition \ref{prop:C1}, in Appendix \ref{Subsectionapp1}, for a full proof of the regularity statement for the commutator ${\sf D}$. We simply mention here that in the special case of ${\sf D}^0$, the regularity estimate comes from the following identity
\begin{align}
\Delta_k\big({\sf D}^0(f,g,h)\big) & = \sum_{\ell \geq k-2} \Delta_k\Big(\Delta_\ell(g) S_{\ell-2}(f) \Delta_\ell(h)\Big) - S_{k-2}(f) \Delta_k\Big(\Delta_\ell(g) \Delta_\ell(h)\Big) \nonumber \\
& \simeq \sum_{\ell \geq k-2} \Delta_k\Big(\Delta_\ell(g) S_{\ell-2}(f) \Delta_\ell(h) - S_{k-2}(f)\Delta_\ell(g) \Delta_\ell(h)\Big)  \label{eq:DkD}
\end{align}

\medskip

\textsf{\textbf{(b)}} This case is easy and do not use the 'difference' structure in the commutator. Indeed since $\alpha<0$ then $\alpha+\beta+\gamma>0$ implies $\beta+\gamma>0$ and so by using the boundedness of the paraproducts and those of the resonant part, it directly comes
$$  
{\sf \Pi}\big(\widetilde{\sf P}_f g, h\Big) \in \mcC^{\alpha+\beta+\gamma} \quad \textrm{and} \quad {\sf P}_f\Big({\sf \Pi}(g,h)\Big) \in \mcC^{\alpha+\beta+\gamma}.
$$
\end{Dem}

\bigskip

\textbf{\textsf{3{\boldmath$.$}1{\boldmath$.$}5  Back to the high order paracontrolled expansion \eqref{EqTaylorExpansion}{\boldmath$.$}}} The high order paracontrolled expansion \eqref{EqTaylorExpansion} might seem far from a possibly expected ordinary Taylor-type expansion, such as it appears for instance in regularity structures \cite{HairerRegularity}. The difference is not that big, as we see in this section by looking at the resonant term ${\sf \Pi}(f(u),\zeta)$. This intermezzo has a formal character. Recall we denote by ${\sf C}^\circ, {\sf D}^\circ$ the corrector, commutator and their iterates, built from the resonant operator and the usual paraproduct ${\sf \Pi}, {\sf P}$ instead of the modified paraproduct $\widetilde{\sf P}$; they enjoy the same continuity properties as $\sf C$ and $\sf D$. We consider in this section the case where $u\in\mcC^\alpha$, with $\frac{2}{5}<\alpha<\frac{1}{2}$, and for a function $g\in C^3_b$, we write the second order version of formula \eqref{EqTaylorExpansion}, giving the paracontrolled expansion of $g(u)$,  under the form
$$
g(u) =: {\sf P}_{g'(u)}u + \frac{1}{2}\Big\{{\sf P}_{g^{(2)}(u)}u^2-{\sf P}_{g^{(2)}(u)u}u\Big\} + g(u)^{(1)};
$$
the exponent $(1)$ refering to the fact that $g(u)^{(1)}\in\mcC^{3\alpha}\subset\mcC^1$, so it can be differentiated in the space direction. We work in our general parabolic setting over a Riemannian manifold so the refined corrector involves the additional term $\sum_{\ell=1}^{\ell_0} \gamma_\ell(V_\ell f){\sf \Pi}_{(1)}^\ell(g,h)$, from Definition \ref{DefnRefinedCorrector}. The reader can think of the time-independent flat torus setting, where the additional term in the definition of the refined corrector is simply $f'{\sf \Pi}_{(1)}^0(g,h)$, with $'$ denoting space derivative.

\medskip

\noindent \textsf{\textbf{Notation.}} \textit{Given any $\beta\in\bbR$, denote by $(\beta)$ an element of $\mcC^\beta$ whose precise definition is unimportant for the reasoning, and whose only noticeable property is its regularity. Its definition may change from line to line.}

\medskip

Pick now a function $f\in C^4_b$. Writing
\begin{equation*}
\begin{split}
&u^2 = 2{\sf P}_uu + {\sf \Pi}(u,u),   \\
&u^3 = {\sf P}_{u^2}u + 2{\sf P}_u({\sf P}_uu) + {\sf P}_u\big({\sf \Pi}(u,u)\big) + 2{\sf \Pi}\big(u,{\sf P}_uu\big) + {\sf \Pi}\big(u,{\sf \Pi}(u,u)\big),
\end{split}
\end{equation*}
the third order paracontrolled expansion formula \eqref{EqTaylorExpansion} writes
\begin{equation*}
\begin{split}
f(u) &= {\sf P}_{a_1}u + {\sf P}_{a_2}{\sf P}_uu + {\sf P}_{a_3}\big({\sf P}_u({\sf P}_uu)\big) + {\sf P}_{a_4}\big({\sf P}_{u^2}u\big)   \\
	   &\quad+ {\sf P}_{b_1}{\sf \Pi}(u,u) + {\sf P}_{b_2}\big({\sf P}_u{\sf \Pi}(u,u)\big) + {\sf P}_{b_3}\big({\sf \Pi}(u,{\sf P}_uu)\big)   \\
	   &\quad+ {\sf P}_c\Big({\sf \Pi}\big(u,{\sf \Pi}(u,u)\big)\Big) + f(u)^\sharp,
\end{split}
\end{equation*}
where $f(u)^\sharp\in\mcC^{4\alpha}$, and
{\small \begin{equation*}
\begin{split}
a_1 &:= f'(u) - uf^{(2)}(u)+\frac{1}{2}\,u^2f^{(3)}(u),   \\
a_2 &:= f^{(2)}(u)-uf^{(3)}(u),   \\ 
a_3 &:= \frac{1}{3}\,f^{(3)}(u),   \\
a_4 &:= \frac{1}{6}\, f^{(3)}(u),   \\
b_1 &:= \frac{1}{2}\,\Big(f^{(2)}(u) - f^{(3)}(u)\Big),  \\
b_2 &:=  \frac{1}{6}\, f^{(3)}(u),   \\ 
b_3 &:=  \frac{1}{3}\, f^{(3)}(u),   \\
c &:=  \frac{1}{6}\, f^{(3)}(u).
\end{split}
\end{equation*}   }
Plugging this formula inside ${\sf \Pi}\big(f(u),\zeta\big)$ gives for it the expression
\begin{equation*}
\begin{split}
									 & a_1\,{\sf \Pi}(u,\zeta) + {\sf C}^\circ(a_1,u,\zeta)   \\
									 &\quad+ a_2u\,{\sf \Pi}(u,\zeta) + a_2 {\sf C}^\circ(u,u,\zeta) + u\,{\sf C}^\circ(a_2,u,\zeta) + {\sf C}^\circ(a_2,(u,u),\zeta)   \\
									 &\quad+ a_3u^2\,{\sf \Pi}(u,\zeta) + 2a_3u\,{\sf C}^\circ(u,u,\zeta) + u^2\,{\sf C}^\circ(a_3,u,\zeta) + 2u\,{\sf C}^\circ\big(a_3,(u,u),\zeta\big) + a_3\,{\sf C}^\circ\big(u,(u,u),\zeta\big)   \\
									 &\quad+ a_4u^2\,{\sf \Pi}(u,\zeta) + a_4\,{\sf C}^\circ(u^2,u,\zeta) + u^2\,{\sf C}^\circ(a_4,u,\zeta) + {\sf C}^\circ\big(a_4,(u^2,u),\zeta\big)   \\
									 &\quad+ b_1\,{\sf \Pi}\big({\sf \Pi}(u,u),\zeta\big) + {\sf C}^\circ\big(b_1,{\sf \Pi}(u,u),\zeta\big)   \\
									 &\quad+ b_2u\,{\sf \Pi}\big({\sf \Pi}(u,u),\zeta\big) + b_2\,{\sf C}^\circ\big(u,{\sf \Pi}(u,u),\zeta\big) + u\,{\sf C}^\circ\big(b_2,{\sf \Pi}(u,u),\zeta\big)   \\
									 &\quad+ b_3u\,{\sf \Pi}\big({\sf \Pi}(u,u),\zeta\big) + b_3\,{\sf C}^\circ\big(u,{\sf \Pi}(u,u),\zeta\big) + b_3\,{\sf \Pi}\big({\sf D}^\circ(u,u,u),\zeta\big) + u\,{\sf C}^\circ\big(b_3,{\sf \Pi}(u,u),\zeta\big)   \\
									 &\quad+ c\,{\sf \Pi}\big({\sf \Pi}(u,{\sf \Pi}(u,u)),\zeta\big) + (5\alpha-2).
\end{split}
\end{equation*}
To get this expression in terms only of \textsl{primitive} quantities involving ill-defined terms where $u$ appears rather than a function of $u$, use the identity
\begin{equation*} \begin{split}
{\sf C}^\circ\big(a_i(u),u,\zeta\big) &= a_i'(u)\,{\sf C}^\circ(u,u,\zeta) + a_i^{(2)}(u)\,{\sf C}^\circ\big((u,u),u,\zeta\big) + \frac{1}{2}a_i^{(2)}(u)\,{\sf C}^\circ\big({\sf \Pi}(u,u),u,\zeta\big)   \\
&\quad+ (5\alpha-2) + {\sf C}^\circ\big(a_i(u)^{(1)},u,\zeta\big)   \\
&= a_i'(u)\,{\sf C}^\circ(u,u,\zeta) + a_i^{(2)}(u)\,{\sf C}^\circ\big((u,u),u,\zeta\big) + \frac{1}{2}a_i^{(2)}(u)\,{\sf C}^\circ\big({\sf \Pi}(u,u),u,\zeta\big)   \\
&\quad+ (5\alpha-2) + \sum_{\ell=1}^{\ell_0}\gamma_\ell V_\ell\big(a_i(u)^{(1)}\big)\,{\sf \Pi}_{(1)}^\ell(u,\zeta).
\end{split} \end{equation*}
We expand ${\sf C}^\circ(u^2,u,\zeta)$ using the identity $u^2=2{\sf P}_uu+{\sf \Pi}(u,u)$ and the iterated corrector. After simplification, one gets the following \textit{multiplicative decomposition} for ${\sf \Pi}\big(f(u),\zeta\big)$
\begin{equation}
\label{EqDecompositionPrimitiveResonant}
\begin{split}
&f'(u)\,{\sf \Pi}(u,\zeta) + f^{(2)}(u)\,\left\{{\sf C}^\circ(u,u,\zeta) + \frac{1}{2}\,{\sf \Pi}\big({\sf \Pi}(u,u),\zeta\big)\right\}   \\
&+ f^{(3)}(u)\,\left\{\frac{1}{3}\,{\sf C}^\circ\big(u,(u,u),\zeta\big) + \frac{1}{3}\,{\sf C}^\circ\big((u,u),u,\zeta\big) + \frac{1}{6}\,{\sf C}^\circ\big(u,{\sf \Pi}(u,u),\zeta\big) + \frac{1}{6}\,{\sf C}^\circ\big({\sf \Pi}(u,u),u,\zeta\big) \right.  \\
&\hspace{1.8cm}+ \left. \frac{1}{3}\,{\sf \Pi}\big({\sf D}^\circ(u,u,u),\zeta\big) + \frac{1}{6}\,{\sf \Pi}\big({\sf \Pi}(u,{\sf \Pi}(u,u)),\zeta\big) \right\}   \\
&+ \sum_{\ell=1}^{\ell_0}(\star)_\ell\,{\sf \Pi}^\ell_{(1)}(u,\zeta) + (5\alpha-2),
\end{split} \end{equation}
with 
\begin{equation}   \label{EqStarEll}
(\star)_\ell := \gamma_\ell V_\ell\big(a_1(u)^{(1)}\big) + u\,\gamma_\ell V_\ell\big(a_2(u)^{(1)}\big) + u^2\,\gamma_\ell V_\ell\Big(a_3(u)^{(1)}+a_4(u)^{(1)}\Big) \in \mcC^{3\alpha-1}.
\end{equation}
In the time-independent setting of the flat torus, this reduces to 
$$
\sum_{\ell=1}^{\ell_0}(\star)_\ell\,{\sf \Pi}^\ell_{(1)}(u,\zeta) = \Big\{\big(a_1(u)^{(1)}\big)' + u\big(a_2(u)^{(1)}\big)' + u^2\Big(a_3(u)^{(1)}+a_4(u)^{(1)}\Big)'\Big\}\,{\sf \Pi}_{(1)}(u,\zeta).
$$   
The term ${\sf \Pi}^\ell_{(1)}(u,\zeta)$ has formal regularity $2\alpha-1$, slightly less than $0$, and the product with the above bracket term $\{\cdot\}$ is well-posed provided ${\sf \Pi}^\ell_{(1)}(u,\zeta)$ can be given meaning as an element of regularity $2\alpha-1$.   

\medskip

Write each expression ${\sf \Pi}\big(u^k,\zeta\big)$ in multiplicative form and asign to a product the regularity of its term of lowest regularity. For an expression $B$ in multiplicative form, we can then denote by $\lfloor B\rfloor_{k\alpha-2}$ the part of $B$ that is of formal regularity $(k\alpha-2)$. So we have for instance
$$
\Big\lfloor{\sf \Pi}(u,\zeta)\Big\rfloor_{2\alpha-2} = {\sf \Pi}(u,\zeta),
$$
and 
$$
{\sf \Pi}({\sf P}_uu,\zeta) = u{\sf \Pi}(u,\zeta) + {\sf C}^\circ(u,u,\zeta),
$$
with 
$$
\Big\lfloor {\sf \Pi}({\sf P}_uu,\zeta)\Big\rfloor_{2\alpha-2} = u{\sf \Pi}(u,\zeta)
$$
and 
$$
\Big\lfloor {\sf \Pi}({\sf P}_uu,\zeta)\Big\rfloor_{3\alpha-2} = {\sf C}^\circ(u,u,\zeta).
$$
In those terms, formula \eqref{EqDecompositionPrimitiveResonant} for ${\sf \Pi}(f(u),\zeta)$ takes the following Taylor-type form
{\small \begin{equation}
\label{EqFirstFundamentalFormula}
\begin{split}
{\sf \Pi}\big(f(u),\zeta\big) &= f'(u)\,\Big\lfloor{\sf \Pi}(u,\zeta)\Big\rfloor_{2\alpha-2} + f^{(2)}(u)\,\left\lfloor {\sf \Pi}\Big(\frac{u^2}{2!},\zeta\Big)\right\rfloor_{3\alpha-2} + f^{(3)}(u)\,\left\lfloor {\sf \Pi}\Big(\frac{u^3}{3!},\zeta\Big)\right\rfloor_{4\alpha-2}   \\
&\hspace{0.3cm}+ \sum_{\ell=1}^{\ell_0}\gamma_\ell\big(V_\ell(a_i(u)^{(1)})\big){\sf \Pi}^\ell_{(1)}(u,\zeta) + (5\alpha-2).
\end{split}
\end{equation}   }

\bigskip

\textbf{\textsf{3{\boldmath$.$}1{\boldmath$.$}6  Iterated commutators{\boldmath$.$}}} In addition to the above continuity properties for the iterated correctors, we also need the following continuity result on iterated commutator operators $\sf D$.

\medskip

\begin{prop}   \label{PropIteratedD}   {\it 
\begin{itemize}
   \item[\textcolor{gray}{$\bullet$}] Given positive regularity exponents $\alpha_1,\alpha_2, \gamma, \delta$, the formulas
\begin{equation*}
\begin{split}
\mcC^{\alpha_1}\times\mcC^{\alpha_2}\times\mcC^\gamma\times\mcC^\delta &\rightarrow\mcC^{\alpha_1+\alpha_2+\gamma+\delta}   \\
(a_1,a_2,g,h) &\mapsto {\sf D}\Big(\widetilde{\sf P}_{a_1}a_2, g,h\Big) - {\sf P}_{a_1}{\sf D}(a_2,g,h),
\end{split}
\end{equation*}
for the lower iterated commutator, and 
\begin{equation*}
\begin{split}
\mcC^\gamma\times\mcC^{\alpha_1}\times\mcC^{\alpha_2}\times\mcC^\delta &\rightarrow\mcC^{\gamma+\alpha_1+\alpha_2+\delta}   \\
(f,a_1,a_2,h) &\mapsto {\sf D}\Big(f, \widetilde{\sf P}_{a_1}a_2,h\Big) - {\sf P}_{a_1}{\sf D}(f,a_2,h),
\end{split}
\end{equation*}
for the upper iterated commutator, define continuous operators. The result also holds true if $\alpha_1\in(0,1)$ and $-1<\alpha_2<0$, with $\alpha_2+\gamma+\delta>0$.   \vspace{0.15cm}
   
   \item[\textcolor{gray}{$\bullet$}] Fix $\alpha,\gamma,\delta_1>0$ and $\delta_2\in(0,1)$. The high frequency commutator 
   \begin{equation*}
   \begin{split}
   \mcC^{\alpha}\times\mcC^{\gamma}\times\mcC^{\delta_1}\times \mcC^{\delta_2} &\rightarrow\mcC^{\alpha+\gamma+\delta_1+\delta_2}      \\
   (f,g,h_1,h_2) &\mapsto {\sf D}\Big(f,g, \widetilde{\sf P}_{h_1}h_2\Big) - {\sf P}_{h_1}{\sf D}(f,g,h_2),
   \end{split}
   \end{equation*}
   is bounded. This continuity result also holds true if $\alpha\in(-1,0)$, provided $\alpha+\gamma+\delta_2>0$.
\end{itemize}   }
\end{prop}

\medskip

\begin{Dem}
As in the proof of Theorem \ref{ThmCorrectors}, we analyse in the present proof what happens in the time-independent setting of the $d$-dimensional torus, in the case where we also use $\Pi^0$ instead of $\widetilde{\sf P}$. So we set 
$$
{\sf D}^0(a_1,a_2,g,h) := {\sf D}^0\Big(\Pi^0_{a_1}a_2, g,h\Big) - \Pi^0_{a_1}{\sf D}^0(a_2,g,h)
$$
and have a look at its continuity properties on the spacial H\"older spaces. Using formula \eqref{eq:DkD}, it follows that we have
\begin{align*}
\Delta_k\big(&{\sf D}^0(a_1,a_2,g,h)\big) \simeq \Delta_k\Big({\sf D}^0\big({\sf P}_{a_1}a_2, g,h\big)\Big) - S_{k-2}(a_1)\, \Delta_k \Big({\sf D}^0(a_2,g,h)\Big)   \\
&\hspace{-0.9cm}\simeq \sum_{\ell \geq k-2} \Delta_k\Big\{\Delta_\ell(g) \Delta_\ell(h)\Big(S_{\ell-2}\big({\sf P}_{a_1}a_2\big) - S_{k-2}\big({\sf P}_{a_1}a_2\big) - S_{k-2}(a_1)\big(S_{\ell-2} a_2 - S_{k-2} a_2\big)\Big) \Big\}.
\end{align*}
The quantity inside the brackets is equal to 
\begin{align*}
S_{\ell-2}\big({\sf P}_{a_1}a_2\big) - &S_{k-2}\big({\sf P}_{a_1}a_2\big) - S_{k-2}(a_1)\big(S_{\ell-2}(a_2)-S_{k-2}(a_2)\big)   \\
&= \sum_{j=k-1}^{\ell-2} \Delta_j\big({\sf P}_{a_1}a_2\big) - S_{k-2}(a_1) \Delta_j(a_2) \\
&\simeq \sum_{j=k-1}^{\ell-2} S_{j-2}(a_1) \Delta_j(a_2) - S_{k-2}(a_1) \Delta_j(a_2) \\
&\simeq \sum_{j=k-1}^{\ell-2} \big(S_{j-2} a_1 - S_{k-2} a_1\big)\Delta_j(a_2),
\end{align*}
which is bounded in $L^\infty$ by
$$ 
\sum_{j=k+1}^{\ell} 2^{-k\alpha} \|a_1\|_{\mcC^{\alpha_1}} 2^{-j\beta} \|a_2\|_{\mcC^{\alpha_2}} \lesssim 2^{-k(\alpha_1+\alpha_2)} \|a_1\|_{C^{\alpha_1}}\|a_2\|_{C^{\alpha_2}}.
$$
This estimate allows us to conclude that 
$$
\big\| \Delta_k\big({\sf D}^0(a_1,a_2,g,h)\big) \big\|_\infty \lesssim 2^{-k(\alpha_1+\alpha_2+\gamma+\delta)} \|a_1\|_{C^{\alpha_1}} \|a_2\|_{C^{\alpha_2}}\|g\|_{C^\gamma}\|h\|_{C^\delta},
$$ 
uniformly in $k$, which proves the continuity result for the $4$-linear operator ${\sf D}^0$. A very similar proof gives the continuity of the simplified version of the upper iterated commutator; we leave the details to the reader.

\smallskip

For the second statement of the first item of the proposition, with $\alpha_2 \in (-1,0)$ we follow the same computations and since we have now
$$ 
\sum_{j=k-1}^{\ell-2} 2^{-k\alpha} \|a_1\|_{\mcC^{\alpha_1}} 2^{-j\alpha_2} \|a_2\|_{\mcC^{\alpha_2}} \lesssim 2^{-k\alpha_1 - \ell \alpha_2} \|a_1\|_{C^{\alpha_1}}\|a_2\|_{C^{\alpha_2}},
$$
then
\begin{align*}
\big\| \Delta_k\big({\sf D}^0(a_1,a_2,g,h)\big) \big\|_\infty & \lesssim \sum_{\ell\geq k}  2^{-k\alpha_1} 2^{-\ell(\alpha_2+\gamma+\delta)} \|a_1|_{C^{\alpha_1}} \|a_2\|_{C^{\alpha_2}}\|g\|_{C^\gamma}\|h\|_{C^\delta} \\
& \lesssim 2^{-\ell(\alpha+\beta+ \gamma+\delta)} \|a_1|_{C^{\alpha_1}} \|a_2\|_{C^{\alpha_2}}\|g\|_{C^\gamma}\|h\|_{C^\delta},
\end{align*} 
since $\beta+\gamma+\delta>0$.

\smallskip

For the second item of the proposition, the same reasoning can be applied by observing that now
\begin{align*}
& \Delta_k\Big[{\sf D}\Big(f,g, \widetilde{\sf P}_{h_1}h_2\Big) - {\sf P}_{h_1}{\sf D}(f,g,h_2)\big] \\
& \simeq \sum_{\ell \geq k} \Delta_k\Big[ \Delta_j(g) \big(S_{\ell-2} f-S_{k-2}f\big) \Delta_\ell h_2 \big(S_{\ell-2} h_1 - S_{k-2} h_1\big) \Big];
\end{align*}
we conclude by using the regularity of the four functions.   
\end{Dem}

\bigskip

\subsection[\hspace{-0.3cm} Iterated paraproducts]{Iterated paraproducts}
\label{SubSectionIteratedParaPdcts}

The operators $\sf C$ and $\sf D$ introduced in Section \ref{SectionToolKit} are used to analyse ill-defined products. The operators $\sf R$ and $\sf S$ that we introduce in this section are used to write down the different terms that appear from using the $\sf C$ and $\sf D$ operators in the paraproduct form required to apply the fixed point strategy for the analysis of Equation \eqref{EqCanonicalSPDE} -- Step 3 in Section \ref{SubsectionOverview}. As a motivating example, let us set ourselves the task of writing the paraproduct ${\sf P}_\zeta u$ under the form $\sum_{i=1}^4 {\sf P}_{v_j} Y_j + (5\alpha-2)$, for some $v_j\in\mcC^\alpha$, some $Y_j$ that depend only on the $Z_i$ and the noise, and a remainder $(5\alpha-2)$ of regularity $5\alpha-2$, assuming if necessary that the $u_i$ also have a paracontrolled expansion up to some $u_i$-dependent order. We repeatedly use for that purpose the $F$-type decomposition \eqref{EqFundExpansion2} on the paraproduct map and its iterates. We use here the sign $\sim$ for the equality of functions, up to some term of positive regularity. We have
$$
{\sf P}_\zeta u \sim \sum_{i=1}^3{\sf P}_\zeta {\sf P}_{u_i}Z_i,
$$
and setting
$$
{\sf S}^\circ(\zeta, u_i,Z_i) := {\sf P}_\zeta {\sf P}_{u_i}Z_i - {\sf P}_{u_i} {\sf P}_\zeta Z_i,
$$
we have
$$
{\sf P}_\zeta u \sim \sum_{i=1}^3 \Big({\sf P}_{u_i}{\sf P}_\zeta Z_i + {\sf S}^\circ(\zeta, u_i,Z_i)\Big).
$$
Theorem \ref{ThmCommutatorParaproducts} below shows that the 'swap' operator ${\sf S}^\circ$ sends continuously $\mcC^{\alpha-2}\times\mcC^{\alpha_1}\times\mcC^{\alpha_2}$ into $\mcC^{\alpha-2+\alpha_1+\alpha_2}$. So we have
$$
{\sf P}_\zeta u \sim \sum_{i=1}^3 {\sf P}_{u_i}{\sf P}_\zeta Z_i + {\sf S}^\circ(\zeta, u_1,Z_1) + {\sf S}^\circ(\zeta, u_2,Z_2).
$$
The first three terms on the right hand side have a good form. To analyse the two ${\sf S}^\circ(\zeta,\cdot,\cdot)$ terms, we do an ${\sf F}$-type expansion on each of them, taking profit from the paracontrolled expansion
$$
u_1 = \sum_{j=1}^2 {\sf P}_{u_{1j}}Z_j + (3\alpha), \qquad u_2={\sf P}_{u_{21}}Z_1 + (2\alpha),
$$
of $u_1$ and $u_2$. The map
$$
{\sf S}^\circ_\zeta\big((a_1,a_2),b\big) := {\sf S}^\circ(\zeta, {\sf P}_{a_1}a_2,b) - {\sf P}_{a_1}{\sf S}^\circ(\zeta, a_2,b)
$$
happens indeed to send $\mcC^{\alpha-2}\times\mcC^{\alpha_1}\times\mcC^{\alpha_2}\times\mcC^\beta$ into $\mcC^{\alpha-2+\alpha_1+\alpha_2+\beta}$ continuously -- Theorem \ref{ThmCommutatorParaproducts}, so we have
\begin{equation}
\label{EqGoodFormIteratedParapdct}
{\sf P}_\zeta u \sim \sum_{i=1}^3{\sf P}_{u_i}({\sf P}_\zeta Z_i) + \sum_{j=1}^2 {\sf P}_{u_{1j}}{\sf S}^\circ(\zeta, Z_j,Z_1) + {\sf P}_{u_{21}}{\sf S}^\circ(\zeta, Z_1,Z_2) + {\sf S}^\circ(\zeta, u_{11},Z_1,Z_1).
\end{equation}
A further ${\sf F}$-type expansion on the last term in the above right hand side does the job.

New things happen when we look at the paraproduct ${\sf P}_\zeta u^2$, as we have to deal with a term
$$
{\sf P}_\zeta{\sf P}_uu = {\sf P}_u{\sf P}_\zeta u + {\sf S}^\circ(\zeta, u,u).
$$
Using \eqref{EqGoodFormIteratedParapdct}, we end up with terms of the form ${\sf P}_u{\sf P}_{u_a}Y_a$. A similar thing happens in the analysis of $S^\circ(\zeta, u,u)$ and ${\sf P}_\zeta{\sf \Pi}(u,u)$. To deal with ${\sf P}_u{\sf P}_{u_a}Y_a$, we use the \textit{merging operator} ${\sf R}^\circ$
$$
{\sf P}_u{\sf P}_{u_a}Y_a =: {\sf P}_{uu_a}Y_a + {\sf R}^\circ(u,u_a,Y_a),
$$
and prove some continuity results and some expansion property of ${\sf R}^\circ$ with respect to its first two arguments.

\medskip

\textsl{In this section, we define and state continuity results for the swap and merging operators $\sf S$ and $\sf R$ defined below. We prove here some of the results in the model setting of the time-independent flat torus and refer the reader to Appendix \ref{Subsectionapp2} for the description of how things work in the parabolic setting.}

\bigskip

\textbf{\textsf{3{\boldmath$.$}2{\boldmath$.$}1  Swap operator $\sf S${\boldmath$.$}}} The result stated below in Theorem \ref{ThmCommutatorParaproducts} is fully proved in Appendix \ref{Subsectionapp2} -- see Proposition \ref{prop:R1App} and Proposition \ref{prop:R2}. Given H\"older distributions $f,g_1,g_2,g_3,g,h$, we define the \textsf{\textbf{modified commutator on paraproducts}}, and its iterates, by the formulas
$$
{\sf S}(f,g,h) := {\sf P}_f\left( \widetilde {\sf P}_gh\right) - {\sf P}_g\big({\sf P}_f h\big),
$$
and 
\begin{equation}
\label{Eq4LinearSwap}
{\sf S}\big(f, (g_1, g_2),h\big) := {\sf S}\Big(f, \widetilde{\sf P}_{g_1}g_2,h\Big) - {\sf P}_{g_1}\Big({\sf S}(f, g_2,h)\Big)
\end{equation}
and 
\begin{equation}
\label{Eq5LinearSwap}
{\sf S}\Big(f, \big((g_1,g_2),g_3\big),h\Big) := {\sf S}\Big(f, \big(\widetilde{\sf P}_{g_1}g_2,g_3\big), h\Big) - {\sf P}_{g_1}\Big({\sf S}(f, (g_2,g_3), h)\Big).
\end{equation}
The continuity properties of these operators are given in the following statement. The range $(-3,3)$ in the next statement is purely technical and can be replaced by any other interval by an adequate choice of constant $b$ in the definition of the paraproduct ${\sf P} = {\sf P}^{(b)}$. Note that no gain of regularity comes from the first argument in the regularity statements below.

\medskip

\begin{thm}  {\it 
\label{ThmCommutatorParaproducts}
Pick $\gamma\in\bbR$.
\begin{itemize}
   \item[\textbf{\textsf{(a)}}] Pick $\alpha\in (0,1)$ and $\beta\in\bbR$ such that $\gamma+\alpha+\beta\in (-3,3)$. The map
\begin{equation*}
\begin{split}
\mcC^{\gamma}\times\mcC^{\alpha}\times\mcC^\beta &\rightarrow\mcC^{(\gamma\wedge 0)+\alpha+\beta}   \\
(f,g,h) &\mapsto {\sf S}(f,g,h)
\end{split}
\end{equation*}
   is continuous.   \vspace{0.15cm}
   
   \item[\textbf{\textsf{(b)}}] Pick $\alpha_1,\alpha_2\in (0,1)$ and $\beta\in\bbR$ such that $\gamma+\alpha_1+\alpha_2+\beta\in (-3,3)$. The map
\begin{equation*}
\begin{split}
\mcC^{\gamma}\times\mcC^{\alpha_1}\times\mcC^{\alpha_2}\times\mcC^\beta &\rightarrow\mcC^{(\gamma\wedge 0)+\alpha_1+\alpha_2+\beta}   \\
\big(f, g_1,g_2,h\big) &\mapsto {\sf S}\big(f, (g_1,g_2),h\big)
\end{split}
\end{equation*}
   is continuous.   \vspace{0.15cm}
   
   \item[\textbf{\textsf{(c)}}] Pick $\alpha_1,\alpha_2, \alpha_3 \in (0,1)$ and $\beta\in\bbR$ such that $\gamma+\alpha_1+\alpha_2+\alpha_3+\beta\in (-3,3)$. The map
\begin{equation*}
\begin{split}
\mcC^{\gamma}\times\mcC^{\alpha_1}\times\mcC^{\alpha_2}\times\mcC^{\alpha_3}\times\mcC^\beta &\rightarrow\mcC^{(\gamma\wedge 0)+\alpha_1+\alpha_2+\alpha_3+\beta}   \\
\big(f, g_1,g_2, g_3,h\big) &\mapsto {\sf S}\big(f, \big((g_1,g_2), g_3\big),h\big)
\end{split}
\end{equation*}
   is continuous. 
\end{itemize}   }
\end{thm}

\bigskip

\textbf{\textsf{3{\boldmath$.$}2{\boldmath$.$}2  Merging operator $\sf R$ and inner difference operator $\mathscr{D}${\boldmath$.$}}} The value at $x\in {\bfT}^d$ of some paraproduct ${\sf P}_u v$ is a sum over the integers $i$ of terms of the form
$$
\Big(\Pi^{0,(i)}_u v\Big)(x) := \iint K_i(x,y) K_{\leq i-1}(x,z) \, u(z) v(y)\,dzdy.
$$
We thus have for $f\in L^\infty, g\in\mcC^\alpha$ with $\alpha\in(0,1)$, and $h\in \mcC^\nu$,
\begin{equation}
\label{EqDefnInnerDifference}
\begin{split}
\Big(\Pi^{0,(i)}_f\big(\Pi^0_gh\big) - \Pi^{0,(i)}_{fg}h\Big)(x) &= \iint K_i(x,y) K_{\leq i-1}(x,z) \,f(z) \Big(\Pi^0_{g-g(z)} h\Big)(y)\,dzdy   \\
&=: \iint K_i(x,y) K_{\leq i-1}(x,z) \,f(z) \big(\Pi^0_{\mathscr{D}g}h\big)(y)\,dzdy.
\end{split}
\end{equation}
The above identity defines the \textsf{\textbf{inner difference operator}} $\mathscr{D}\,\big(=\mathscr{D}_z\big)$. In those terms, we have
$$
{\sf R}^0(f,g,h) := \Pi^0_f\big(\Pi^0_a g\big) - \Pi^0_{fa}g = \Pi^0_f\Big(\Pi^0_{\mathscr{D}a} g\Big)
$$
and, given the definition of the inner difference operator in the parabolic setting of Section \ref{SectionAppendixContinuity}, we have more generally,
\begin{equation}
\label{eq:*}
{\sf R}(f,g,h) := {\sf P}_f\Big(\widetilde{\sf P}_g h \Big) - {\sf P}_{fg} h = {\sf P}_f\Big(\widetilde{\sf P}_{\mathscr{D}g} h \Big).
\end{equation}
We use the same letter $\mathscr{D}$ in the spatial and parabolic settings. Compare this expression with the formal multiple integral, where we use the same letters to make it more stricking,
$$
\int f(z)d\left(\int^\cdot gdh\right)(z) = \int fg dh + \int f(z)d\left(\int^\cdot \big(g-g(z)\big)dh\right)(z).
$$
A similar reasoning as in the proof of continuity of ${\sf C}(f,g,h) = {\sf \Pi}\big({\sf P}_{\mathscr{C}f}g,h\big)$, shows that $\sf R$ sends continuously $L^\infty\times\mcC^\alpha\times\mcC^\nu$ into $\mcC^{\alpha+\nu}$, as soon as $\alpha\in(0,1)$. A subtle thing happens here, though, as one has actually a refined continuity estimate on $\sf R^\circ$ that also takes into account the regularity of $f$ as well. It is related to the fact that ${\sf R}^\circ({\bf 1},g,h) = {\sf R}^\circ(f,{\bf 1},h) = 0$, for all $f,g,h$, whereas ${\sf R}({\bf 1},g,h)\neq 0$.

\medskip

\begin{prop}   \label{PropOuf}
\begin{itemize}
   \item For $\alpha, \beta\in [0,1)$ and $\gamma\in(-3,3)$, we have 
\begin{equation}
\label{eq:R1tera}
\big\| {\sf R}^\circ(f,g,h) \big\|_{\mcC^{\alpha+\beta+\gamma}} \lesssim \|f\|_{\mcC^\alpha} \, \|g\|_{\mcC^\beta} \, \|h\|_{\mcC^\gamma}.
\end{equation}   
   
   \item For $a_1,a_2,h\in\mcC^\alpha$, one has
   $$
   {\sf R}\big(1,{\sf P}_{a_1}a_2,h\big) - {\sf P}_{a_1}{\sf R}(1,a_2,h) \in\mcC^{3\alpha}.   
   $$  

   \item If $a_2,g\in \mcC^\alpha$, with $0<\alpha<1/2$, and $a_1\in \mcC^\nu$, with $\nu\in(0,1/2)$, and $h\in \mcC^\gamma$ for $\gamma\in(-3,3)$, then we have
\begin{equation} 
\label{prop:ouf2}
{\sf R}^\circ \big({\sf P}_{a_1}a_2,g,h\big) - {\sf P}_{a_1}{\sf R}^\circ(a_2,g,h) \in \mcC^{2\alpha+\gamma+\nu}. 
\end{equation}   
\end{itemize}
\end{prop}   

\medskip

Note that we need both $f$ and $g$ to be in $\mcC^\alpha$ in \eqref{eq:R1tera}; this is required by the method of proof, based on an interpolation argument. 

\medskip

\begin{Dem} The first statement \eqref{eq:R1tera} is proved in details in the third statement, inequality \eqref{eq:R1ter}, of Proposition \ref{prop:R1App} in Appendix \ref{SectionAppendixContinuity}. The second statement is a particular ($u=1$) case of Proposition \ref{prop:R}, also in the same appendix. Let us check the last statement \eqref{prop:ouf2}. Because of the symmetry character that the first two arguments of ${\sf R}^\circ$ play in the proof of its regularity properties, Proposition \ref{prop:R} implies the boundedness of the operator 
$$
\big(a_1,a_2,g,h\big) \rightarrow {\sf R}^\circ \big({\sf P}_{a_1}a_2,g,h\big) - {\sf P}_{a_1}{\sf R}^\circ(a_2,g,h)
$$ 
from $\mcC^{2\nu} \times L^\infty \times \mcC^{2\alpha} \times \mcC^\gamma$ into $\mcC^{2\nu+2\alpha+\gamma}$, because $2\nu,2\alpha<1$. We also have the boundedness of that operator from $L^\infty \times \mcC^{2\alpha} \times L^\infty \times \mcC^\gamma$ into $\mcC^{2\alpha+\gamma}$, because $2\alpha<1$ -- we only use \eqref{eq:R1tera} to estimate the two quantities and do not use the difference structure in the commutator. We conclude by interpolating between these two estimates.
\end{Dem}

\medskip

\begin{thm}
\label{ThmComposingParaproducts}   {\it 
\begin{itemize}
   \item Let $f,g\in L^\infty$ and $h\in\mcC^\nu$ be given for $\nu\in(-3,3)$. Let also $a_1\in\mcC^{\alpha_1}$ and $a_2\in\mcC^{\alpha_2}$ be given, with $\alpha_1, \alpha_2\in(0,1)$. Then both
\begin{equation}
\label{EqRHigh}
{\sf R}^\circ\big(f,(a_1,a_2),h\big) := {\sf R}^\circ\Big(f,{\sf P}_{a_1}a_2,h\Big) - {\sf R}^\circ(fa_1,a_2,h)
\end{equation}
and 
\begin{equation}
\label{EqRLow}
{\sf R}^\circ\big((a_1,a_2), g,h\big) := {\sf R}^\circ\Big(\widetilde{\sf P}_{a_1}a_2,g,h\Big) - {\sf P}_{a_1}{\sf R}^\circ(a_2,g,h)
\end{equation}
are elements of $\mcC^{\alpha_1+\alpha_2+\nu}$.   \vspace{0.1cm}

   \item If $f\in \mcC^{\mu}$, with $\mu\in(0,1)$, and $a_1\in\mcC^{\alpha_1}$ and $a_2\in\mcC^{\alpha_2}$, and $h\in\mcC^\nu$ with $\nu\in(-3,3)$, then we have
\begin{equation} 
\label{prop:ouf}
{\sf R}^\circ\big(f, (a_1,a_2),h\big) - {\sf P}_f\Big({\sf R}^\circ\big({\bf 1}, (a_1,a_2),h\big)\Big) \in \mcC^{\mu+\alpha_1+\alpha_2+\nu}. 
\end{equation}   
\end{itemize}   }
\end{thm}

\medskip

The range $(0,1)$ for the exponent $\alpha,\beta$ and $\gamma$, is dictated by the operator ${\mathscr{D}}$, which involves a first order increment and so can only encode regularity at order at most $1$. 

\medskip

\begin{Dem} The second estimate \eqref{EqRLow} is proved in the Appendix, see Proposition \ref{prop:R}, because the first two arguments of ${\sf R}^\circ$ play a ``symmetric'' role in the proof of the continuity estimates; see Remark \ref{rem:R}.

\ssk

We prove the two other corresponding statements \eqref{EqRHigh} and \eqref{prop:ouf} in the model time-independent setting of the flat torus. We prove first \eqref{EqRHigh}. We have
$$
\Pi^0_f\Big(\Pi^0_{\Pi^0_{a_1}a_2} g\Big) - \Pi^0_{f\Pi^0_{a_1}a_2} g - \Pi^0_{fa_1}\Big(\Pi^0_{\mathscr{D}a_2} g\big) = \Pi^0_f\Big(\Pi^0_{\mathscr{D}\Pi^0_{a_1}a_2} g\Big) - \Pi^0_{fa_1}\Big(\Pi^0_{\mathscr{D}a_2} g\Big) 
$$
is a sum over $i$ of double integrals
\begin{equation*}
\begin{split}
\iint K_i(x,y) K_{\leq i-2}(x,z)& \,f(z) \big(\Pi^0_{\mathscr{D}(\Pi^0_{a_1}a_2 -a_1(z)a_2)} g\big)(y)\,dzdy   \\
											   &= \iint K_i(x,y) K_{\leq i-2}(x,z) \,f(z) \big(\Pi^0_{\mathscr{D}\Pi^0_{\mathscr{D}a_1} a_2 } g\big)(y)\,dzdy
\end{split}
\end{equation*}
on which we read off that their $L^\infty$ norm is of order $2^{-i(\alpha_1+\alpha_2+\nu)}$. This point finishes the proof since the last quantity corresponds to the dyadic blocks $\Delta_i\big({\sf R}(f,(a_1,a_2),g)\big)$.

\smallskip

We prove \eqref{prop:ouf}. On the other hand, we have that
\begin{equation*}
\begin{split}
\Delta_i\Big({\sf P}_f\big({\sf R}\big({\bf 1}\,;(a_1,a_1)&\,;g\big)\big)\Big)(x)   \\
&= \iiint K_i(x,y) K_{\leq i-2}(x,u) \,f(u) K_{\leq i-2}(x,z) \big(\Pi^0_{\mathscr{D}\Pi^0_{\mathscr{D}a_1} a_2 } g\big)(y)\,dzdy du.
\end{split}
\end{equation*}
So using that the kernel $K_{\leq i-1}(x,\cdot)$ has an integral equal to $1$, we deduce that the difference of the two dyadic blocks is given by
\begin{equation*}
\iiint K_i(x,y) K_{\leq i-2}(x,z) K_{\leq i-2}(x,u) \, \big(f(z)-f(u)\big) \big(\Pi^0_{\mathscr{D}\Pi^0_{\mathscr{D}a_1} a_2 } g\big)(y)\,dzdydu,
\end{equation*}
on which we read off that their $L^\infty$ norm is of order $2^{-i(\mu+\alpha_1+\alpha_2+\nu)}$. 
\end{Dem}

\medskip

\begin{rem} \label{rem:R} {\it The operator ${\sf R}^\circ$ is not symmetric in a proper sense, but rather in terms of how the regularity properties of its first two arguments are taken into account. The frequency analysis is exactly the same, if we assume a regularity on $f$ or on $g$. Since ${\sf R}^\circ$ involves only ${\sf P}$-paraproducts and not $\widetilde{\sf P}$-paraproducts, we can exploit the difference structure between ${\sf P}_{fg}$  and ${\sf P}_f$ or ${\sf P}_g$. In the operator ${\sf R}$, one can only exploit the regularity on $g$ in ${\sf R}(f,g,h)$ because of the $\widetilde {\sf P}$-operator on $g$. This phenomenon is illustrated by the observation that 
$$
{\sf R}^\circ(1,g,h)={\sf R}^\circ(f,1,h)=0={\sf R}(f,1,h),
$$ 
while ${\sf R}(1,g,h)=\widetilde{\sf P}_gh-{\sf P}_gh\neq 0$.   }
\end{rem}

\bigskip

\textbf{\textsf{3{\boldmath$.$}2{\boldmath$.$}3  Back again to the high order paracontrolled expansion \eqref{EqTaylorExpansion}{\boldmath$.$}}} Recall from Proposition \ref{PropOuf} that the operator ${\sf R}^\circ(a,b,c) = {\sf P}_a{\sf P}_bc-{\sf P}_{ab}c$, has better regularity properties than the merging operator $\sf R$. We can use the continuity properties of ${\sf R}^\circ$ to rewrite the high order paracontrolled expansion \eqref{EqTaylorExpansion} and clarify the cancellations that happen in each of its brackets. We have first that
\begin{equation*}
\begin{split}
{\sf P}_{f^{(2)}(u)}u^2 - 2{\sf P}_{f^{(2)}(u)u}u &= {\sf P}_{f^{(2)}(u)}{\sf \Pi}(u,u) + 2{\sf P}_{f^{(2)}(u)}({\sf P}_uu) - 2{\sf P}_{f^{(2)}(u)u}u   \\
&= {\sf P}_{f^{(2)}(u)}{\sf \Pi}(u,u) + 2{\sf R}^\circ\big(f^{(2)}(u),u,u\big)   \\
&= {\sf P}_{f^{(2)}(u)}{\sf \Pi}(u,u) + 2{\sf P}_{f^{(3)}(u)}{\sf R}^\circ\big(u,u,u\big) + (2\alpha\wedge 1 + 2\alpha),
\end{split}
\end{equation*}
from Bony's paralinearisation $f^{(2)}(u) = {\sf P}_{f^{(3)}(u)}u + (2\alpha)$ and the continuity properties of the merging operator $\sf R$. Write ${\sf D}^\circ(a,b,c)$ for ${\sf \Pi}({\sf P}_ab,c)-{\sf P}_a{\sf \Pi}(b,c)$; this operator has the same regularity properties as $\sf D$. The third order expansion formula is only needed when $\alpha<1/2$. In that case, an elementary computation shows that 
\begin{equation*}
\begin{split}
&{\sf P}_{f^{(3)}(u)}u^3 - 3{\sf P}_{uf^{(3)}(u)}u^2 + 3{\sf P}_{u^2f^{(3)}(u)}u   \\
&= \Big(3{\sf R}^\circ\big(f^{(3)}(u),u^2,u\big) + 2{\sf P}_{f^{(3)}(u)}{\sf R}^\circ(u,u,u) + {\sf P}_{f^{(3)}(u)}\big(2{\sf D}^\circ(u,u,u) + {\sf \Pi}(u,{\sf \Pi}(u,u))\big)    \\
&\quad- 6{\sf R}^\circ\big(uf^{(3)}(u),u,u\big)\Big) + 3{\sf R}^\circ\big(f^{(3)}(u),u,{\sf \Pi}(u,u)\big)   \\
&= {\sf P}_{f^{(3)}(u)}\Big(2{\sf R}^\circ(u,u,u)+2{\sf D}^\circ(u,u,u) + {\sf \Pi}(u,{\sf \Pi}(u,u))\Big)   \\
&\quad+ \Big\{3{\sf R}^\circ\big(f^{(3)}(u),{\sf \Pi}(u,u),u\big) + 6{\sf R}^\circ\big(f^{(3)}(u),(u,u),u\big) + 3{\sf R}^\circ\big(f^{(3)}(u),u,{\sf \Pi}(u,u)\big)\Big\}
\end{split}
\end{equation*}
is the sum of $3\alpha$-terms and a $4\alpha$-term, as one can use in that case the refined continuity estimates \eqref{prop:ouf2} on $\sf R$. 

\medskip

\begin{cor}
\label{PropTaylorPractical}
Let $f\in C^5_b$. If $\alpha<1/2$, one has
\begin{equation}
\label{EqTaylorPractical}
\begin{split}
f(u) &= {\sf P}_{f'(u)}u + \frac{1}{2!}\,{\sf P}_{f^{(2)}(u)}{\sf \Pi}(u,u)   \\
	   &\quad+ \frac{1}{3!}\,{\sf P}_{f^{(3)}(u)}\Big(8\,{\sf R}^\circ(u,u,u)+2\,{\sf D}^\circ(u,u,u) + {\sf \Pi}(u,{\sf \Pi}(u,u))\Big) + f(u)^\sharp,
\end{split}
\end{equation}
for a remainder $f(u)^\sharp\in\mcC^{4\alpha}$.
\end{cor}

\medskip

\begin{prop}
\label{CorDecompositionFuZeta}
Let a noise $\zeta\in\mcC^{\alpha-2}$ be given, together with a function $u$ given by a paracontrolled system with reference functions $\bfZ$, and $f\in C^5_b$. Then one can write the product $f(u)\zeta$ under the form
\begin{equation}
\label{EqDecompositionFuZeta}
\begin{split}
f(u)\zeta &= {\sf P}_{f(u)}\zeta + \sum_{i=1}^2{\sf P}_{f'(u)u_i}Y_i + {\sf P}_{f'(u)u_{11}}Y_{11} + {\sf P}_{f^{(2)}(u)u_1^2}Y_{(1,1)} + (\sharp)
\end{split}
\end{equation}
for explicit reference distributions $Y$'s that depend only on the noise $\zeta$ and $\bfZ$, and a remainder term 
$$
(\sharp) \in \mathscr{L}(\mcC^{4\alpha}) \subset \mcC^{4\alpha-2},
$$
that depends continuously on $\widehat{u}$ and the $Y$'s.
\end{prop}

\medskip

The proper interpretation of the above statement is that the remainder $(\sharp) \in \mathscr{L}(\mcC^{4\alpha})$, provided the noise $\zeta$ is enhanced into $\widehat{\zeta}$, with components belonging to $\mathscr{L}(\mcC^{k\alpha-2})$-spaces; the remainder $(\sharp)$ also depends on $\widehat{\zeta}$; see Remark 1 after the proof. We first give a detailled proof of the statement, in which the reader will see that we keep repeating the same kind of computations. The mechanics at play here will be emphasized after the proof; this is nothing but the mechanics of Section \ref{SubsectionMechanicsComputations}.

\medskip

\begin{Dem}
We start from the identity
$$
f(u)\zeta = {\sf P}_{f(u)}\zeta + {\sf P}_{\zeta}f(u) + {\sf \Pi}\big(f(u),\zeta\big),
$$
and treat the second and third terms separately. We shall denote by $(\flat)$ a term in $\mathscr{L}(\mcC^{4\alpha})$ that may change from place to place. A term of the form ${\sf P}_\zeta b$, with $b\in\mcC^{3\alpha}$, is a $(\flat)$-term; so is a term ${\sf P}_a (\flat)$, if $a$ has non-negative parabolic H\"older regularity. We use below the following fact proved in Proposition \ref{prop:R1App} - \eqref{eq:R1} and Proposition \ref{PropR2} of Appendix \ref{SectionAppendixContinuity}.

\medskip

\begin{lem}\label{LemRemainder}
Let $a_1\in\mcC^{\alpha_1}, a_2\in\mcC^{\alpha_2}, a_3\in\mcC^{\alpha_3}$ be given, with $\alpha_i\in(0,1)$. The quantities 
$$
{\sf S}^\circ(\zeta,a_1,a_2), \;{\sf S}(\zeta,a_1,a_2), \;{\sf S}\big(\zeta, (a_1,a_2),a_3\big),
$$ 
belong to $\mathscr{L}(\mcC^{4\alpha})$, if $\alpha_1+\alpha_2\geq 3\alpha$ for the first two quantities, and $\alpha_1+\alpha_2+\alpha_3\geq 3\alpha$, for the third quantity. These operators are continuous functions of the $a_i$ under these assumptions.
\end{lem}

\medskip

\textsf{\textbf{1.}} We have
\begin{equation*}\begin{split}
{\sf P}_\zeta f(u) &=  {\sf P}_\zeta{\sf P}_{f'(u)}u + \frac{1}{2}\,{\sf P}_\zeta{\sf P}_{f^{(2)}(u)}{\sf \Pi}(u,u) + (\flat)   \\
					   &= {\sf P}_{f'(u)}{\sf P}_\zeta u + {\sf S}^\circ\big(\zeta,f'(u),u\big) + \frac{1}{2}\,{\sf P}_{f^{(2)}(u)}{\sf P}_\zeta{\sf \Pi}(u,u) + (\flat),
\end{split}\end{equation*}
from Lemma \ref{LemRemainder}. So ${\sf P}_\zeta f(u)$ is equal to 
\begin{equation*}\begin{split}
&{\sf P}_{f'(u)}{\sf P}_{u_i}({\sf P}_\zeta Z_i) + {\sf P}_{f'(u)}{\sf S}(\zeta,u_i,Z_i) + {\sf S}^\circ\big(\zeta,f'(u),u\big) + \frac{1}{2}\,{\sf P}_{f^{(2)}(u)}{\sf P}_\zeta{\sf P}_{u_1^2}{\sf \Pi}(Z_1,Z_1) + (\flat)   \\
&= \Big({\sf P}_{f'(u)u_i}{\sf P}_\zeta Z_i + {\sf R}^\circ\big(f'(u),u_i,{\sf P}_\zeta Z_i\big)\Big) + \Big({\sf P}_{f'(u)}{\sf P}_{u_{ij}} {\sf S}(\zeta,Z_j,Z_i) + (\flat)\Big)   \\
&\quad + \Big({\sf P}_{f^{(2)}(u)u_1^2}{\sf S}^\circ(\zeta,Z_1,Z_1) + (\flat)\Big) + \left(\frac{1}{2}\,{\sf P}_ {f^{(2)}(u)}{\sf P}_{u_1^2}\big({\sf P}_\zeta {\sf \Pi}(Z_1,Z_1)\big) + (\flat)\right) + (\flat),
\end{split}\end{equation*}
from Lemma \ref{LemRemainder} again. Note that it follows from the sharp continuity estimate \eqref{PropOuf} from Proposition \ref{PropOuf} that ${\sf R}^\circ\big(f'(u),u_i,{\sf P}_\zeta Z_i\big)$ is a $(\flat)$-term. Also, one has from Lemma \ref{LemRemainder} that ${\sf S}(\zeta, Z_j,Z_i)\in\mathscr{L}(\mcC^{(i+j+1)\alpha})$ and ${\sf S}^\circ(\zeta, Z_j,Z_i)\in\mathscr{L}(\mcC^{(i+j+1)\alpha})$, while ${\sf P}_\zeta {\sf \Pi}(Z_1,Z_1)\in\mathscr{L}(\mcC^{3\alpha})$. This allows us to write ${\sf P}_\zeta f(u)$ under the required form
\begin{equation*}\begin{split}
{\sf P}_\zeta f(u) &= {\sf P}_{f'(u)u_i}({\sf P}_\zeta Z_i) + {\sf P}_{f'(u)u_{11}} {\sf S}(\zeta,Z_1,Z_1) + {\sf P}_{f^{(2)}(u)u_1^2}{\sf S}^\circ(\zeta,Z_1,Z_1)   \\
&\quad+  \frac{1}{2}\,{\sf P}_ {f^{(2)}(u)u_1^2}\big({\sf P}_\zeta {\sf \Pi}(Z_1,Z_1)\big) + (\flat).
\end{split}\end{equation*}

\medskip

\textsf{\textbf{2.}} Consider now the resonant term ${\sf \Pi}\big(f(u),\zeta\big)$, and start for that purpose from the multiplicative formula \eqref{EqDecompositionPrimitiveResonant} 
\begin{equation*}\begin{split}
&f'(u)\,{\sf \Pi}(u,\zeta) + f^{(2)}(u)\,\left\{{\sf C}^\circ(u,u,\zeta) + \frac{1}{2}\,{\sf \Pi}\big({\sf \Pi}(u,u),\zeta\big)\right\}   \\
&+ f^{(3)}(u)\,\left\{\frac{1}{3}\,{\sf C}^\circ\big(u,(u,u),\zeta\big) + \frac{1}{3}\,{\sf C}^\circ\big((u,u),u,\zeta\big) + \frac{1}{6}\,{\sf C}^\circ\big(u,{\sf \Pi}(u,u),\zeta\big) + \frac{1}{6}\,{\sf C}^\circ\big({\sf \Pi}(u,u),u,\zeta\big) \right.  \\
&\hspace{1.8cm}+ \left. \frac{1}{3}\,{\sf \Pi}\big({\sf D}^\circ(u,u,u),\zeta\big) + \frac{1}{6}\,{\sf \Pi}\big({\sf \Pi}(u,{\sf \Pi}(u,u)),\zeta\big) \right\}   \\
&+ \sum_{\ell=1}^{\ell_0}(\star)_\ell\,{\sf \Pi}_{(1)}^\ell(u,\zeta) + (5\alpha-2)
\end{split}\end{equation*}
for ${\sf \Pi}\big(f(u),\zeta\big)$ -- recall $(\star)_\ell$ is defined in \eqref{EqStarEll}. It suffices to see that each term in this decomposition has the form \eqref{EqDecompositionFuZeta} of the statement; we proceed from the more to the less regular terms. 

\medskip

\textsf{\textbf{2.1.}} Use Theorem \ref{ThmPi1} on the expansion rule for the operator ${\sf \Pi}_{(1)}$, and the first order paracontrolled structure of $u$, to write
$$
{\sf \Pi}_{(1)}^\ell(u,\zeta) = u_1{\sf \Pi}^\ell_{(1)}(Z_1,\zeta) + (3\alpha-1).
$$
This gives the paracontrolled representation
\begin{equation}   \label{EqContributionStarEll}
(\star)_\ell\,{\sf \Pi}_{(1)}^\ell(u,\zeta) = {\sf P}_{(\star)_\ell u_1}{\sf \Pi}_{(1)}^\ell(Z_1,\zeta) + (5\alpha-2),
\end{equation}
since $\alpha<1/2$.   

\medskip

\textsf{\textbf{2.2.}} Write $\{u\}$ for the bracket term after $f^{(3)}(u)$. The first order paracontrolled structure of $u$ and the continuity properties of the correctors and commutators give 
$$
f^{(3)}(u)\{u\} = {\sf P}_{f^{(3)}(u)u_1^3}\{Z_1\} + (5\alpha-2),
$$
for a remainder term $(5\alpha-2)\in\mcC^{5\alpha-2}$. The fact that $\{Z_1\}\in\mathscr{L}(\mcC^{4\alpha})$ is part of the data $\widehat{\zeta}$. 

\medskip

\textsf{\textbf{2.3.}} We use the paracontrolled structure of $u$ and the continuity results on correctors to write the term ${\sf \Pi}\big({\sf \Pi}(u,u),\zeta\big)$ in multiplicative form. This gives
\begin{equation*}\begin{split}
{\sf \Pi}\big({\sf \Pi}(u,u),\zeta\big) &= {\sf \Pi}\big({\sf P}_{u_i}{\sf \Pi}(Z_i,u),\zeta\big) + {\sf \Pi}\big({\sf D}(u_i,Z_i,u),\zeta\big) + (5\alpha-2)   \\
										   &= u_i{\sf \Pi}\big({\sf \Pi}(Z_i,u),\zeta\big) + {\sf C}^\circ\big(u_i,{\sf \Pi}(Z_i,u),\zeta\big) + u_{11}u_1{\sf \Pi}\Big({\sf D}(Z_1,Z_1,Z_1),\zeta\Big) + (5\alpha-2)  \\
										   &= u_iu_j{\sf \Pi}\big({\sf \Pi}(Z_i,Z_j),\zeta\big) + u_{11}u_1\,(\star)+ (5\alpha-2)  \\
\end{split}\end{equation*}
with
$$
(\star) := 2{\sf C}^\circ\big(Z_1,{\sf \Pi}(Z_1,Z_1),\zeta\big) + 2{\sf \Pi}\Big({\sf D}(Z_1,Z_1,Z_1),\zeta\Big),
$$
after elementary computations. The fact that ${\sf \Pi}\big({\sf \Pi}(Z_1,Z_2),\zeta\big)$ and $(\star)$ are elements of $\mathscr{L}(\mcC^{4\alpha})$ is part of the data $\widehat{\zeta}$. Now, 
$$
f^{(2)}(u)u_1^2 = {\sf P}_{2f^{(2)}(u)u_1u_{11}+f^{(3)}(u)u_1^3}Z_1+(2\alpha),
$$
so, recalling that an element in $\mcC^{5\alpha-2}$ is of $(\flat)$-type, we have
\begin{equation*}\begin{split}
&f^{(2)}(u)u_1^2 \,{\sf \Pi}\big({\sf \Pi}(Z_1,Z_1),\zeta\big)   \\
&= {\sf P}_{f^{(2)}(u)u_1^2}{\sf \Pi}\big({\sf \Pi}(Z_1,Z_1),\zeta\big) + (\flat) + {\sf \Pi}\Big(f^{(2)}(u)u_1^2, {\sf \Pi}\big({\sf \Pi}(Z_1,Z_1),\zeta\big)\Big)   \\
&= {\sf P}_{f^{(2)}(u)u_1^2}{\sf \Pi}\big({\sf \Pi}(Z_1,Z_1),\zeta\big) + {\sf P}_{2f^{(2)}(u)u_1u_{11}+f^{(3)}(u)u_1^3}{\sf \Pi}\Big(Z_1, {\sf \Pi}\big({\sf \Pi}(Z_1,Z_1),\zeta\big)\Big) + (\flat).
\end{split}\end{equation*}
Here again, the fact that ${\sf \Pi}\big(Z_1, {\sf \Pi}\big({\sf \Pi}(Z_1,Z_1),\zeta\big)\big)$ is of $(\flat)$-type is part of the data $\widehat{\zeta}$, so the term $f^{(2)}(u)u_1^2 {\sf \Pi}\big({\sf \Pi}(Z_1,Z_1),\zeta\big)$ has indeed the right form. Similar computations show that $f^{(2)}(u)u_1u_2 {\sf \Pi}\big({\sf \Pi}(Z_1,Z_2),\zeta\big)$ and $u_{11}u_1(\star)$ also have the right form.

\smallskip

Very similar computations give the right decomposition of the term $f^{(2)}(u)\,{\sf C}^\circ(u,u,\zeta)$, once ${\sf C}^\circ(u,u,\zeta)$ has been put in multiplicative form. We start from the identity
$$
{\sf C}^\circ(u,u,\zeta) = {\sf C}^\circ\big({\sf P}_{u_i}Z_i,u,\zeta\big) + {\sf C}^\circ\Big({\sf R}(1,u_1,Z_1),u,\zeta\Big) + {\sf C}^\circ\big(u^{(1)},u,\zeta\big) + (5\alpha-2)
$$
with
$$
u^{(1)} := u-\widetilde{\sf P}_{u_1}Z_1 - \widetilde{\sf P}_{u_2}Z_2
$$
an element of $\mcC^{3\alpha}\subset \mcC^1$, since $3\alpha>1$. Recall Definition \ref{DefnRefinedCorrector} of the refined corrector, and observe that because of Theorem \ref{ThmUpper}
$$ 
{\sf C}^\circ\big(u^{(1)},u,\zeta\big)= u_1{\sf C}^\circ\big(u^{(1)},Z_1,\zeta\big) + (5\alpha-2),
$$ 
since $3\alpha-1>5\alpha-2$. Then using the continuity result on the refined corrector ${\sf C}_{(1)}$ from Theorem \ref{ThmCorrectors1}, and the fact that $\alpha>5\alpha-2$, one has
\begin{equation}   \label{EqRSPC1}
{\sf C}^\circ\big(u^{(1)},Z_1,\zeta\big) = \sum_{\ell=1}^{\ell_0} \gamma_\ell V_\ell\big(u^{(1)}\big)\, {\sf \Pi}_{(1)}^\ell(Z_1,\zeta) + (5\alpha-2).
\end{equation} 
The fact that ${\sf \Pi}_{(1)}^\ell(Z_1,\zeta)\in\mathscr{L}(\mcC^{2\alpha+1})$ is part of the data $\widehat{\zeta}$, and it follows from the fact that $5\alpha-2>0$, that the preceding sum is a $(\flat)$-term. This gives, from Proposition \ref{PropOuf}, the identity
$$
{\sf C}^\circ(u,u,\zeta) = u_i{\sf C}^\circ(Z_i,u,\zeta) + {\sf C}^\circ\big((u_i,Z_i),u,\zeta\big) + u_{11}{\sf C}^\circ\big({\sf R}(1,Z_1,Z_1),u,\zeta\big) + (\flat) + (5\alpha-2).
$$
The multiplicative decomposition of ${\sf C}^\circ(u,u,\zeta)$ follows by an elementary computation.

\medskip

\textsf{\textbf{2.4.}} Last, we have the term $f'(u){\sf \Pi}(u,\zeta)$. We first put ${\sf \Pi}(u,\zeta)$ in multiplicative form
\begin{equation}   \label{EqRSPC2}
{\sf \Pi}(u,\zeta) = u_i{\sf \Pi}(Z_i,\zeta) + u_{ij}{\sf C}(Z_j,Z_i,\zeta) + u_{111}{\sf C}\big((Z_1,Z_1),Z_1,\zeta\big) + {\sf C}\big(u_1^\sharp,Z_1,\zeta\big) + (5\alpha-2).
\end{equation}
The term ${\sf C}\big(u_1^\sharp,Z_1,\zeta\big)$ is treated using the refined corrector, as above, which gives a contribution 
\begin{equation}   \label{EqRSPC3}
\sum_{\ell=1}^{\ell_0} f'(u)\gamma_\ell V_\ell(u_1^\sharp){\sf \Pi}^\ell_{(1)}(Z_1,\zeta) + (5\alpha-2),
\end{equation}
for that term, in the analysis of $f'(u){\sf \Pi}(u,\zeta)$; this is a remainder term of $(\flat)$-type. Look at the term $f'(u)u_1{\sf \Pi}(Z_1,\zeta)$; the other terms are dealt with similarly. One has
\begin{equation*}\begin{split}
f'(u)u_1{\sf \Pi}(Z_1,\zeta) &= {\sf P}_{f'(u)u_1}{\sf \Pi}(Z_1,\zeta) + {\sf P}_{{\sf \Pi}(Z_1,\zeta) }f'(u)u_1 + {\sf \Pi}\big(f'(u)u_1, {\sf \Pi}(Z_1,\zeta)\big)  \\
									 &= {\sf P}_{f'(u)u_1}{\sf \Pi}(Z_1,\zeta) + {\sf P}_{{\sf \Pi}(Z_1,\zeta) }\big({\sf P}_{f'(u)}u_1+{\sf P}_{u_1}f'(u)\big) +  {\sf P}_{{\sf \Pi}(Z_1,\zeta) }{\sf \Pi}\big(f'(u),u_1\big)   \\
									 &\quad+ f'(u){\sf \Pi}\big(u_1, {\sf \Pi}(Z_1,\zeta)\big) + {\sf C}^\circ\big(f'(u),u_1,{\sf \Pi}(Z_1,\zeta)\big)   \\
									 &\quad+ u_1{\sf \Pi}\big(f'(u), {\sf \Pi}(Z_1,\zeta)\big) + {\sf C}^\circ\big(u_1,f'(u),{\sf \Pi}(Z_1,\zeta)\big)   \\
									 &\quad+ {\sf \Pi}\big({\sf \Pi}(f'(u),u_1),{\sf \Pi}(Z_1,\zeta)\Big)
\end{split}\end{equation*}
The last six terms in the right hand side are of regularity $(4\alpha-2)$; a first order expansion of $f'(u)$ and $u_1$ allows to put the in the right form \eqref{EqDecompositionFuZeta}. For the second and third terms in the right hand side, simply write
$$
{\sf P}_{{\sf \Pi}(Z_1,\zeta)}{\sf P}_{f'(u)}u_1= {\sf P}_{f'(u)}{\sf P}_{{\sf \Pi}(Z_1,\zeta)}u_1 + {\sf S}^\circ\big({\sf \Pi}(Z_1,\zeta), f'(u),u_1\big).
$$
The last term has regularity $(4\alpha-2)$, and a first order expansion of $f'(u)$ and $u_1$ allows to put the in the right form. Also we have
\begin{equation*}\begin{split}
{\sf P}_{f'(u)}{\sf P}_{{\sf \Pi}(Z_1,\zeta)}u_1 &= {\sf P}_{f'(u)}{\sf P}_{{\sf \Pi}(Z_1,\zeta)}{\sf P}_{u_{1j}}Z_j + (5\alpha-2)   \\
														&= {\sf P}_{f'(u)}{\sf P}_{u_{1j}}{\sf P}_{{\sf \Pi}(Z_1,\zeta)}Z_j + {\sf P}_{f'(u)}{\sf S}^\circ\big({\sf \Pi}(Z_1,\zeta), u_{1j},Z_j\big)  + (5\alpha-2)   \\
														&= {\sf P}_{f'(u)u_{1j}}{\sf P}_{{\sf \Pi}(Z_1,\zeta)}Z_j + {\sf P}_{f'(u)u_{111}}{\sf S}^\circ\big({\sf \Pi}(Z_1,\zeta), Z_1,Z_1\big)  + (\flat)
\end{split}\end{equation*}
from Proposition \ref{PropOuf}. Once again, the fact that ${\sf S}^\circ\big({\sf \Pi}(Z_1,\zeta), Z_1,Z_1\big)$ is of $(\flat)$-type is part of the data on $\widehat{\zeta}$; this eventually gives
$$
{\sf P}_{f'(u)}{\sf P}_{{\sf \Pi}(Z_1,\zeta)}u_1 = {\sf P}_{f'(u)u_{1j}}{\sf P}_{{\sf \Pi}(Z_1,\zeta)}Z_1 + (\flat).
$$
The term ${\sf P}_{{\sf \Pi}(Z_1,\zeta) }{\sf P}_{u_1}f'(u)$ is treated similarly.
\end{Dem}

\bigskip

\noindent \textsf{\textbf{Remarks --}} \textsf{\textbf{1.}} {\sl The proof of Proposition \ref{CorDecompositionFuZeta} gives the structure
$$
(\sharp) = (\flat) + {\sf P}_{h_n(\widehat{u})}Y^n_{4\alpha-2} + (5\alpha-2)
$$
for the remainder $(\sharp)$ in \eqref{EqDecompositionFuZeta}, with an implicit sum in $n$ in the right hand side, explicit functions $h_n(\widehat{u})\in\mcC^\alpha$ depending continuously on $\widehat{u}$, and reference distributions $Y^n_{4\alpha-2}$ assumed to be in $\mathscr{L}(\mcC^{4\alpha})$, along the way; the latter are components of $\widehat{\zeta}$ -- see Section \ref{SubsectionEnhancedDistribution}.   }

\medskip

\textsf{\textbf{2.}} {\sl The above computations have a simple structure, that can be summarized using the ${\sf E}$ and ${\sf F}$ notations from Section \ref{SubsectionMechanicsComputations}. Denote by ${\sf F}^\beta(\cdot)$, or ${\sf F}^\beta(\cdot,\cdot)$, a well-defined function on parabolic H\"older spaces that sends $\mcC^{k\alpha}$ into $\mcC^{k\alpha + \beta}$, respectively $\mcC^{k\alpha}\times\mcC^{\ell\alpha}$ into $\mcC^{(k+\ell)\alpha + \beta}$, and enjoying the ${\sf F}$-type expansion property. We allow ourselves to write identities like 
\begin{equation}\label{EqFRule}
{\sf F}^\beta({\sf P}_{a_1}a_2) = {\sf P}_{a_1}{\sf F}^\beta(a_2) + {\sf F}^{\vert a_2\vert+\beta}(a_1),
\end{equation}
for $a_2\in\mcC^{\vert a_2\vert}$. Also, denote by ${\sf E}^\beta(\cdot)$, or ${\sf E}^\beta(\cdot,\cdot)$, an operator satisfying an ${\sf E}$-type expansion formula, sending formally $\mcC^\alpha$ into $\mcC^{\alpha+\beta}$, respectively $\mcC^{k\alpha}\times\mcC^{\ell\alpha}$ into $\mcC^{(k+\ell)\alpha+\beta}$. In those terms, we have 
\begin{equation}\label{EqERule}
{\sf E}^\beta({\sf P}_{a_1}a_2) = {\sf P}_{a_1}{\sf E}^\beta(a_2) + {\sf F}^{\beta+\vert a_2\vert}(a_1) + {\sf E}^{\beta+\vert a_2\vert}(a_1).
\end{equation}
The above proof of decomposition \eqref{EqDecompositionFuZeta} starts from the identity
$$
f(u)\zeta = {\sf P}_{f(u)}\zeta + {\sf F}^{\alpha-2}\big(f(u)\big) + {\sf E}^{\alpha-2}\big(f(u)\big),
$$
and proceeds by writing
\begin{equation*}\begin{split}
{\sf F}^{\alpha-2}\big(f(u)\big) &= {\sf F}^{\alpha-2}(u) + {\sf F}^{\alpha-2}\big(f'(u),u\big) + {\sf F}^{\alpha-2}\big({\sf \Pi}(u,u)\big) + (\flat)   \\
										   &= {\sf F}^{\alpha-2}(u) + {\sf F}^{\alpha-2}\big(u,u\big) + (\flat),
\end{split}\end{equation*}
and 
\begin{equation*}\begin{split}
{\sf E}^{\alpha-2}\big(f(u)\big) &= f'(u)\,{\sf E}^{\alpha-2}(u) + f^{(2)}(u)\,{\sf E}^{\alpha-2}(u,u) + f^{(3)}(u)\,{\sf E}^{\alpha-2}(u,u,u) + (\flat).
\end{split}\end{equation*}
One then uses the paracontrolled structure of $u$ and the expansion rules \eqref{EqFRule} and \eqref{EqERule} to run the computations.   }

\medskip

\textsf{\textbf{3.}} {\sl Note that Proposition \ref{CorDecompositionFuZeta} makes sense from a regularity structures point of view. Let us work in the regularity structure of the (gPAM) equation on the $3$-dimensional torus, together with a model on it. Let $u$ be represented by the modelled distribution
$$
\underline{u} = u{\bf 1} + u'X + u_\tau\underline{\tau} =: \sum_{a\in\mathscr{A}} u_a\,\underline{Z}_a ,
$$
 in $\mcD^\gamma$, for $\gamma=(3/2)^+$. Denote by $\circ$ the noise symbol in the regularity structure. Then one would have
$$
f(\underline{u})\circ = f(u)\circ + f'(u)u_a\,\underline{Z}_a\circ + \frac{1}{2!}\, f^{(2)}(u)u_au_b\,\underline{Z}_a\underline{Z}_b\circ + \frac{1}{3!}\, f^{(3)}(u)u_au_bu_c\,\underline{Z}_a\underline{Z}_b\underline{Z}_c\circ,
$$
with sums in $\mathscr{A}$ restricted to $\vert a\vert, \vert a\vert+\vert b\vert, \vert a\vert+\vert b\vert+\vert c\vert\leq 3\alpha$. Using Theorem 1 in Bailleul and Hoshino's work \cite{BailleulHoshino}, we would have for $f(u)\zeta$ the paracontrolled representation 
\begin{equation}\label{EqSyntheticFUZeta}\begin{split}
&{\sf P}_{f(u)}\zeta + {\sf P}_{f'(u)u_a}[\underline{Z}_a\circ] + \frac{1}{2!}\,{\sf P}_{f^{(2)}(u)u_au_b}[\underline{Z}_a\underline{Z}_b\circ] + \frac{1}{3!}\,{\sf P}_{f^{(3)}(u)u_au_bu_c}[\underline{Z}_a\underline{Z}_b\underline{Z}_c\circ] + (5\alpha-2)   \\
&=: {\sf P}_{f(u)}\zeta + {\sf P}_{f'(u)u_a}[\underline{Z}_a\circ] + \frac{1}{2!}\,{\sf P}_{f^{(2)}(u)u_au_b}[\underline{Z}_a\underline{Z}_b\circ] + (\flat),
\end{split}\end{equation}
for distributions $[\underline{Z}_a\circ], [\underline{Z}_a\underline{Z}_b\circ], [\underline{Z}_a\underline{Z}_b\underline{Z}_c\circ]$, built from the regularity structure and the model. This is the content of identity \eqref{EqDecompositionFuZeta}, and one reads off the functions $h_\ell(\widehat{u})$ of point 1 in the preceding formula. The term $f'(u)u' X\circ$ from the regularity structures picture appears in the above paracontrolled analysis under the form $f'(u)u_3{\sf \Pi}(Z_3,\zeta)$, in identity \eqref{EqRSPC2}. This makes perfect sense if one considers that the piece of $u$ that is differentiable is in its paracontrolled representation is given by $\overline{\sf \Pi}_{u_3}Z_3 + u^\sharp$, and we recall the pointwise first order expansion for paraproducts from Corollary \ref{CorEstimate}. The term $f'(u)^2u'\circ\mathcal{I}(X\circ)$ from the regularity structures picture appears above under the form
$$
f'(u)\,\gamma_\ell V_\ell(u_1^\sharp)\,{\sf \Pi}_{(1)}^\ell(Z_1,\zeta),
$$
in identity \eqref{EqRSPC3}, while the term $f^{(2)}(u)f(u)u'X\circ\,\mathcal{I}(\circ)$ from the regularity structures analysis appears 
in the form of identities \eqref{EqContributionStarEll} and \eqref{EqRSPC1}. (Recall $\sum_{\ell=1}^{\ell_0}\gamma_\ell V_\ell(v) = v'$, on the flat torus.)   }

\bigskip

\subsection{Dealing with derivatives}
\label{SubsectionPPdctDerivative}

We work in this section in the one dimensional torus $\bfT$, with $x$ as canonical coordinate and $L=\partial_x^2=:\partial^2$. The study of the generalised KPZ equation requires the analysis of quantities of the form ${\sf P}_{\partial {\sf P}_fg}\partial h$, or similar quantities where $f,g$ or $h$ is itself given by a paraproduct. The following remark provides the key to the analysis of such terms. 

Recall the notations $\mcP_t$ and $\mcQ_t$ for the parabolic projectors from the standard collection of operators with cancellation -- Definition \ref{DefnStandardCollection} in Appendix \ref{SubsectionAppendixApproxOperators}, and the notation $\nu = dt\otimes dx$ for the parabolic volume measure.

\medskip

\begin{lem}
\label{LmPiPartial}
Both ${\sf P}_{\partial f}\partial g$ and ${\sf \Pi}(\partial f,\partial g)$ can be written under the form 
$$ 
\int_0^1 \frac{1}{t}\,\mcP_t^{\bullet} \big(\widetilde \mcQ^1_tf \, \widetilde \mcQ^2_tg \big) \, \frac{dt}{t},
$$
with $\widetilde{\mcQ}^1,\widetilde{\mcQ}^2\in{\sf GC}^1$, and $\mcP_t\in{\sf StGC}^{[0,b]}$.
\end{lem}

\medskip

\begin{Dem}
Consider $\mcP_t$, one of our localization operator at the parabolic scale $t^{1/2}$, then by integrating by parts in space, we see that for $e:=(\tau,x)$
\begin{align*}
(\mcP_tf)(e) & = \int_{\mcM} K_{\mcP_t}(e,e') \partial f(e')  \nu(de') \\
             & = - \int_{\mcM} \partial_{x'} K_{\mcP_t}(e,e') f(e') \nu(de') \\ 
             & = - t^{-1/2} \int_{\mcM} t^{1/2} \partial_{x'} K_{\mcP_t}(e,e') f(e') \nu(de').
\end{align*}
Then since we assume regularity estimates on the heat kernel, it follows that $t^{1/2} \partial_{x'} K_{\mcP_t}(e,e')$ satisfies the same kind of pointwise estimates as $K_{\mcP_t}$. Moreover its first momentum is null, which is a cancellation property of ordre $1$
$$ 
\int_{\mcM} t^{1/2} \partial_{x'} K_{\mcP_t}(e,e') \nu(de') = \int_{\mcM} t^{1/2} \partial_{x'} K_{\mcP_t}(e,e') \nu(de)=0.
$$
In terms of the notation introduced and described in Appendix \ref{SectionAppendix}, the collection $(t^{1/2} \partial \mcP_t)_{t>0}$ belongs to the class ${\sf GC}^1$. That legitimates to use the notation
$$ 
\widetilde \mcQ_t:=t^{1/2} \partial \mcP_t,
$$ 
from which the representation of the statement follows. A similar observation holds for the operator ${\sf \Pi}(\partial f,\partial g)$. 
\end{Dem} 

\medskip

\noindent \textbf{\textsf{Remark --}} \textsl{In the model setting of spatial paraproducts on the one dimensional torus
$$ 
{\sf P}_{\partial f}\partial g = \sum_{k} S_{k-2}(\partial f) \Delta_k(\partial g),
$$
so an integration by parts shows that 
$$
S_{k-2}(\partial f) = 2^{k-2} {\widetilde\Delta}_{k-2}f,
$$
for some Fourier multiplier $\widetilde\Delta_{k-2}$ acting on a distribution $f(\cdot) = \sum c_n e^{in\cdot}$, as
$$
\big(\widetilde\Delta_{k-2}f\big)(x) = 2^{-(k-2)} \sum_{|n|\leq 2^{k-2}} c_n in e^{inx} = \sum_{|n|\leq 2^{k-2}} \frac{in}{2^{k-2}} c_n e^{inx},
$$
with symbol $\frac{in}{2^{k-2}} {\bf 1}_{|n|\leq 2^{k-2}}$. This symbol is not exactly supported on the annulus at scale $2^{k-2}$, as it is the case for the Fourier projector $\Delta_k$, but it satisfies some decay property at $0$ and at infinity, so it still encodes some cancellation property. We have
$$ 
{\sf P}_{\partial f}\partial g = \sum_{k} 2^{2k-2} \big({\widetilde\Delta}_kf\big) (\Delta'_kg),
$$
for operators $\Delta'_k$ perfectly localized at frequencies of scale $2^k$. The resonant operator ${\sf \Pi}(\partial f,\partial g)$ has the same structure
$$ 
{\sf \Pi}(\partial f,\partial g) = \sum_{k} 2^{2k} \widehat \Delta_k(f) \Delta'_k(g),
$$
for operators $\widehat \Delta_k$ perfectly localized at frequencies of scale $2^k$.   }
 
\medskip 

It follows from that lemma that all the continuity results of Section \ref{SubSectionCorrectors} on the corrector {\sf C} and its iterates have direct counterparts in terms of the operator $(f,g)\mapsto {\sf P}_{\partial f}g$. We single out three of them here to make that point clear. Define on the space of bounded measurable functions on the parabolic space $\mcM$ the correctors
\begin{equation*}
\begin{split}
{\sf C}^{<}_\partial\big((f_1,f_2),g\big) &:= {\sf P}_{\partial \widetilde{\sf P}_{f_1}f_2}\partial g - f_1 {\sf P}_{\partial f_2}\partial g,   \\
{\sf C}^{<}_\partial\big(f,(g_1,g_2)\big) &:= {\sf P}_{\partial f}(\widetilde{\sf P}_{g_1}g_2) - g_1 {\sf P}_{\partial f}\partial g_2,   \\
{\sf C}^{=}_\partial\big((f_1,f_2),g\big) &:= {\sf \Pi}\big(\partial {\sf P}_{f_1}f_2,\partial g\big) - f_1{\sf \Pi}(\partial f_2,\partial g).
\end{split}
\end{equation*}
We use the exponent $<$ in the notation to remind the reader from the fact that paraproducts ${\sf P}_fg$ are defined in a Fourier setting with frequencies of $f$ strictly smaller than those of $g$, while the resonant operator involves frequencies that are essentially equal. 

\medskip

\begin{thm}
\label{ThmPiPartial}
\begin{itemize}
   \item[\textcolor{gray}{$\bullet$}] Let $\alpha,\beta,\gamma$ be regularity exponents, with $\alpha\in (0,1)$ and $\alpha+\beta\leq 1$. If
   $$
   \beta+\gamma-2 < 0 < \alpha + \beta + \gamma - 2,
   $$
   then the maps
   \begin{equation*}
   \begin{split}
   \mcC^\alpha\times\mcC^\beta\times\mcC^\gamma&\rightarrow\mcC^{\alpha+\beta+\gamma-2}   \\
   (f,g,h) &\mapsto {\sf C}^<_\partial\big((f,g),h\big), {\sf C}^<_\partial\big(f,(g,h)\big)
   \end{split}
   \end{equation*}
   are continuous.   \vspace{0.15cm}
   
   \item[\textcolor{gray}{$\bullet$}] Let $\alpha,\beta,\gamma_1,\gamma_2$ be regularity exponents, with $\alpha, \gamma_1 \in (0,1)$ and $\alpha+\beta\leq 1$. If
   $$
   \beta+\gamma_1+\gamma_2 - 2 < 0 < \alpha + \beta + \gamma_1 + \gamma_2 - 2,
   $$
   then the map
   \begin{equation*}
   \begin{split}
   \mcC^\alpha\times\mcC^\beta\times\mcC^{\gamma_1}\times\mcC^{\gamma_2}&\rightarrow\mcC^{\alpha+\beta+\gamma_1+\gamma_2-2}   \\
   \big((f,g),(h_1,h_2)\big) &\mapsto {\sf C}^<_\partial\big((f,g),\widetilde{\sf P}_{h_1}h_2\big) - h_1 \, {\sf C}_\partial\big((f,g),h_2\big)
   \end{split}
   \end{equation*}
   is continuous.
\end{itemize}
\end{thm}

\medskip

\begin{Dem}
Let us concentrate first on the first statement, in the model case of the flat torus, where
$$ 
{\sf C}^<_\partial\big((f,g),h\big)= \sum 2^{2k-2} \left(\widetilde \Delta_{k-2}(\Pi^0_fg) -f \big(\widetilde \Delta_{k-2}g\big)\right) \Delta'_k(h).
$$
Note that since
$$
{\widetilde \Delta}_{k-2}\big(\Pi^0_fg\big) - f\widetilde \Delta_{k-2}g = \sum_{\ell\leq k-2} 2^{-(k-2-\ell)} \big(S_\ell f - f\big) \Delta_\ell g,
$$
we have the estimate
\begin{align*}
\Big\| \widetilde \Delta_{k-2}\big(\Pi^0_fg\big) - f\big(\widetilde \Delta_{k-2}g\big) \Big\|_\infty &\lesssim \sum_{\ell \leq k-2} 2^{-(k-2-\ell)} 2^{-\ell (\alpha+\beta)} \|f\|_{C^\alpha} \|g\|_{C^\beta}   \\
&\lesssim 2^{-k(\alpha+\beta)} \|f\|_{C^\alpha} \|g\|_{C^\beta}, 
\end{align*}
since $\alpha+\beta\leq 1$. We see here the importance of the decay of the symbol of the operator $\widetilde\Delta_{k-2}$, encoded in the factor $2^{-(k-2-\ell)}$. The conclusion follows then from the estimate
\begin{align*}
\left\|\Delta_n\Big({\sf C}^<_\partial\big((f,g),h\big)\Big)\right\|_\infty & \lesssim \left(\sum_{k\geq n-1} 2^{-k(\alpha+\beta+\gamma-2)} + \sum_{k\leq n-1} 2^{-n\alpha} 2^{-k(\beta+\gamma-2)} \right) \|f\|_{C^\alpha} \|g\|_{C^\beta} \|h\|_{C^\gamma}   \\
& \lesssim 2^{-n(\alpha+\beta+\gamma-2)} \|f\|_{C^\alpha} \|g\|_{C^\beta} \|h\|_{C^\gamma}.
\end{align*}

\ssk

If now $h = \Pi^0_{h_1}h_2$, then $\Delta_k(h) \simeq (S_{k-2}h_1) (\Delta_k h_2)$, and we have
\begin{align*}
{\sf C}^<_\partial\big((f,g),{\sf P}_{h_1}h_2\big) - h_1\,{\sf C}_\partial\big((f,g),h_2\big) = \sum 2^{k-2} \left(\widetilde \Delta_{k-2}({\sf P}_fg) -f \widetilde \Delta_{k-2}g\right) \Delta_k(h_2) \big(h_1 - S_{k-2}h_1\big)
\end{align*}
and we may conclude by the same reasoning as above, with the extra exponents coming from the positive regularity of $h_1$, since
$$ 
\big\|h_1-S_{k-2}h_1\big\|_\infty \lesssim 2^{-k\gamma_1} \, \|h_1\|_{C^{\gamma_1}}.
$$
\end{Dem}

\medskip

We let the reader state and prove the other continuity results for the iterated versions of ${\sf C}^<_\partial$ and ${\sf C}^=_\partial$. Set 
\begin{equation}\label{EqIndexSetA}
\mathscr{A}:=\{i,ij,ijk\}_{i,(i+j),(i+j+k)\leq 3},
\end{equation}
and set 
$$
\vert i\vert := i, \quad\vert ij\vert := i+j, \quad\vert ijk\vert := i+j+k.
$$

\medskip

\begin{prop}
\label{PropDecompositionGuDuDu}
Let $u$ be given by a paracontrolled system with reference functions $\bfZ$, together with a function $g\in C^4_b$. Then one can write the product $g(u)(\partial u)^2$ under the form 
\begin{equation}
\label{EqDecompositionGuDuDu}
\begin{split}
g(u)(\partial u)^2 &= \underset{\vert c_1\vert +\vert c_2\vert\leq 4}{\sum_{c_1,c_2\in\mathscr{A}}} {\sf P}_{g(u)u_{c_1}u_{c_2}}X_{c_1c_2} + \underset{\vert a\vert + \vert c_1\vert +\vert c_2\vert\leq 4}{\sum_{a, c_1,c_2\in\mathscr{A}}} {\sf P}_{g'(u)u_au_{c_1}u_{c_2}}X_{ac_1c_2}   \\
&\quad+ \frac{1}{2!}\,{\sf P}_{g^{(2)}(u)u_1^4}X_{1111} + (5\alpha-2)
\end{split}
\end{equation}   
for some remainder term of regularity $5\alpha-2$, and some explicit reference distributions $X$'s that depend only on the noise $\zeta$ and $\bfZ$.
\end{prop}

\medskip

As in Corollary \ref{CorDecompositionFuZeta} on the paracontrolled representation of $f(u)\zeta$, the proper interpretation of the above statement is that the remainder is of parabolic regularity $5\alpha-2$, provided the $X$'s are well-defined as elements of their natural spaces. In the present work, these $X$'s are given by the enhancement $\widehat{\zeta}$ of the noise.

\medskip

\begin{Dem}
If one is not interested in the precise form of the $X_a$ in \eqref{EqDecompositionGuDuDu}, one can proceed very efficiently making only computations with ${\sf E}$ and ${\sf F}$ notations, as in the proof of Corollary \ref{CorDecompositionFuZeta}. We first provide a multiplicative decomposition for 
$$
(\partial u)^2 = (\partial u)^2 = 2{\sf P}_{\partial u}\partial u + {\sf \Pi}(\partial u,\partial u) =: {\sf E}^{-2}(u,u)
$$ 
as follows. The function ${\sf E}^{-2}$ is a function of ${\sf E}$-type with respect to its two arguments. Recall that two ${\sf E}$ operators in the same identity may mean different ${\sf E}$-type operators. Distributions of regularity $\beta$ that do not depend on $u$ are denoted by $X^\beta$. Write $\simeq$ to mean equality up to a remainder term of regularity $5\alpha-2$. One has
\begin{equation*}\begin{split}
(\partial u)^2 &= {\sf E}^{-2}(u,u) \\
					  &= u_i {\sf E}^{i\alpha-2}(u) + {\sf E}^{i\alpha-2}(u_i,u)   \\
					  &= u_iu_jX^{(i+j)\alpha-2} + {\sf E}^{(i+j)\alpha-2}(u_j) + u_{ij}{\sf E}^{(i+j)\alpha-2}(u)   \\
					  &= u_iu_jX^{(i+j)\alpha-2} + u_iu_{jk}X^{(i+j+k)\alpha-2} + u_i{\sf E}^{(i+j+k)\alpha-2}(u_{jk})   \\
					  &\quad+u_{ij}u_k X^{(i+j+k)\alpha-2} + u_{ij}{\sf E}^{(i+j+k)\alpha-2}(u_k)   \\
					  &\simeq u_iu_jX^{(i+j)\alpha-2} +  u_iu_{jk}X^{(i+j+k)\alpha-2} + u_iu_{jk\ell} X^{(i+j+k+\ell)\alpha-2}   \\
					  &\quad+u_{ij}u_k X^{(i+j+k)\alpha-2} + u_{ij}u_{k\ell}{\sf E}^{(i+j+k+\ell)\alpha-2}   \\
					  &\simeq u_iu_jX^{(i+j)\alpha-2} +  u_iu_{jk}X^{(i+j+k)\alpha-2} + \big(u_1u_{111}+u_{11}^2\big) X^{4\alpha-2}   \\
\end{split}\end{equation*}
The above implicit sums are restricted to $(i+j), (i+j+k), (i+j+k+\ell)\leq 4$. Look now at the term $g(u)u_1^2 X^{2\alpha-2}$ and show that it can be written in the form \eqref{EqDecompositionGuDuDu}; the other terms are easier to deal with. We have
\begin{equation}\label{EqDecompositionProofGuDuDu}
g(u)u_1^2 X^{2\alpha-2} = {\sf P}_{g(u)u_1^2}X^{2\alpha-2} + {\sf F}^{2\alpha-2}\big(g(u)u_1^2\big) + {\sf E}^{2\alpha-2}\big(g(u)u_1^2\big),
\end{equation}
with 
\begin{equation*}\begin{split}
{\sf F}^{2\alpha-2}\big(g(u)u_1^2\big) &= {\sf P}_{g(u)}{\sf F}^{2\alpha-2}\big(u_1^2\big) + {\sf P}_{u_1^2}{\sf F}^{2\alpha-2}\big(g(u)\big) + {\sf F}^{2\alpha-2}(g(u),u_1^2)   \\
&= {\sf P}_{g(u)}{\sf P}_{u_1}{\sf F}^{2\alpha-2}\big(u_1\big) + {\sf P}_{g(u)}{\sf F}^{2\alpha-2}(u_1,u_1) + {\sf P}_{u_1^2}{\sf F}^{2\alpha-2}\big(g(u)\big)   \\
&\quad+ {\sf P}_{g'(u)u_1^2u_{11}}X^{4\alpha-2} + (5\alpha-2)   \\
&= {\sf P}_{g(u)u_1}{\sf F}^{2\alpha-2}\big(u_1\big) + {\sf P}_{g(u)u_{11}^2}X^{4\alpha-2} + {\sf P}_{u_1^2}{\sf F}^{2\alpha-2}\big(g(u)\big)   \\
&\quad+ {\sf P}_{g'(u)u_1^2u_{11}}X^{4\alpha-2} + (5\alpha-2).
\end{split}\end{equation*}
One writes the first term in the right hand side in the good form 
\begin{equation*}\begin{split}
{\sf P}_{g(u)u_1}{\sf F}^{2\alpha-2}\big(u_1\big) &= {\sf P}_{g(u)u_1}{\sf P}_{u_{1j}}X^{(2+j)\alpha-2} + {\sf F}^{3\alpha-2}(u_{11})   \\
																	&= {\sf P}_{g(u)u_1u_{1j}}X^{(2+j)\alpha-2} + {\sf P}_{u_{111}}X^{4\alpha-2} + (5\alpha-2).
\end{split}\end{equation*}
We use the paracontrolled expansion 
$$
g(u) = {\sf P}_{g'(u)u_i}Z_i + {\sf P}_{g'(u)u_{11}}{\sf R}(1,Z_1,Z_1) + \frac{1}{2}\,{\sf P}_{g^{(2)}(u)u_1^2}{\sf \Pi}(Z_1,Z_1) + (3\alpha),
$$
of $g(u)$ to order $2$ to get
\begin{equation*}\begin{split}
{\sf P}_{u_1^2}{\sf F}^{2\alpha-2}\big(g(u)\big) &= {\sf P}_{g'(u)u_1^3}X^{3\alpha-2} + {\sf P}_{g^{(2)}(u)u_1^3u_{11}}X^{4\alpha-2} + {\sf P}_{g'(u)u_2u_1^2}X^{4\alpha-2}   \\
&\quad+ {\sf P}_{g'(u)u_1^2u_{11}}X^{4\alpha-2} + {\sf P}_{g^{(2)}(u)u_1^4}X^{4\alpha-2} + (5\alpha-2).
\end{split}\end{equation*}
The term ${\sf E}^{2\alpha-2}\big(g(u)u_1^2\big)$ in \eqref{EqDecompositionProofGuDuDu} is dealt with similarly, using the ${\sf E}$-type expansion rule
\begin{equation*}\begin{split}
{\sf E}^{2\alpha-2}\big(g(u)u_1^2\big) &= g(u){\sf E}^{2\alpha-2}\big(u_1^2\big) +  u_1^2{\sf E}^{2\alpha-2}\big(g(u)\big) + {\sf E}^{2\alpha-2}\big(g(u),u_1^2\big) + {\sf E}^{2\alpha-2}\big({\sf \Pi}(g(u),u_1^2)\big)  
\end{split}\end{equation*}
and the above second order paracontrolled expansion for $g(u)$, to put ${\sf E}^{2\alpha-2}\big(g(u)u_1^2\big)$ in multiplicative form first, and then use the $F$-type expansion to put it in the form of Equation \eqref{EqDecompositionGuDuDu}. Details are left to the reader.
\end{Dem}

\medskip

\noindent \textsf{\textbf{Remark --}} {\sl Similarly to Proposition \ref{CorDecompositionFuZeta}, Proposition \ref{PropDecompositionGuDuDu}  makes perfect sense from a regularity structures point of view. Note also that the above proof implicitely uses a refined version of the ${\sf C}_\partial^<$ corrector to deal with a first argument $u_1^\sharp$ of regularity strictly bigger than $1$. This term corresponds to the term $2g(u)f(u)f'(u)u'\partial\mathcal{I}(X\circ)\partial\mathcal{I}(\circ)$, that appears in the regularity structures analysis.   }

\vspace{0.6cm}

We summarize here the notations introduced in this section.
\begin{equation*}
\begin{split}
{\sf C}(f,g,h) &= {\sf \Pi}\big(\widetilde{\sf P}_f g , h\Big) - f \, {\sf \Pi}(g,h),   \\
{\sf D}(f,g,h) &= {\sf \Pi}\big(\widetilde{\sf P}_f g , h\Big) - {\sf P}_f{\sf \Pi}(g,h),   \\
{\sf S}(f,g,h) &= {\sf P}_f\Big(\widetilde{\sf P}_g h\Big) - {\sf P}_g\Big({\sf P}_f h\Big),   \\
{\sf R}(f,g,h) &= {\sf P}_f\Big(\widetilde{\sf P}_gh\Big) - {\sf P}_{fg}h.
\end{split}
\end{equation*}
We use the same letters for the iterates of these operators, such as they were defined above. The ${}^\circ$ operators are defined by the corresponding formulas where $\widetilde{\sf P}$ is replaced by ${\sf P}$. To deal specifically with the derivatives, we have the operators 
\begin{equation*}
\begin{split}
{\sf C}_\partial^<\big((f_1,f_2),g\big) &:= {\sf P}_{\partial \widetilde{\sf P}_{f_1}f_2}\partial g - f_1 {\sf P}_{\partial f_2}\partial g,   \\
{\sf C}^<_\partial\big(f,(g_1,g_2)\big) &:= {\sf P}_{\partial f}\big(\partial\widetilde{\sf P}_{g_1}g_2\big) - g_1 {\sf P}_{\partial f}\partial g_2,   \\
{\sf C}^=_\partial(f_1,f_2,g) &:= {\sf \Pi}\big(\partial {\sf P}_{f_1}f_2,\partial g\big) - f_1{\sf \Pi}(\partial f_2,\partial g),
\end{split}
\end{equation*}
and their iterates. The only thing that matters from a computational point of view is to identify which operators are of ${\sf E}$-type, and which operators are of ${\sf F}$-type.

\bigskip

\section[\hspace{0.7cm} Nonlinear singular PDEs]{Nonlinear singular PDEs}
\label{SectionNonlinearPDEs}

We choose to illustrate the use of paracontrolled calculus for the study of singular partial differential equations on the examples of the $3$-dimensional generalised parabolic Anderson model equation
$$
\mathscr{L}u = f(u)\zeta,
$$
and the generalised (KPZ) equation
$$
\mathscr{L}u = f(u)\zeta + g(u)(\partial u)^2.
$$
on the one-dimensional torus. We introduce the notion of enhanced noise, and consistent enhancement, in Section \ref{SubsectionEnhancedDistribution}. We define in Section \ref{SubsectionFixedPoint} the maps that are used to give a fixed point formulation of the (gPAM) and (gKPZ) equations in Section \ref{SubsectionSolvingEquation}.

\bigskip

\subsection[\hspace{-0.3cm} Enhanced noise]{Enhanced noise}
\label{SubsectionEnhancedDistribution}

The archetype of singular equation is given by the controlled ordinary differential equation 
\begin{equation}
\label{EqControlledODE}
dx_t = V(x_t)dh_t,
\end{equation}
where $h$ is a non-differentiable $\bfR^\ell$-valued control and $V$ an $\textrm{L}(\bfR^\ell,\bfR^d)$-valued one form on $\bfR^d$, say. Think of a Brownian path for the control $h$. One of the deepest insights of T. Lyons in his theory of rough paths \cite{Lyons98} was to understand that one needs to change the notion of control to make sense of such an equation, and that this enhanced control takes values in a very specific universal algebraic structure. In simple terms, the enhanced control consists of $h$ and the collection of a number of objects playing the role of the non-existing iterated integrals $\int_{s\leq s_1\leq \cdots\leq s_k\leq t} dh_{s_1}\otimes\,\cdots\,\otimes dh_{s_1}$ -- such iterated integrals cannot be defined as continuous functions of their integrands, here $(h,\dots,h)$, if $h$ is not sufficiently regular; see Proposition 1.29 in \cite{LyonsStFlour}. Once given these extra data, one can make sense of, and solve uniquely, the controlled ordinary differential equation \eqref{EqControlledODE} under some appropriate regularity conditions on the one form $V$, and the solution path happens to be a continuous function of the enhanced control, in some appropriate topology. The enhancement of the control cannot be made on a purely analytic basis and requires some extra input, typically the use of probabilistic methods when the control $h$ is random.

\medskip

Hairer's theory of regularity structures provides a conceptually close framework for the study of a large class of singular partial differential equations containing the generalised parabolic Anderson model equation 
\begin{equation}\label{EqGPAM}
\mathscr{L}u = f(u)\zeta.
\end{equation}
as a particular case. To make sense of equation \eqref{EqGPAM}, one needs to enhance the distribution $\zeta$ with the a priori datum of a number of other distributions. Contrary to the case of the controlled ordinary differential equation \eqref{EqControlledODE}, this enhanced 'control' takes values in an equation-dependent algebraic structure. The resolution process is also different, as the equation is first recast in some abstract space of jets of solutions, where it can be solved under appropriate conditions. This corresponds to looking for a solution in a specific space of distributions where one can actually make sense of all the terms in the equation, especially some a priori ill-defined products. A fundamental tool, the reconstruction operator, allows then to associate to this abstract solution a classical distribution. The equation-dependent algebraic structure in which the enhanced distribution lives also allows to give sense to this solution distribution as a limit of solutions to some family of classically well-posed equations in which the distribution $\zeta$ has been smoothened. The latter point is related to renormalisation matters.

\bigskip 

Recall now Proposition \ref{PropDecompositionGuDuDu}, and Remark 1 following its proof, giving the paracontrolled expansion \eqref{EqDecompositionGuDuDu} of $f(u)\zeta$, under the assumption that a number of quantities $(Y_{k\alpha-2})_{2\leq k\leq 4}$ are given a priori as elements of the $\mathscr{L}(\mcC^{k\alpha})$-spaces, as measurable functions of $\zeta$ and $\bfZ$; write $Y=:Y(\zeta,\bfZ)$. A choice of $\frak{Z}_k\in\mcC^{k\alpha}$, with $2\leq k\leq 4$, and 
$$
\mathscr{L}\frak{Z}_k = Y_{k\alpha-2},
$$
defines an \textsf{\textbf{enhancement $\widehat{\zeta} = (\zeta, \frak{Z})$ of the noise}} $\zeta$ for the generalised parabolic Anderson model equation
$$
\mathscr{L}u = f(u)\zeta.
$$
Enlarge the finite collection $\frak{Z}_2$ by adding ${\sf \Pi}(Z_1,Z_1)$ and ${\sf R}(1,Z_1,Z_1)$ into it. This defines the collection $\overline{\frak{Z}}_2$. Set $\overline{\frak{Z}}_k := \frak{Z}_k$, for $3\leq k\leq 4$. We define this way the finite collection $\overline{\frak{Z}}$. An enhancement $\widehat{\zeta}$ is said to be \textsf{\textbf{coherent}} if 
$$
\mathscr{L}\frak{Z}_j = Y_{j\alpha-2}\big(\zeta,(\overline{\frak{Z}}_\ell)_{2\leq \ell\leq 3}\big), \qquad 2\leq j\leq 3.
$$
A coherent enhancement of the noise $\zeta$ for the (gPAM) equation can be used to work with ${\bfZ}=\overline{\frak{Z}}$ in a paracontrolled setting. The introduction of ${\sf \Pi}(Z_1,Z_1)$ and ${\sf R}(1,Z_1,Z_1)$ in the system is necessary not only to give a paracontrolled represention of the solution $u$ of the equation but also a lower order paracontrolled representation of its derivatives, amongst which $f(u)$. This is reminiscent of the use of two spaces of trees in the regularity structures setting.

\medskip

The study of the generalised (KPZ) equation 
$$
\mathscr{L}u = f(u)\zeta + g(u)(\partial u)^2,
$$
requires the introduction of further quantities ($X_{c_1c_2}, X_{ac_1c_2}, X_{1111},\dots$) that appear in \eqref{EqDecompositionGuDuDu} and the proof of this decomposition. They are assumed to be elements of $\mathscr{L}(\mcC^{k\alpha})$-spaces, and measurable functions of $\zeta$ and $\bfZ$. Write
$$
X_{k\alpha-2} = X_{k\alpha-2}(\zeta,\bfZ),
$$ 
and recall that $Z_2, Z_3$ may stand for tuples $(Z_2^{n_2}),( Z_3^{n_3})$. The joint choice of $\overline{\frak{Z}}$ and elements $\frak{X}^k\in\mcC^{k\alpha}$ such that 
$$
\mathscr{L}\frak{X}_k = X_{k\alpha-2}(\zeta,\bfZ),
$$
defines an \textsf{\textbf{enhancement $\widehat{\zeta} = (\zeta, \frak{Z}, \frak{X})$ of the noise}} $\zeta$ for the generalised (KPZ) equation. This enhancement is said to be coherent if 
$$
\mathscr{L}\frak{X}_k = X_{k\alpha-2}\Big(\zeta,(\overline{\frak{Z}}^\ell)_{2\leq \ell\leq 3}, (\frak{X}^\ell)_{2\leq \ell\leq 3})\Big).
$$

\medskip

In accordance with the regularity structures picture, the different terms that form the enhanced noise correspond to the different pieces of a paracontrolled/resonant expansion of all the formal products that appear in the tree of negative homogeneity in the regularity structures expansions of the right hand sides of the generalised (PAM) and (KPZ) equations. Constructing coherent enhancements is the task of renormalisation of stochastic singular partial differential equations. This has been implemented in the setting of regularity structures in the groundbreaking works \cite{BrunedHairerZambotti, ChandraHairer} of Hairer and co-authors. This work has not been done yet in the paracontrolled setting, but the use of Bailleul and Hoshino's results in \cite{BailleulHoshino} allows to construct coherent enhancements from the renormalised models built from \cite{ChandraHairer}.

\bigskip

\subsection[\hspace{-0.3cm} Fixed point equation]{Fixed point equation}
\label{SubsectionFixedPoint}

We give in this section the fixed point formulations of the generalised (PAM) and (KPZ) equations. In both cases, we assume a coherent enhancement $\widehat{\zeta}$ of the noise $\zeta$ is given, and work with the associated reference system $\bfZ$. Recall $\frac{2}{5}<\alpha\leq\frac{1}{2}$, so $5\alpha-2>0$. To avoid working with time-weighted H\"older spaces with exploding weights, we assume the initial conditions $u_0\in C^{4\alpha}$, to have $Pu_0\in\mcC^{4\alpha}$, and treat this term in the integral formulation of the equation as a remainder term -- recall $P$ stands for the heat propagator. See \cite{GIP, BailleulDebusscheHofmanova} or \cite{LNPCNicolas} for a sample of works where $u_0\in C^\alpha$, in a first order paracontrolled setting.

\bigskip

\textbf{\textsf{4{\boldmath$.$}2{\boldmath$.$}1  Generalised (PAM) equation{\boldmath$.$}}} Recall from \eqref{EqIndexSetA} the definition of the index set $\mathscr{A}$. We obtained in Proposition \ref{CorDecompositionFuZeta} the paracontrolled decomposition
\begin{equation*}\begin{split}
f(u)\zeta &= {\sf P}_{f(u)}\zeta + {\sf P}_{f'(u)u_a}Y_a + \frac{1}{2!}\,{\sf P}_{f^{(2)}(u)u_1^2}Y_{11} + (\flat) + {\sf P}_{h_\ell(\widehat{u})}Y^{(\ell)}_{4\alpha-2}  + (5\alpha-2)(\widehat{u}),
\end{split}\end{equation*}
with implicit sums restricted to $\vert a\vert\leq 2\alpha$, and $(\flat)\in\mathscr{L}(\mcC^{4\alpha})$. The function $(5\alpha-2)(\widehat{u})\in\mcC^{5\alpha-2}$, is a locally Lipschitz function of $\widehat{u}$. Define $v\in\mcC^\alpha$, setting, with obvious notations,
\begin{equation}\label{EqDefnV}\begin{split}
v &:= Pu_0 + \mathscr{L}^{-1}\big(f(u)\zeta\big)   \\
   &\;= \widetilde{\sf P}_{f(u)}Z_1 + \widetilde{\sf P}_{f'(u)u_a}Z_{\vert a\vert+1} + \frac{1}{2!}\,\widetilde{\sf P}_{f^{(2)}(u)u_1^2}Z_3^{(1)} + \mathscr{L}^{-1}(\flat)  \\
   &\quad+ \widetilde{\sf P}_{h_\ell(\widehat{u})}Z_4^{(\ell)} + \mathscr{L}^{-1}\big((5\alpha-2)(\widehat{u})\big) + Pu_0.
\end{split}\end{equation}
Note that $f(u)$ has a second order paracontrolled expansion
\begin{equation}\label{EqPCFu}\begin{split}
f(u) &= \widetilde{\sf P}_{f'(u)u_1}Z_1 + \widetilde{\sf P}_{f'(u)u_2}Z_2 + \widetilde{\sf P}_{f^{(2)}(u)u_1^2+2f'(u)u_{11}}{\sf R}(1,Z_1,Z_1)   \\
&\quad+ \frac{1}{2!}\,\widetilde{\sf P}_{f^{(2)}(u)u_1^2}{\sf \Pi}(Z_1,Z_1) + f(u)^\sharp(\widehat{u}),   \\
f'(u)u_1 &= \widetilde{\sf P}_{f^{(2)}(u)u_1^2+f'(u)u_{11}}Z_1 + \big(f'(u)u_1\big)^\sharp(\widehat{u}).
\end{split}\end{equation}
Equations \eqref{EqDefnV} and \eqref{EqPCFu} show that $v$ is the first component of a paracontrolled system $\widehat{v}$. Write $\widehat{u}^\sharp$ for the collection of all the remainders that define $\widehat{u}$, and $\widehat{v}^\sharp$ for the collection of all the remainders that define $\widehat{v}$. Set
\begin{equation}\label{EqDefnPhi}
\Phi(\widehat{u}^\sharp) := \widehat{v}^\sharp.
\end{equation}

\bigskip

\textbf{\textsf{4{\boldmath$.$}2{\boldmath$.$}2  Generalised (KPZ) equation{\boldmath$.$}}} Set 
\begin{equation}\label{EqDefnMathscrB}
\mathscr{B} := \Big\{(c_1c_2), (c_1c_2c_3), (c_1c_2c_3c_4)\,;\, c_i\in\mathscr{A}, \vert c_1\vert+\vert c_2\vert, \vert c_1\vert+\vert c_2\vert+\vert c_3\vert, \vert c_1\vert+\vert c_2\vert+\vert c_3\vert+\vert c_4\vert \leq 4\alpha \Big\}.
\end{equation}
Proposition \ref{CorDecompositionFuZeta} and Proposition \ref{PropDecompositionGuDuDu} together give the decomposition 
\begin{equation}\label{EqPCDecompositiongKPZ}\begin{split}
f(u)\zeta + g(u)(\partial u)^2 &= {\sf P}_{f(u)}\zeta + \sum_{\vert a\vert\leq 2\alpha}{\sf P}_{f'(u)u_a}Y_a + \frac{1}{2!}\,{\sf P}_{f^{(2)}(u)u_1^2}Y_{11} + (\flat) + {\sf P}_{h_\ell(\widehat{u})}Y_\ell^{4\alpha-2}   \\
												&\quad+ \sum_{b\in\mathscr{B}} {\sf P}_{g_b(\widehat{u})} X_b + 	(5\alpha-2)(\widehat{u}),
\end{split}\end{equation}
for the explicit functions $g_b$ of $\widehat{u}$ that appear in formula \eqref{EqDecompositionGuDuDu}. We have for instance $h_{11}(\widehat{u})=g(u)u_1^2$, and one checks that $h_{11}(\widehat{u})$ has a second order paracontrolled expansion, using the paracontrolled expansion of $g(u)$. One has
\begin{equation}\label{EqExpansionGuU_12}
g(u)u_1^2 = \widetilde{\sf P}_{(\star)_1}Z_1 + \widetilde{\sf P}_{g'(u)u_1^2u_2}Z_2 + \widetilde{\sf P}_{(\star)_2}{\sf \Pi}(Z_1,Z_1) + \widetilde{\sf P}_{(\star)_3}{\sf R}(1, Z_1,Z_1) + (3\alpha)(\widehat{u}),
\end{equation}
with 
\begin{equation*}\begin{split}
&(\star)_1 := 2g(u)u_1u_{11}+g'(u)u_1^3,   \\
&(\star)_2 := g(u)u_{11}^2+\frac{1}{2}\,g^{(2)}(u)u_1^4+2g'(u)u_1^2u_{11}   \\
&(\star)_3 := f^{(2)}(u)u_1^4+f'(u)u_1^2u_{11} + g^{(2)}(u)u_1^4+g'(u)u_1^2u_{11}   \\
&\hspace{1cm}+ 2g'(u)u_1^2u_{11}+2g(u)u_{11}^2+2g(u)u_1u_{111}+g^{(2)}(u)u_1^4+3g'(u)u_1^2u_{11},
\end{split}\end{equation*}
and a remainder in $\mcC^{3\alpha}$ that is a locally Lipscthiz function of $\widehat{u}$. (Note that the two reference functions ${\sf \Pi}(Z_1,Z_1)$ and ${\sf R}(1,Z_1,Z_1)$ are already in the $Z_2$-collection, since we are working with a coherent enhancement of the noise. They appear separetely in \eqref{EqExpansionGuU_12} as a result of computations.) The other functions $g_b(\widehat{u})$ that appear in formula \eqref{EqPCDecompositiongKPZ} all have a first order paracontrolled expansion, obtained by elementary means. Setting 
$$
v := Pu_0 + \mathscr{L}^{-1}\big(f(u)\zeta + g(u)(\partial u)^2\big)
$$
defines the first component of a paracontrolled system. (Recall we assume $u_0\in C^{4\alpha}$, so $Pu_0\in\mcC^{4\alpha}$, and one can treat it as a remainder term.) Write $\widehat{v}^\sharp$ for the collection of all the remainders that define $\widehat{v}$, and set
\begin{equation}\label{EqDefnPhiKPZ}
\Phi(\widehat{u}^\sharp) := \widehat{v}^\sharp.
\end{equation}

\bigskip

\subsection[\hspace{-0.3cm} Solving the equation]{Solving the equation}
\label{SubsectionSolvingEquation}

Recall we work on the parabolic space $[0,T)\times M$, for some possibly infinite positive time horizon $T$.

\bigskip

\textbf{\textsf{4{\boldmath$.$}3{\boldmath$.$}1  Generalised (PAM) equation{\boldmath$.$}}} Let an initial condition $u_0\in C^{4\alpha}(M)$ for (gPAM) equation be given. Assume a coherent enhancement $\widehat{\zeta}$ of the noise $\zeta$ is given. We work with the associated reference system ${\bfZ} = (Z_1,Z_2,Z_3)\in\prod_{k=1}^3\mcC^{k\alpha}$. Pick regularity exponents $(\beta_a)_{a\in\mathscr{A}}$ such that
\begin{equation}\label{EqChoiceExponentsBeta}
\frac{2}{5}<\beta_i<\beta_{ij}<\beta_{ijk}<\alpha, \quad (i), (ij), (ijk)\in\mathscr{A},
\end{equation}
and set $\beta := \min_{1\leq i\leq 3}\beta_i$. Define a map 
\begin{equation*}\begin{split}
\prod_{a\in\mathscr{A}}\mcC^{(3-\vert a\vert)\alpha+\beta_a}\times\mcC^{3\alpha+\beta} &\rightarrow \mcC^\alpha\times\prod_{a\in\mathscr{A}}\mcC^\alpha   \\
\big((u_a^\sharp)_{a\in\mathscr{A}},u^\sharp\big) &\mapsto \widehat{u}:=\big(u,(u_a)_{a\in\mathscr{A}}\big),
\end{split}\end{equation*}
given by the paracontrolled system 
\begin{equation*}
\begin{split}
&u = \sum_{i=1..3}\widetilde{\sf P}_{u_i}Z_i + u^\sharp,   \\
&u_i = \sum_{i+j=1..3} \widetilde{\sf P}_{u_{ij}}Z_j + u_i^\sharp,   \\
&u_{ij} = \sum_{i+j+k=1..3} \widetilde{\sf P}_{u_{ijk}}Z_k + u_{ij}^\sharp,   \\
&u_{ijk} = u_{ijk}^\sharp.
\end{split}
\end{equation*}
Write $\widehat{u}^\sharp$ for $\big((u_a^\sharp)_{a\in\mathscr{A}},u^\sharp\big)$. Recall the synthetic form \eqref{EqSyntheticFUZeta} obtained in Section \textbf{\textsf{3{\boldmath$.$}2{\boldmath$.$}3}} for $f(u)\zeta$, and the second order paracontrolled expansion \eqref{EqPCFu} for $f(u)$. One reads on these formulas the fact that a 'solution' of the equation needs to satisfy the constraint
$$
u_a = h_a(u),
$$
for explicit functions $h_a$. (A proper definition of a solution to the (gPAM) equation is given in Definition \ref{DefnSolutionGPAM} below.) One has for instance 
$$
u_1=f(u), u_2^{(1)}=(f'f)(u), u_3^{(1)}=f'(u)u_2=\big((f')^2f\big)(u), u_3^{(2)}=\frac{1}{2}\,f^{(2)}(u)u_1^2=\frac{1}{2}\,\big(f^{(2)}(f')^2\big)(u),
$$ 
and so on. Note that a 'solution' $u$ has null derivatives in the ${\sf \Pi}(Z_1,Z_1)$ and ${\sf R}(1,Z_1,Z_1)$ directions. Define 
$$
\mcS^\textrm{{\tiny pam}}_T(u_0) := \Big\{\widehat{u}^\sharp\,;\,{u^\sharp_a}_{\vert \tau=0}=h_a(u_0),\,u^\sharp_{\vert \tau=0}=u_0\Big\};
$$
equipped with the natural norm induced from $\prod_{a\in\mathscr{A}}\mcC^{(3-\vert a\vert)\alpha+\beta_a}\times\mcC^{3\alpha+\beta}$, the space $\mcS^\textrm{{\tiny pam}}_T(u_0)$ is complete. Recall the definition of the map $\Phi$ defined in \eqref{EqDefnPhi}.

\medskip

\begin{prop}
\label{PropContractiongPAM}
A choice of time horizon $T$ sufficiently small ensures that $\Phi$ is a contraction of $\mcS^\textrm{{\tiny pam}}_T(u_0)$.
\end{prop}

\medskip

\begin{Dem}
The formulae defining the family of remainders $\widehat{v}^\sharp$ are actually explicit -- see the end of the proof for a sample, and one can read on them that 
$$
\Phi\big(\mcS^\textrm{{\tiny pam}}_T(u_0)\big) \subset \prod_{a\in\mathscr{A}}\mcC^{(3-\vert a\vert)\alpha+\beta'_a}\times\mcC^{3\alpha+\beta'},
$$
with $\beta_a'>\beta_a$ and $\beta'>\beta$, as a consequence of the choice of exponents \eqref{EqChoiceExponentsBeta} -- we use the regularity assumptions on the components of the enhanced noise $\widehat\zeta$ and the classical Schauder estimate on H\"older space of positive regularity. So not only do we have that $\Phi$ sends $\mcS^\textrm{{\tiny pam}}_T(u_0)$ into itself, but we also have that $\Phi$ is locally Lipschitz from $\mcS^\textrm{{\tiny pam}}_T(u_0)$ into $ \prod_{a\in\mathscr{A}}\mcC^{(3-\vert a\vert)\alpha+\beta'_a}\times\mcC^{3\alpha+\beta'}$, as a consequence of the locally Lipschitz character of the corrector and commutator and their iterates, and the refined corrector. The contraction property of $\Phi : \mcS^\textrm{{\tiny pam}}_T(u_0)\rightarrow\mcS^\textrm{{\tiny pam}}_T(u_0)$, for $T$ sufficiently small, follows then from the elementary estimate
\begin{equation}
\|w\|_{\mcC^{\delta_1}} \lesssim T^{\frac{\delta_2-\delta_1}{2}}\,\|w\|_{\mcC^{\delta_2}}, \quad w_{\vert \tau=0}=0 \label{eq:ffff}
\end{equation}
that holds for any $\delta_2>\delta_1>0$.

\medskip

Here are explicit formulae for the components $v_1^\sharp$ and $v^\sharp_2=v_{11}^\sharp$ of $\widehat{v}^\sharp$. The term $v_1^\sharp$ is the remainder in the second order paracontrolled expansion of $f(u)$. Denote by $(3\alpha)_f(u)$ the $3\alpha$-remainder in the second order paracontrolled expansion for $f(u)$, defined by 
$$
f(u) = {\sf P}_{f'(u)}u + \frac{1}{2}\,{\sf P}_{f^{(2)}(u)}{\sf \Pi}(u,u) + (3\alpha)_f(u).
$$
The function $(3\alpha)_f$ sends continuously $\mcC^\alpha$ into $\mcC^{3\alpha}$. One has
\begin{equation*}\begin{split}
v_1^\sharp &= \sum_{i=1}^3{\sf R}^o\big(f'(u),u_i,Z_i\big) + \sum_{j=2}^3{\sf P}_{f'(u)}{\sf R}(1,u_j,Z_j) + {\sf P}_{f'(u)}{\sf R}\big(1,\widetilde{\sf P}_{u_{12}}Z_2+u_1^\sharp,Z_1\big)   \\
&\quad+ {\sf R}\big(1,(u_{11},Z_1),Z_1\big) + {\sf P}_{f'(u)}u^\sharp + \frac{1}{2}\,{\sf P}_{f^{(2)}(u)}{\sf D}(u_i,Z_i,u) + (3\alpha)_f(u)   \\
&\quad+ \frac{1}{2}\,{\sf R}\big(f^{(2)}(u),u_i,{\sf \Pi}(Z_i,u)\big) + {\sf R}\big(f'(u),u_{11},{\sf R}(1,Z_1,Z_1)\big) + \frac{1}{2}\sum_{j=2}^3{\sf P}_{f^{(2)}(u)u_j}{\sf \Pi}(Z_j,u)   \\
&\quad+ \frac{1}{2}\sum_{j=2}^3{\sf P}_{f^{(2)}(u)u_1}{\sf \Pi}\big(Z_1, \widetilde{\sf P}_{u_k}Z_k+u^\sharp\big) + \frac{1}{2}\,{\sf P}_{f^{(2)}(u)u_1}{\sf D}(u_1,Z_1,Z_1)   \\
&\quad+ \frac{1}{2}\,{\sf R}\Big(f^{(2)}(u)u_1,u_1,{\sf \Pi}(Z_1,Z_1)\Big) + {\sf R}\big(1,f'(u)u_2,Z_2\big) + {\sf R}\big(1,f'(u)u_{11},{\sf R}(1,Z_1,Z_1)\big)   \\
&\quad + \frac{1}{2}\,{\sf R}\big(1,f^{(2)}(u)u_1^2,{\sf \Pi}(Z_1,Z_1)\big) + {\sf R}\big(1,v_{11}^\sharp,Z_1\big) + {\sf R}\Big(1,\big(f^{(2)}(u)u_1^2+f'(u)u_{11},Z_1\big),Z_1\Big)   \\
&\quad + {\sf R}\big(1,f^{(2)}(u)u_1^2,{\sf \Pi}(Z_1,Z_1)\big),Z_1\big).
\end{split}\end{equation*}
The term $v_2^\sharp=v_{11}^\sharp$ is the $2\alpha$-remainder in the first order paracontrolled expansion of $f'(u)u_1$. Denote by $(2\alpha)_{f'}(u)$ the $2\alpha$-remainder in the paralinearisation formula for $f'(u)$, defined by 
$$
f'(u) = {\sf P}_{f ^{(2)}(u)}u + (2\alpha)_{f'}(u).
$$ 
One has
\begin{equation*}\begin{split}
v_2^\sharp=v_{11}^\sharp &= {\sf \Pi}\big(f'(u),u_1\big) + {\sf P}_{f'(u)}\big(\widetilde{\sf P}_{u_{12}}Z_1+u_1^\sharp\big) + (2\alpha)_{f'}(u) + {\sf R}\big(f'(u),u_{11},Z_1\big)   \\
&\quad+ {\sf R}^o\big(u_1,f^{(2)}(u),u\big) + {\sf R}\big(1,f'(u)u_{11},Z_1\big) + \sum_{k=2}^3{\sf P}_{f^{(2)}(u)}\big(\widetilde{\sf P}_{u_k}Z_k+u^\sharp\big)   \\
&\quad+ {\sf R}^o\big(f^{(2)}(u)u_1,u_1,Z_1\big) + {\sf R}\big(1,f^{(2)}(u)u_1^2,Z_1\big).
\end{split}\end{equation*}
The explicit expressions for the other terms of $\widehat{v} ^\sharp$ are similar or easier to obtain. These expression are simpler in the case where $f(u)=u$, as $f'(u)=1$ and $f^{(2)}(u)=0$, and ${\sf R}^o(1,\dots)=0$ and ${\sf P}_1(*)=(*)$.
\end{Dem}

\medskip

\begin{definition}   \label{DefnSolutionGPAM}
A \textbf{\textsf{solution to the generalised (PAM) equation}} is a fixed point of $\Phi$ in $\mcS^\textrm{{\tiny pam}}_T(u_0)$.
\end{definition}

\medskip

We obtain the following well-posedness result as a direct consequence of Proposition \ref{PropContractiongPAM}.

\medskip

\begin{thm} \label{thm:PAM1}
The generalised \emph{(PAM)} equation has a unique solution in $\mcS^\textrm{{\tiny pam}}_T(u_0)$; it depends continuously on the coherent enhancement $\widehat{\zeta}$ of the noise $\zeta$.
\end{thm}

\medskip

\noindent \textsf{\textbf{Remarks.}} 
\begin{itemize}
	\item \textsl{So far, the theory of regularity structures has not been developed in a manifold setting. The recent work \cite{DDD17} of Dahlqvist-Diehl-Driver shows how this can be done in the simplest case where the noise is not too rough, corresponding in our setting to a regularity exponent $\alpha>\frac{2}{3}$. A first order description of the objects is sufficient in that setting, as was the case in our previous work \cite{BB15}, whose content covers partly their results.   }   \vspace{0.15cm}

	\item \textsl{We assume here that the initial condition is in $C^{4\alpha}$. We use that fact to put the term $P(u_0)$ in the remainder. One can improve upon this constraint on $u_0$ and only require that $u_0\in C^\alpha$, at the price of working with weighted H\"older spaces with a temporal weight, explosive at $\tau=0$, for example a space equipped with the norm
$$ 
\sup_{0<\tau \leq T} \tau^\gamma \,\|u(\tau)\|_{C^\alpha}
$$
for some $\gamma>0$.  See \cite{GIP, BailleulDebusscheHofmanova} or \cite{LNPCNicolas} for a sample of works where $u_0\in C^\alpha$, in a first order paracontrolled setting.   }
\end{itemize}

\medskip

The next statement is about the {\sl linear} (PAM) equation
\begin{equation}\label{EqPAM}
\mathscr{L}u=u\zeta.
\end{equation}

\medskip

\begin{thm} \label{thm:PAM2}
Assume the noise $\zeta$ in the \emph{(PAM)} equation is a space white noise on $M$, and the components $\frak{Z}_k, \frak{X}_k$ of the enhanced noise $\widehat{\zeta}$ are such that both $\mathscr{L}\frak{Z}_k$ and $\mathscr{L}\frak{X}_k$ take values in $L^\infty_\tau C_x^{k\alpha-2}$. Then equation \eqref{EqPAM} has a unique, global in time, solution in $\mcS^\textrm{{\tiny pam}}_\infty(u_0)$. Its restriction to any finite time interval depends continuously on the coherent enhancement $\widehat{\zeta}$ of the noise $\zeta$.
\end{thm}

\medskip

\begin{Dem} 
Define a temporal weight
\begin{equation}
w(\tau) := e^{\kappa\tau}, \label{eq:weight2}
\end{equation}
for some non-negative constant $\kappa$ to be fixed later, and work in weighted parabolic H\"older space $\mcC^\gamma_w$, described in Appendix \ref{SectionAppendix}. We let the reader to check that all the proofs of continuity for the paraproduct, resonant, correctors and commutators, still hold in the setting of weighted parabolic H\"older spaces, with estimates that are uniform in $\kappa\geq 0$. This uniform character of the continuity estimates comes from the important fact that all our approximation operators $\mcP_t,\mcQ_t, ...$ are supported only on ``past time''-half spaces. Indeed, if $T$ is an operator acting on the time variable with a kernel $K(\tau,\sigma)$ supported on the past time half space $\{(\tau,\sigma),0\leq \sigma\leq \tau\}$, then 
\begin{align*}
 w(\tau)^{-1}T(f)(\tau) = e^{-\kappa \tau} T(f)(\tau) & = \int K(\tau,\sigma) e^{-\kappa \tau} f(\sigma) \, d\sigma \\
 & = \int_{0}^{\tau} K(\tau,\sigma) e^{-\kappa (\tau-\sigma)} e^{-\kappa \sigma} f(\sigma) \, d\sigma,
 \end{align*}
 so we are reduced to the case of $T(w^{-1}f)$, up to an extra coefficient $e^{-\kappa (\tau-\sigma)}$ which is $\kappa$-uniformly bounded by $1$, because of the time-support of the kernel $K$.

\smallskip

With this point in mind, we repeat the previous computations in the weighted H\"older spaces. The assumptions on the components of the enhanced noise allow us to use Schauder estimates, Theorem \ref{ThmSchauder}, and gain an extra factor $\kappa^{-(\delta_2-\delta_1)/2}$. This is the analogue of \eqref{eq:ffff} in the above unweighted proof. Taking $\kappa$ large enough allows to compensate any 'implicit' constants in the different estimates and get the contraction property of the map $\Phi$. We then conclude the proof of the existence and uniqueness of the fixed point, for an arbitrary horizon time. This approach using weighted spaces only works for a linear application $f$.
\end{Dem}

\bigskip

\textbf{\textsf{4{\boldmath$.$}3{\boldmath$.$}2  Generalised (KPZ) equation{\boldmath$.$}}} We proceed exactly as in the preceding section, working with an initial condition $u_0\in C^{4\alpha}(\bfT)$. We assume a coherent enhancement $\widehat{\zeta}$ of the (gKPZ) noise $\zeta$ is given, and work with the associated reference system $\bfZ$. Recall from \eqref{EqDefnMathscrB} that the components of $\bfZ$ are naturally indexed by the index sets $\mathscr{A}$ and $\mathscr{B}_{\leq 3\alpha} := \big\{b\in\mathscr{B}\,;\,\vert b\vert \leq 3\alpha\big\}$. Here again, a 'formal' solution of the equation needs to satisfy the constraint
$$
u_b=h_b(u),
$$
for $b\in\mathscr{B}_{\leq 3\alpha}$ and explicit functions $h_b$, in addition to the constraints $u_a=h_a(u), \,a\in\mathscr{A}$, on the components that are common with the (gPAM) equation. One has for instance
$$
h_{(1)(1)}(u) = (gf^2)(u),\; h_{(1)(2)}(u) = \big(gf^2f'\big)(u),\;h_{(1)((1)(1))}(u)=\big(g^2f^2\big)(u).
$$
(Our index notation becomes a bit messy on concrete examples, and the tree-indexed notation of regularity structures seems more appropriate to get concise notations.) Set 
$$
\mcS^\textrm{{\tiny kpz}}_T(u_0) :=\Big\{\widehat{u}^\sharp\,;\,{u^\sharp_a}_{\vert \tau=0}=h_a(u_0),\,u^\sharp_{\vert \tau=0}=u_0\Big\}.
$$
Equipped with the natural norm induced from $\prod_{a\in\mathscr{A}\cup\mathscr{B}_{\leq 3\alpha}} \mcC^{(3-\vert a\vert)\alpha+\beta_a}\times\mcC^{3\alpha+\beta}$, the space $\mcS^\textrm{{\tiny kpz}}_T(u_0)$ is complete. Recall the definition of the map $\Psi$ defined in \eqref{EqDefnPhi}. 

\medskip

\begin{definition*}
A solution of the \emph{(gKPZ)} equation is a fixed point of $\Psi$ in $\mcS^\textrm{{\tiny kpz}}_T(u_0)$.
\end{definition*}

\medskip

The very same reasoning as above provides the local in time well-posed character of the (gKPZ) equation; details are left to the reader.

\medskip

\begin{thm} \label{thm:KPZ}
Given $u_0\in C^{4\alpha}(\bfT)$, one can choose the time horizon $T$ sufficiently small for the generalised \emph{(KPZ)} equation to have a unique solution in $\mcS^\textrm{{\tiny kpz}}_T(u_0)$. This solution depends continuously on the coherent enhancement $\widehat{\zeta}$ of the noise $\zeta$.
\end{thm}

\bigskip

\appendix

\section[\hspace{0.7cm} Parabolic setting]{Details on the parabolic setting}
\label{SectionAppendix}

For the reader's convenience, we recall in this Appendix a number of notions/facts introduced and studied in detail in our previous work \cite{BBF15}, with the hope that this will make the reading of the present work self-contained. We refer the reader to \cite{BBF15} for the proofs of the different statements given here. We describe in Section \ref{SubsectionAppendixApproxOperators} a class of operators with some cancellation property; they play here the role played by the Fourier projectors $\Delta_i$ on dyadic blocks in Littlewood-Paley theory. Parabolic H\"older spaces are described in Section \ref{SubsectionParabolicHolderSpaces}, together with Schauder estimates in this scale of spaces. We introduce the pair $\big({\sf P},\widetilde{\sf P}\big)$ of paraproducts in Section \ref{SubsectionParaproducts}. The statements given here are explicitly used in the proofs of the continuity results of Section \ref{SectionToolKit}, to be found in Appendix \ref{SectionAppendixContinuity}.

\ssk

We use the notations introduced in Section \ref{SubsectionSettingResults} and assume the operator $L$ satisfies the assumption stated there. Recall we denote by $e$ a generic element of the parabolic space $\mcM$.

\bigskip

\subsection{Approximation operators}
\label{SubsectionAppendixApproxOperators}

The use of paraproducts and other kind of singular operators involve the fundamental notion of  approximation operators, some aspects of which we discuss in this section. Recall typical space/time points are denoted by $(\tau,x)$ and $(\sigma,y)$.

\ssk

The following parabolic Gaussian-like kernels $(\mcG_t)_{0<t\leq 1}$ will be used as reference kernels. For $0<t\leq 1$ and $\sigma\leq\tau$, set 
$$
\mcG_t\big((\tau,x),(\sigma,y)\big):= \nu\Big(B_\mcM \big((\tau,x),\sqrt{t}\big) \Big)^{-1} \left(1+c\,\frac{\rho \big((\tau,x), (\sigma,y)\big)^2}{t} \right)^{-\ell_1}  
$$
and set $\mcG_t \equiv 0$ if $\tau\leq \sigma$. We do not emphasize the dependence of $\mcG$ on the positive constant $c$ in the above definition, and we shall allow ourselves to abuse notations and write $\mcG_t$ for two functions corresponding to two different values of that constant. So we have for instance, for $s,t\in (0,1)$, the estimate
\begin{equation}
\label{EqIteratedG}
\int_\mcM \mcG_{t}\big((\tau,x),(\sigma,y)\big) \, \mcG_{s}\big((\sigma,y),(\lambda,z)\big) \, \nu(d\sigma dy) \lesssim \mcG_{t+s}\big((\tau,x),(\lambda,z)\big).
\end{equation}
Presently, note that a choice of  large enough constant $\ell_1$ in the definition of $\mcG_t$ ensures that we have
\begin{equation*} 
\sup_{t\in(0,1]} \sup_{(\tau,x)\in \mcM} \ \int_\mcM \mcG_{t}\big((\tau,x),(\sigma,y)\big)\, \nu(d\sigma dy) <\infty,
\end{equation*}
so any linear operator on $\mcM$, with a kernel pointwisely bounded by some $\mcG_t$ is bounded in $L^{p}(\nu)$ for every $p\in[1,\infty]$.

\medskip

\begin{definition*}   {\it
We shall denote throughout by ${\sf G}$ the set of families $(\mcP_t)_{0<t\leq 1}$ of linear operators on $\mcM$, with kernels pointwisely bounded by
$$ 
\Big| K_{\mcP_t}(e,e') \Big| \lesssim  \mcG_{t}(e,e').
$$   }
\end{definition*}

\medskip

Given a real-valued integrable function $\phi$ on $\bfR$, set 
$$
\phi_t(\cdot) := \frac{1}{t}\, \phi\Big(\frac{\cdot}{t}\Big);
$$
the family $(\phi_t)_{0<t\leq 1}$ is uniformly bounded in $L^1(\bfR)$. We also define the ``convolution'' operator $\phi^\star$ associated with $\phi$ via the formula
$$
\phi^\star(f)(\tau) := \int_{0}^\infty \phi(\tau-\sigma)f(\sigma)d\sigma.
$$
Note that if $\phi$ has support in $\bfR_+$, then the operator $\phi^\star$ has a kernel supported on the same set $\big\{(\sigma,\tau)\,; \sigma\leq \tau\big\}$ as our Gaussian-like kernel. Moreover, we let the reader check that if $\phi_1,\phi_2$ are two functions in $L^1$, with $\phi_2$ supported on $[0,\infty)$ then 
$$ 
\big(\phi_1 \ast \phi_2\big)^\star = \phi_1^\star \circ \phi_2^\star,
$$
where $\phi_1 \ast \phi_2$ stands for the usual convolution of $\phi_1$ and $\phi_2$. 

Given an integer $b\geq 1$, we define a family of operators on $L^2(M)$ setting
$$ 
Q_t^{(b)} := \gamma_b^{-1} (tL)^b e^{-tL} \qquad \textrm{and} \qquad  -t\partial_t P_t^{(b)} = Q_t^{(b)},
$$
with $\gamma_b := (b-1)!$; so $P_t^{(b)}$ is an operator of the form $p_b(tL)e^{-tL}$, for some polynomial $p_b$ of degree $b-1$, with value $1$ in $0$. Recall from Section \ref{SubsectionSettingResults} the definition of the differential operators $V_I$, for $I$ a multiindex. Under the assumptions on $L$ stated in Section \ref{SubsectionSettingResults}, the operators $P_t^{(b)}$ and $Q_t^{(b)}$ both satisfy, for any multi-index $I$, the Gaussian regularity estimates
\begin{equation*}
\left| K_{t^\frac{| I |}{2}V_IR} (x,y) \right| \vee \left| K_{t^\frac{| I |}{2}RV_I} (x,y) \right| \lesssim  \frac{1}{\mu\big(B(x,\sqrt{t})\big)} \, e^{-c\,\frac{d(x,y)^2}{t}},
\end{equation*}
with $R$ standing here for $P_t^{(b)}$ or $Q_t^{(b)}$, as well as the following pointwise regularity estimates. For $d(x,z)\leq \sqrt{t}$, we have
\begin{equation*}
\Big| K(x,y) - K(z,y)\Big| \lesssim \frac{d(x,z)}{\sqrt{t}}  \frac{1}{V\big(x,\sqrt{t}\big)} e^{-c\,\frac{d(x,y)^2}{t}},
\end{equation*}
where $K$ is the kernel of either $t^\frac{| I |}{2}V_IR$ or $t^\frac{| I |}{2}RV_I$.
 
\medskip

The parameters $b$ and $\ell_1$ are chosen large enough and fixed once and for all -- see \cite{BBF15} to see how to choose these parameters. 
The reader should simply keep in mind that the higher $b$ and $\ell_1$ are, the higher order of regularity we can deal with. In our applications, we need all the objects to have a regularity order in the range $(-3,3)$, so $b$ and $\ell_1$ are chosen big enough to allow for this range in all the following continuity results.

\medskip

\begin{defn}
\label{DefnStandardCollection} {\it 
Let an integer $a\in\llbracket 0,2b \rrbracket$ be given. The following collection of families of operators is called the \textbf{standard collection of operators with cancellation of order $a$}, denoted by ${\sf StGC}^a$. It is made up of all the space-time operators
$$ 
\Big( \big(t^\frac{|J|}{2} V_J\big) (tL)^\frac{a-|J|-2k}{2} P_t^{(c)}  \otimes m^\star_t \Big)_{0<t\leq 1} 
$$
where $k$ is an integer with $2k+|J| \leq a$, and $c\in \llbracket 1,b\rrbracket$, and $m$ is any smooth function supported on $\big[\frac{1}{2},2\big]$ such that
\begin{equation} 
\label{eq:moment} 
\int \tau^i m(\tau) \, d\tau=0,  
\end{equation}
for all $0\leq i\leq k-1$, with its first $b$ derivatives bounded by $1$.  We also set 
$$
{\sf StGC}^{[0,2b]} := \bigcup_{0\leq a \leq 2b} {\sf StGC}^a.
$$   }
\end{defn}

\ssk

The above mentioned cancellation effect is quantified by the property \eqref{prop:cancellation} stated in Proposition \ref{prop:SO} below. Note here that it makes sense at an intuitive level to say that $L^\frac{a-|J|-2k}{2}$ encodes cancellation in the space-variable of order $a-|J|-2k$, that $V_J$ encodes a cancellation in space of order $|J|$  and that the moment condition \eqref{eq:moment} encodes a cancellation property in the time-variable of order $k$ for the convolution operator $m^\star_t$. Since we are in the parabolic scaling, a cancellation of order $k$ in time corresponds to a cancellation of order $2k$ in space, so the operator $V_J L^\frac{a-|J|-2k}{2} P_t^{(c)}  \otimes m^\star_t$ is expected to have a space-time cancellation property of order $a$. 

\ssk

\begin{definition*} {\it 
Given an operator $Q := V_I\,\phi(L)$, with $|I|\geq 1$, defined by functional calculus from some appropriate function $\phi$, we write $Q^\bullet$ for the \textbf{formal dual operator}
$$ 
Q^{\bullet} :=  \phi(L) V_I.
$$
For $I= \emptyset$, and $Q = \phi(L)$, we set $Q^\bullet := Q$. For an operator $Q$ as above we set
$$
\big(Q\otimes m^\star\big)^\bullet := Q^\bullet\otimes m^\star.
$$   }
\end{definition*}

\ssk

Note that the above definition is \textit{not} related to any classical notion of duality, and emphasize that we do \textit{not assume} that $L$ is self-adjoint in $L^2(\mu)$. This notation is only used to indicate that a $Q_t$ operator , resp. a $Q_t^\bullet$ operator, can be composed on the right, resp. on the left, by another operator $\psi(L)$, for a suitable function $\psi$, due to the functional calculus for $L$. $L$ is supposed to be sectorial and to have a holomorphic functional calculus, so for example it is known to give a sense to $\phi(L)$ for every holomorphic function $\phi$ (or function which are holomorphic and bounded in a small sector of the complex plane around $(0,\infty)$, as $\phi(z)=z^k e^{-z}$ for example).

\ssk

\begin{prop} \label{prop:SO} {\it 
Consider $\mcQ^1\in {\sf StGC}^{a_1}$ and $\mcQ^2\in {\sf StGC}^{a_2}$ two standard collections with cancellation, and set $a:=\min(a_1,a_2)$. Then for every $s,t\in(0,1]$, the operator $\mcQ^1_{s} \circ \mcQ^{2\bullet}_{t}$ has a kernel pointwisely bounded by
\begin{equation}
\label{prop:cancellation}
\left| K_{\mcQ^1_{s}\circ\mcQ^{2\bullet}_{t}} (e,e') \right| \lesssim  \left(\frac{ts}{(s+t)^2}\right)^\frac{a}{2} \mcG_{t+s}(e,e'). 
\end{equation}   }
\end{prop}

\ssk

The above mentioned {\it orthogonality} property of standard operators with cancellation is encoded in the factor $\left(\frac{ts}{(s+t)^2}\right)^\frac{a}{2}$ that appears in the above estimate. This factor is small as soon as $s$ or $t$ is small compared to the other.

\medskip

\begin{definition*} {\it
Let $a$ be an integer in $[0,2b]$. We define the subset ${\sf GC}^a$ of ${\sf G}$ of \textbf{families of operators with the cancellation property of order $a$} as the set of elements $\mcQ$ of ${\sf G}$ with the following cancellation property. For every $0<s,t\leq 1$ and every standard family $\mcS\in {\sf StGC}^{a'}$, with $a'\in\llbracket a,2b\rrbracket$, the operator $\mcQ_t\circ\mcS_s^\bullet$ has a kernel pointwisely bounded by
\begin{equation} 
\label{eq:mcOa}
\Big| K_{\mcQ_t \circ\mcS_s^\bullet}(e,e') \Big| \lesssim \left(\frac{st}{(s+t)^2}\right)^\frac{a}{2} \mcG_{t+s}(e,e'). 
\end{equation}   }
\end{definition*}

\ssk

We introduced above the operators $Q^{(b)}_t$ and $P^{(b)}_t$ acting on functions/distributions on $M$; we now introduce their parabolic counterpart. Choose arbitrarily a smooth real-valued function $\varphi$ on $\bfR$, with support in $\big[\frac{1}{2},2\big]$, unit integral and such that for every integer $k=1,..,b$ 
$$ 
\int \tau^k \varphi(\tau) \, d\tau =0.
$$ 
Set 
$$
\mcP_t^{(b)}:= P^{(b)}_t \otimes \varphi_t^\star \qquad \textrm{and} \qquad \mcQ_t^{(b)}:= -t \partial_t \mcP^{(b)}_t.
$$
An easy computation yields that
$$ 
\mcQ_t^{(b)} = Q_t^{(b)} \otimes \varphi_t^\star + P_t^{(b)} \otimes \psi_t^\star
$$
where $\psi(\sigma) = \varphi(\sigma) + \sigma\varphi'(\sigma)$. Note that, from its very definition, a parabolic operator $\mcQ_t^{(b)}$ belongs at least to ${\sf GC}^{2}$, for $b\geq 2$.  Note also that due to the normalization of $\varphi$, then for every $f\in L^p ({\mathbb R})$ supported on $[0,\infty)$, we have the convergence in $L^p$
$$
\varphi^\star_t(f) \xrightarrow[t\to 0]{} f.
$$
So, the operators $\mcP_t$ tend weakly as $t$ goes to $0$ to the identity on $L^p_{[0,\infty)}(\mcM)$, the set of functions $f\in L^p(\mcM)$ with time-support included in $[0,\infty)$, with $p \in [1,\infty)$. The same convergence holds on the set of functions $f\in C^0(\mcM)$ with time-support included in $[0,\infty)$. The following {\bf Calder\'on reproducing formula} follows as a consequence. For every continuous function $f\in L^\infty(\mcM)$ with time-support in $[0,\infty)$, then 
\begin{equation}
\label{eq:calderon}
 f = \int_0^1 \mcQ^{(b)}_t f \, \frac{dt}{t} + \mcP^{(b)}_1f.
\end{equation}
Noting that the measure $\frac{dt}{t}$ gives unit mass to intervals of the form $\big[2^{-i-1},2^{-i}\big]$, and considering the operator $\mcQ^{(b)}_t $ as a kind of multiplier roughly localized at ''frequencies'' of size $t^{-\frac{1}{2}}$, Calder\'on's formula appears as nothing else but a continuous time analogue of the Littlewood-Paley decomposition of $f$, with  $\frac{dt}{t}$ in the role of the counting measure.

\bigskip

\subsection{Parabolic H\"older spaces and Schauder estimates}
\label{SubsectionParabolicHolderSpaces}

We recall in this section the definitions and basic properties of the space and space-time weighted H\"older spaces, with possibly negative regularity index. We also recall the fundamental regularization properties of the heat operator, quantified by Schauder estimates.

\medskip

Let us start recalling the following well-known facts about H\"older space on $M$, and single out a good class of weights on the parabolic space $\mcM$. Given $0<\alpha\leq 1$, the classical metric H\"older space $H^\alpha$ is defined as the set of real-valued functions $f$ on $M$ with finite $H^\alpha$-norm, defined by the formula
$$
\|f\|_{H^\alpha} := \big\|f\big\|_{L^\infty(M)} + \sup_{x\neq y\in M} \frac{\big| f(x) - f(y) \big|}{d(x,y)^\alpha} <\infty.
$$

\begin{definition*} \label{def:C} {\it 
For $\alpha\in (-3,3)$, define $C^\alpha := C^\alpha(M)$ as the closure of the set of bounded and smooth functions for the $C^\alpha$-norm, defined by the formula
$$ 
\|f\|_{C^\alpha} := \big\|e^{-L}f\big\|_{L^\infty(M)} + \sup_{0<t\leq 1} t^{-\frac{\alpha}{2}} \big\|Q_t^{(b)}f \big\|_{L^\infty(M)}.
$$   }
\end{definition*}

\medskip

This norm does not depend on the integer $b > \frac{|\alpha|}{2}$, and the two spaces $H^\alpha$ and $C^\alpha$ coincide and have equivalent norms when $0<\alpha<1$ -- see for instance Proposition 2.5 in \cite{BB15}. These notions have parabolic counterparts which we now introduce. Recall we work with the parabolic space $\mcM = [0,T]\times M$, for a finite time horizon $T$. The introduction of a time weight in the next definition thus has no effect on the space involved, nor on its topology. Its introduction happens however to be a convenient freedom which allows to simplify a number of arguments. Let then a non-negative parameter $\kappa$ be given and denote by $w$ the weight
\begin{equation}
w(\tau) := e^{\kappa\tau}. \label{eq:weight}
\end{equation}
For $0<\alpha\leq 1$, the metric parabolic H\"older space $\mcH^\alpha = \mcH^\alpha(\mcM)$ is defined as the set of all functions on $\mcM$ with finite $\mcH^\alpha$-norm, defined by the formula
$$ 
\|f\|_{\mcH^\alpha} := \big\| w^{-1}f\big\|_{L^\infty(\mcM)} + \sup_{0<\rho((\tau,x),(\sigma,y))\leq 1;\,\tau\geq \sigma} \;\; w^{-1}(\tau)\, \frac{|f(\tau,x)-f(\sigma,y)|}{\rho\big((\tau,x),(\sigma,y)\big)^\alpha}.
$$
As in the above space setting one can recast this definition in a more functional setting, using the parabolic standard operators. A set of distributions was introduced in \cite{BBF15}, whose precise definition is irrelevant here.

\medskip

\begin{definition*} \label{def:calC} {\it 
For $\alpha\in (-3,3)$, we define the parabolic H\"older space 
$$
\mcC^\alpha = \mcC^\alpha(\mcM) = \mcC^\alpha_w(\mcM) = \mcC^\alpha_w
$$ 
as the closure, in the set of distributions, of the set of bounded and continuous functions on $\mcM$ for the $\mcC^\alpha-w$-norm, defined by 
$$ 
\|f\|_{\mcC^\alpha_w} := \sup_{\genfrac{}{}{0pt}{}{\mcQ \in {\sf StGC}^{k}}{0\leq k\leq 2b}} \big\|w^{-1}\mcQ_1(f)\big\|_{L^\infty(\mcM)} +  \sup_{\genfrac{}{}{0pt}{}{\mcQ \in {\sf StGC}^{k}}{ |\alpha|< k\leq 2b}} \ \sup_{0<t\leq 1} t^{-\frac{\alpha}{2}} \big\|w^{-1}\mcQ_t(f) \big\|_{L^\infty(\mcM)}.
$$   }
\end{definition*}

\medskip

We write $\mcC^\alpha_w$ if we want to emphasize the dependence of the norm on $w$. The following result was proved in \cite{BBF15}, building on Calder\'on's formula \eqref{eq:calderon}.

\medskip

\begin{prop} \label{prop:equivalence} {\it 
Choose any finite non-negative parameter $\kappa$ in the definition \eqref{eq:weight} of the weight $w$. Given $\alpha\in (0,2)$, set
$$ 
\mcE^\alpha := \Big(C_\tau^{\alpha/2} L^\infty_x\Big) \cap \Big(L^\infty_\tau C^\alpha_x\Big),
$$
and endow this space with its natural norm. Then $\mcE^\alpha$ is continuously embedded into $\mcC^\alpha_w$. Furthermore, if $\alpha \in (0,1)$, the spaces $\mcE^\alpha, \mcC^\alpha_w$ and $\mcH^\alpha$ are equal, with equivalent norms.  }
\end{prop}

\medskip

The weighted version $\Big(L^\infty_\tau C^\alpha_x\Big)_w$ of $L^\infty_\tau C^\alpha_x$ is the same space, equipped with the norm 
$$
\|f\|_{\Big(L^\infty_\tau C^\alpha_x\Big)_w} := \sup_{0\leq \tau \leq T}\,e^{-\kappa\tau}\,\big\| f(\tau,\cdot)\big\|_{C^\alpha}.
$$
We use the following regularization properties of the heat operator associated with $L$, in the proof of global in time well-posedness for the (PAM) equation. It is proved under this form in Section 3.4 of \cite{BBF15}.

\medskip

\begin{thm}[\textbf{Schauder estimates}]
\label{ThmSchauder}   {\it 
\begin{itemize}
   \item For any choice of parameters $\beta$ and $\epsilon>0$, such that $-2+2\epsilon <\beta<0$, we have
$$ 
\big\| \mathscr{L}^{-1} (v) \big\|_{\mcC^{\beta+2 - 2\epsilon}_w} \lesssim_T \kappa^{-\epsilon} \big\|v\big\|_{\big(L^\infty_TC^\beta_x\big)_w}.
$$
   
   \item Given $\beta\in(0,2)$ and $\epsilon\in[0,1)$, we have the continuity estimate
$$ 
\| \mathscr{L}^{-1} (v) \|_{\mcC^{\beta+2-2a-2\epsilon}_w} \lesssim \kappa^{-\epsilon} \big\|v\big\|_{\mcC^\beta_w},
$$
for an implicit constant in the inequality ndependent of $\kappa$.
\end{itemize}   }
\end{thm}

\medskip

Before turning to the definition of an intertwined pair of parabolic paraproducts we close this section with another useful continuity property involving the H\"older spaces $\mcC^\sigma_\omega$; it is used in the proof of the continuity properties of the swap an merging operators, Proposition \ref{prop:R1App} in Appendix \ref{SectionAppendixContinuity}. Recall the manifold $M$ is compact. 

\medskip

\begin{prop}  \label{prop:hol}   {\it  
Given $\alpha\in(0,1)$, a space-time weight $\omega$, some integer $a\geq 0$ and a standard family $\mcP\in {\sf StGC}^a$, there exists a constant $c$ depending only on the weight $\omega$, such that 
$$
\omega(\tau)^{-1} \Big|\big(\mcP_t f\big)(e) - \big(\mcP_sf\big)(e') \Big| \leq c\, \left(s+t + \rho(e,e')^2\right)^\frac{\alpha}{2}\, \big\|f\big\|_{\mcC^\alpha_\omega},
$$
uniformly in $s,t\in(0,1]$ and $e= (\tau,x)$ and $e'=(\sigma,y)\in\mcM$, with $\tau\geq \sigma$.   }
\end{prop}

\bigskip

It is possible, and necessary for our purpose here, to make the link between the notions of regularity defined in terms of the operator $L$ and the usual notion of regularity given by the differentiable structure of the manifold, for regularity indices in the range $(1,2)$. Since the collection of vector fields $(V_i)_{1\leq i\leq \ell_0}$ spans smoothly each tangent space, for every function $f\in C^1(M)$ and every $x\in M$
$$ 
\nabla f(x) = \sum_{1\leq i\leq \ell_0} \gamma_i(x) (V_if)(x)\,V_i(x),
$$
for a collection $(\gamma_i)_{1\leq i\leq \ell_0}$ of smooth coefficients. For two points $x,y\in M$ and $f\in C^2(M)$ we have
$$ 
\big|f(x)-f(y)-\langle \nabla f(x), \pi_{x,y}\rangle_{T_xM}\big| \lesssim_f d(x,y)^2,
$$
where $T_xM$ is the canonical tangent space of $M$ at the point $x\in M$, and $\pi_{x,y}$ is a tangent vector of $T_xM$ of length $d(x,y)$, whose associated geodesic reaches $y$ at time $1$. Such a tangent vector $\pi_{xy}$ is unique if $d(x,y)$ is no greater than $r_M$, the injectivity radius of the compact manifold $M$. Combining these two facts, we have
$$ 
\Big|f(x)-f(y)- \sum_{\ell=1}^{\ell_0} \gamma_\ell(x) V_\ell(f)(x) \langle V_\ell(x), \pi_{x,y}\rangle_{T_xM}\Big| \lesssim_f d(x,y)^2.
$$
We can then define for $\alpha\in(1,2)$, the space $H^\alpha$ defined by the norm
\begin{equation}
\|f\|_{H^\alpha} := \big\|f\big\|_{H^1(M)} + \sup_{\genfrac{}{}{0pt}{}{x\neq y\in M}{d(x,y)< r_M}} \frac{\Big| f(x)-f(y)- \sum_{\ell=1}^{\ell_0} \gamma_\ell(x) V_\ell(f)(x) \langle V_\ell(x), \pi_{x,y}\rangle_{T_xM} \Big|}{d(x,y)^\alpha} <\infty.
\label{eq:Calpha}
\end{equation}
Following Proposition 2.5 in \cite{BB15}, it can be easily proved that for $\alpha\in(1,2)$ then $C^\alpha$ is continuously embedded into $H^\alpha$: Uniformly in $x,y\in M$ with $d(x,y)\leq r_M$, one has
\begin{equation}   \label{eq:Calpha1bis}
\Big| f(x)-f(y)- \sum_{\ell=1}^{\ell_0} \gamma_\ell(x) V_\ell(f)(x) \langle V_\ell(x), \pi_{x,y}\rangle_{T_xM} \Big|\lesssim d(x,y)^\alpha \|f\|_{C^\alpha}.
\end{equation}
The parabolic counterpart goes as follows, taking into account the fact that because of the parabolic scaling, a regularity of order $\alpha<2$ can be encoded in finite increments in time, with no need of a higher order expansion, and a first order expansion in space. The precise statement takes the following form.

\medskip

\begin{prop}  \label{prop:Calpha2}   {\it  
Given $\alpha\in(1,2)$ and $f\in\mcC^\alpha$, there exists a positive implicit constant such that for every $e=(\tau,x)$, $e'=(\sigma,y)$ in $\mcM$, with $\rho(e,e')\leq r_M$ then
$$ 
\Big| f(e)-f(e')- \sum_{\ell=1}^{\ell_0} \gamma_\ell(e) V_\ell(f)(e) \langle V_\ell(x), \pi_{x,y}\rangle_{T_xM} \Big| \lesssim \|f\|_{\mcC^\alpha}\,\rho(e,e')^\alpha.
$$
}
\end{prop}

\medskip

\begin{Dem}
We do not give all the details of the proof since it follows exactly the proof of Proposition 2.5 in \cite{BB15}. Here is a guideline.   \vspace{0.15cm}

We decompose $f$ at the scale $r=\rho(e,e')$ with 
$$ 
f = \big(f-\mcP_{r^2}^2f\big) + c \int_{r^2}^\infty \mcP_t \mcQ_t(f) \, \frac{dt}{t}.
$$
We plug this decomposition in the left hand side of the desired inequality. One gets estimate on the contribution of the first part $(f-\mcP_{r^2}^2f)$ using that
$$ 
\big\|f-\mcP_{r^2}^2f\big\|_\infty \lesssim r^{\alpha} \|f\|_{\mcC^\alpha}, \qquad \big\|V_i(f-\mcP_{r^2}^2f)\big\|_\infty \lesssim r^{\alpha-1} \|f\|_{\mcC^\alpha}
$$
and recalling that $\pi_{x,y}$ has length $d(x,y)\leq r$. For the second part, we integrate the contributions along $t>r^2$, applying  \eqref{eq:Calpha1bis} at order $2$ to the kernel $K_{\mcP_t}$ to obtain, for $e=(\tau,x)$ and $e'=(\sigma,y)$, the estimate
\begin{align*}
\Big|K_{\mcP_t}(e,e'')-K_{\mcP_t}(e',e'') - \sum_{\ell=1}^{\ell_0} \gamma_\ell(x) K_{V_\ell\mcP_t}(e,e'') \langle V_\ell(x), \pi_{x,y}\rangle_{T_xM}\Big| \lesssim \left(\frac{r}{\sqrt{t}}\right)^2 \mcG_t(e,e''),
\end{align*}
where $\mcG_t$ is the Gaussian kernel. Integrating this estimate against $\mcQ_t(f)$ gives a factor $\left(\frac{r}{\sqrt{t}}\right)^2 t^{\alpha/2}$ which can then be integrated along $t\in(r^2,\infty)$ since $\alpha\in(1,2)$.
\end{Dem}

\bigskip

\subsection{Parabolic paraproducts}
\label{SubsectionParaproducts}

We give here a quick presentation of the pair of intertwined paraproducts introduced in \cite{BBF15}, following the semigroup approach developed first in \cite{BB15}. The starting point for the introduction of the operator $\Pi$ is Calder\'on's reproducing formula \eqref{eq:calderon}. Using iteratively the Leibniz rule for the differentiation operators $V_i$ or $\partial_\tau$, 
we have the following decomposition 
$$ 
fg = \sum_{\mcI_b} \, a^{I,J}_{k,\ell} \int_0^1\Big(\mcA^{I,J}_{k,\ell}(f,g) + \mcA^{I,J}_{k,\ell}(g,f)\Big) \, \frac{dt}{t} + \sum_{\mcI_b} \, b^{I,J}_{k,\ell} \int_0^1 \mcB^{I,J}_{k,\ell}(f,g) \, \frac{dt}{t},
$$
where 
\begin{itemize}
   \item[\textcolor{gray}{$\bullet$}] $\mcI_b$ is the set of all tuples $(I,J,k,\ell)$ with the tuples $I,J$ and the integers $k,\ell$ satisfying the constraint 
$$
\frac{|I|+|J|}{2} + k+\ell = \frac{b}{2};
$$

   \item[\textcolor{gray}{$\bullet$}] $a^{I,J}_{k,\ell},b^{I,J}_{k,\ell}$ are bounded sequences of numerical coefficients;   \vspace{0.15cm}
   
   \item[\textcolor{gray}{$\bullet$}] for $(I,J,k,\ell)\in \mcI_b$, $\mcA^{I,J}_{k,\ell}(f,g)$ has the form
$$ 
\mcA^{I,J}_{k,\ell}(f,g) := \mcP_t^{(b)} \Big(t^{\frac{|I|}{2}+k} V_I \partial_\tau^k\Big) \left( \mcS_t^{(b/2)} f \cdot \big(t^{\frac{|J|}{2}+\ell} V_J\partial_\tau^\ell\big) \mcP_t^{(b)} g\right) 
$$  
with $\mcS^{(b/2)} \in {\sf GC}^{b/2}$;   \vspace{0.15cm}

   \item[\textcolor{gray}{$\bullet$}]  for $(I,J,k,\ell)\in \mcI_b$, $\mcB^{I,J}_{k,\ell}(f,g)$ has the form
$$ 
\mcB^{I,J}_{k,\ell}(f,g) := \mcS_t^{(b/2)} \left( \Big\{\big(t^{\frac{|I|}{2}+k} V_I\partial_\tau^k\big) \mcP_t^{(b)} f\Big\} \cdot \Big\{\big(t^{\frac{|J|}{2}+\ell} V_J\partial_\tau^\ell\big)\mcP_t^{(b)} g\Big\}\right)
$$
with $\mcS^{(b/2)} \in {\sf GC}^{b/2}$.
\end{itemize}

\medskip

\begin{definition*} 
\label{def:paraproduit} {\it
Given $f$ in $\bigcup_{s\in (0,1)} \mcC^s$ and $g\in L^\infty(\mcM)$, we define the \textbf{\textsf{paraproduct}} ${\sf P}^{(b)}_g f$ by the formula
\begin{align*}
{\sf P}^{(b)}_g f &:= \int_0^1 \Bigg\{  \sum_{\mcI_b ; \frac{|I|}{2} + k > \frac{b}{4}} \, a^{I,J}_{k,\ell} \, \mcA^{I,J}_{k,\ell}(f,g) +  \sum_{\mcI_b ; \frac{|I|}{2}+k> \frac{b}{4}} \, b^{I,J}_{k,\ell} \, \mcB^{I,J}_{k,\ell}(f,g) \Bigg\} \, \frac{dt}{t},
\end{align*}
and the \textbf{\textsf{resonant term}} $\Pi^{(b)}(f,g)$ by the formula
\begin{equation*}
\begin{split}
&\Pi^{(b)}(f,g) :=   \\
&\int_0^1 \left\{   \sum_{\mcI_b ; \frac{|I|}{2} + k \leq  \frac{b}{4}} \, a^{I,J}_{k,\ell} \Big( \mcA^{I,J}_{k,\ell}(f,g) + \mcA^{I,J}_{k,\ell}(g,f)\Big) + \sum_{\mcI_b ; \frac{|I|}{2}+k = \frac{|J|}{2}+\ell = \frac{b}{4}} \, b^{I,J}_{k,\ell} \, \mcB^{I,J}_{k,\ell}(f,g) \right\} \, \frac{dt}{t}.
\end{split}
\end{equation*}   }
\end{definition*}

\medskip

With these notations, Calder\'on's formula becomes
$$ 
fg= {\sf P}^{(b)}_g f +{\sf P}^{(b)}_f g + \Pi^{(b)}(f,g)+\Delta_{-1}(f,g)
$$
with the ``low-frequency part''
$$ 
\Delta_{-1}(f,g):= \mcP_1^{(b)}\left(\mcP_1^{(b)} f \cdot \mcP_1^{(b)} g \right).
$$

\medskip

If $b$ is chosen large enough, then all the operators involved in the paraproduct and resonant terms have a kernel pointwisely bounded by a kernel $\mcG_t$ at the right scaling. Moreover, 
\begin{enumerate}
   \item[\textsf{(a)}] the paraproduct term ${\sf P}^{(b)}_g f$ is a finite linear combination of operators of the form
$$ 
\int_0^1 \mcQ^{1\bullet}_t\Big( \mcQ^{2}_t f \cdot \mcP^1_t g\Big) \, \frac{dt}{t}
$$ 
with $\mcQ^1,\mcQ^2\in {\sf StGC}^\frac{b}{4}$, and $\mcP^1\in {\sf StGC}^{[0,b]}$,  \vspace{0.1cm}
   
   \item[\textsf{(b)}] the resonant term $\Pi^{(b)}(f,g)$ is a finite linear combination of operators of the form
$$ 
\int_0^1 \mcP^{1\bullet}_t\Big( \mcQ^{1}_t f \cdot \mcQ^{2}_t g\Big) \, \frac{dt}{t}
$$ 
with $\mcQ^1, \mcQ^2\in {\sf StGC}^\frac{b}{4}$ and $\mcP^1\in {\sf StGC}^{[0,b]}$.
\end{enumerate}

\medskip

We invite the reader to see what happens of all this when working in the flat torus with its associated Laplacian. Note also that
${\sf P}^{(b)}_f({\bf 1}) = \Pi^{(b)}(f,{\bf 1}) = 0$, and that we have the identity
$$ 
{\sf P}^{(b)}_{{\bf 1}} f = f - \mcP_{\bf 1}^{(b)}\Big(\mcP_{\bf 1}^{(b)} f\Big),
$$
as a consequence of our choice of renormalizing constant. Therefore the paraproduct with the constant function ${\bf 1}$ is equal to the identity operator, up to the strongly regularizing operator $\mcP_{\bf 1}^{(b)}\circ\mcP_{\bf 1}^{(b)}$. The regularity properties of the paraproduct and resonant operators can be described as follows; it behaves as its classical, Fourier-based, counterpart \eqref{EqBonyParaProduct}.

\medskip

\begin{prop}
\label{PropRegularityParaproduct}   {\it 
\begin{itemize}
   \item[{\bf (a)}] For every real-valued regularity exponent $\alpha,\beta$, and every positive regularity exponent $\gamma$, we have
   $$
   \big\| \Delta_{-1}(f,g)\big\|_{\mcC^\gamma} \lesssim  \|f\|_{\mcC^\alpha} \|g\|_{\mcC^\beta}.
   $$  
   
   \item[{\bf (b)}] For every $\alpha\in (-3,3)$ and $f\in \mcC^\alpha$, we have
   \begin{equation*} 
      \Big\| {\sf P}^{(b)}_g f\Big\|_{\mcC^\alpha} \lesssim \big\| g \big\|_\infty \|f\|_{\mcC^\alpha}
      \end{equation*}   
   for every $g\in L^\infty$, and 
   \begin{equation*} 
      \Big\| {\sf P}^{(b)}_g f\Big\|_{\mcC^{\alpha+\beta}} \lesssim \|g\|_{\mcC^\beta} \|f\|_{\mcC^\alpha}
      \end{equation*}
   for every $g\in\mcC^\beta$ with $\beta<0$ and $\alpha+\beta\in (-3,3)$.  \vspace{0.2cm}
   
   \item[{\bf (c)}] For every $\alpha,\beta\in (-3,3)$ with $\alpha+\beta>0$, we have the continuity estimate
$$ 
\Big\| \Pi^{(b)}(f,g) \Big\|_{\mcC^{\alpha+\beta}} \lesssim \|f\|_{\mcC^\alpha} \|g\|_{\mcC^\beta}.
$$  
\end{itemize}   }
\end{prop}

\medskip

\begin{definition*}  {\it 
We define a \textsf{\textbf{modified paraproduct}} $\widetilde{\sf P}^{(b)}$ setting 
$$
\widetilde{\sf P}^{(b)}_g f  := \mathscr{L}^{-1} \Big({\sf P}^{(b)}_g \big( \mathscr{L} f \big)\Big).
$$    }
\end{definition*}

\medskip

The next proposition shows that if one chooses the parameter $\ell_1$ that appears in the reference kernels $\mcG_t$, and the exponent $b$ in the definition of the paraproduct, both large enough, then the modified paraproduct $\widetilde{\sf P}^{(b)}$ has the same algebraic/analytic properties as ${\sf P}^{(b)}$.

\medskip

\begin{prop}
\label{PropStructureTildePi}   {\it 
\begin{itemize}
   \item[\textcolor{gray}{$\bullet$}] For a choice of large enough constants $\ell_1$ and $b$, the modified paraproduct $\widetilde {\sf P}_g f$ is a finite linear combination of operators of the form
$$ 
\int_0^1 \mcQ^{1\bullet}_t\Big( \mcQ^{2}_t f \cdot \mcP^1_t g\Big) \, \frac{dt}{t}
$$ 
with $\mcQ^1 \in {\sf GC}^{\frac{b}{8}-2}$, $\mcQ^2\in {\sf StGC}^\frac{b}{4}$ and $\mcP^1\in {\sf StGC}^{[0,b]}$.   \vspace{0.15cm}
   
   \item[\textcolor{gray}{$\bullet$}] For every $\alpha\in (-3,3)$ and $\epsilon\in (0,1)$ with $\alpha-\epsilon\in(-3,3)$ and $f\in \mcC^\alpha$, we have
\begin{equation*} 
\Big\| \widetilde{\sf P}^{(b)}_g f\Big\|_{\mcC_w^{\alpha-\epsilon}} \lesssim \kappa^{-\epsilon} \big\| w^{-1}g\big\|_\infty \|f\|_{\mcC^\alpha},
\end{equation*}
for every $g\in L^\infty$.
\end{itemize}   }
\end{prop}

\medskip

The notation $\mcQ^{1\bullet}$ does not make sense for a generic element $\mcQ^1$ of ${\sf GC}$. The operator that appears in this formula is actually of the form $\mcQ^1_t = \mcQ_tt^{-1}\mathscr{L}^{-1}$, with $\mcQ\in {\sf StGC}$, for which the notation $\mcQ^{1\bullet}$ makes sense. Note that the norm $\|f\|_{\mcC^\alpha}$ above has no weight. Note here the normalization identity 
$$
\widetilde{\sf P}^{(b)}_{\bf 1}f = f - \mathscr{L}^{-1}\circ\mcP^{(b)}_{\bf 1}\circ\mcP^{(b)}_{\bf 1}(\mathscr{L} f)
$$
for every distribution in $f\in\mcS'_o$; it reduces to 
$$
\widetilde{\sf P}^{(b)}_{\bf 1}f = f - \mcP^{(b)}_{\bf 1}\mcP^{(b)}_{\bf 1}(f)
$$
if $f_{\big| \tau =0} = 0$.

\medskip

Following the definition of the \textsf{\textbf{inner difference operator}} $\mathscr{D}$ given in Subsection \ref{SubSectionIteratedParaPdcts}, we extend it to the parabolic setting defining $\mathscr{D}\,\big(=\mathscr{D}_e\big)$ by the formula
$$
\iint_{\mcM^2} \big(\mathscr{D} f\big)(e') g(e)\,\nu(de)\nu(de') := \iint_{\mcM^2} \big(f(e')-f(e)\big)g(e)\,\nu(de)\nu(de');
$$
with this notation, the crucial motivating relation
$$
{\sf P}_f\Big(\widetilde{\sf P}_a g\Big) - {\sf P}_{fa} g = {\sf P}_f\Big(\widetilde{\sf P}_{\mathscr{D}a} g\Big)
$$
holds true.

\medskip

Last, we prove an elementary property of the modified paraproduct that provides some pointwise information on the solutions to singular partial differential equations. 

\medskip 

\begin{prop} 
\label{prop:derive}   {\it 
Let $\alpha,\beta$ be positive regularity exponents, and let $u,Z\in \mcC^\alpha$ with $v\in \mcC^\beta$ be given, with $Z(0,\cdot)=0$. Assume that 
$$ 
u-\widetilde {\sf P}_v Z \in \mcC^{\alpha+\beta},
$$
and define $\gamma:=\min(\alpha+\beta,1)$. If $\alpha+\beta \neq 1$, we have
$$ 
\left| u(e)-u(e') - v(e)\big(Z(e)-Z(e')\big)\right|\lesssim  \rho(e,e')^\gamma,
$$
uniformly in $e,e'\in \mcM$ with $\rho(e,e')\leq 1$. If $\alpha +\beta = 1$, we have an additional logarithmic loss 
$$ 
\left| u(e)-u(e') - v(e)\big(Z(e)-Z(e')\big)\right|\lesssim  \rho(e,e') \log\Big(1+\rho(e,e')^{-1}\Big).
$$   }
\end{prop}

\medskip

\begin{Dem}
Due to the assumption, one has
$$ \left|u(e)-u(e') - v(e)\big(Z(e)-Z(e')\big)\right| \lesssim \rho(e,e')^\beta +(\star)$$
with
$$
(\star):= \left| \big(\widetilde{\sf P}_vZ\big)(e) - \big(\widetilde{\sf P}_vZ\big)(e') - v(e)\big(Z(e)-Z(e')\big)\right|.
$$
Using Calder\'on reproducing formula, or the normalization, yields 
$$
\widetilde {\sf P}_1 Z=Z
$$ 
since $Z(0,\cdot)=0$, and we see that $(\star)$ is equal to
$$ 
\left| \int_0^1 \mcQ_t^{\bullet}\big( \mcQ_t Z \mcP_tv \big)(e)-\mcQ_t^\bullet\big( \mcQ_t Z \mcP_tv \big)(e') - v(e) \mcQ_t^\bullet\big( \mcQ_t Z \big)(e) + v(e) \mcQ_t^\bullet\big( \mcQ_t Z \big)(e') \, \frac{dt}{t} \right|,
$$   
so
$$ 
(\star) \lesssim  \int_0^1 \left| \int \left(K_{\mcQ_t^\bullet}(e,a) - K_{\mcQ_t^\bullet}(e',a)\right)  \mcQ_t Z(a) \big(\mcP_tv(a)- v(e)\big) \nu(da) \right|\, \frac{dt}{t}.
$$
Using the regularity estimates on $v$ and on the kernel of the approximation operators, one sees that
\begin{align*}
 (\star) & \lesssim  \|v\|_{\mcC^\beta} \int_0^1  \int \min\left\{1,\frac{\rho(e,e')}{\sqrt{t}}\right\} \mcG_t(e,a)  \left|\mcQ_t Z(a)\right| \left(t+\rho(a,e)^2\right)^{\gamma/4}  \, \nu(da)\frac{dt}{t}  \\
  & \lesssim \|v\|_{\mcC^\beta} \|Z\|_{\mcC^\alpha} \int_0^{\rho(e,e')^2}  t^{(2\alpha+\gamma)/4} \,\frac{dt}{t}  + \|v\|_{\mcC^\beta} \|Z\|_{\mcC^\alpha} \int_{\rho^2}^1  \int \frac{\rho(e,e')}{\sqrt{t}} t^{(2\alpha+\gamma)/4} \,\frac{dt}{t} \\
  & \lesssim \|v\|_{\mcC^\beta} \|Z\|_{\mcC^\alpha} \rho(e,e')^{\gamma},
\end{align*}
which concludes the proof.
\end{Dem}

\bigskip

\section[\hspace{0.7cm} Paracontrolled expansion formula]{Paracontrolled expansion formula}
\label{SectionAppendixTaylor}

We give in this section a detailed and rigorous proof of Theorem \ref{ThmTaylorExpansion}. The parameter $b$ is fixed, and we note $\Pi$ for $\Pi^{(b)}$.

\medskip

\begin{thm}[\textsf{\textbf{Higher order Taylor expansion}}]   
\label{ThmTaylorExpansion-bis}   
{\it 
Let $f : \bfR\mapsto\bfR$ be a $C^4$ function, and let $u$ be a real-valued and $\mcC^\alpha$ function on $\mcM$, with $\alpha\in(0,1)$. Then
\begin{equation}
\label{EqTaylorExpansion-bis}
\begin{split}
f(u) &= {\sf P}_{f'(u)}(u) + \frac{1}{2}\,\left\{{\sf P}_{f^{(2)}(u)}(u^2) - 2{\sf P}_{f^{(2)}(u)u}(u) \right\}   \\
&+ \frac{1}{3!}\Big\{ {\sf P}_{f^{(3)}(u)}(u^3) - 3{\sf P}_{f^{(3)}(u)u}(u^2) + 3{\sf P}_{f^{(3)}(u)u^2}(u) \Big\} + f(u)^\sharp
\end{split}
\end{equation}
for some remainder $f(u)^\sharp\in \mcC^{4\alpha}$. If moreover $f$ is of class $C^5$, then the remainder term $f(u)^\sharp$ is locally-Lipschitz with respect to $u$, in the sense that 
\begin{equation}\label{EqStabilityTaylor}
\big\|f(u)^\sharp - f(v)^\sharp\big\|_{\mcC^{4\alpha}} \lesssim \big(1+\|u\|_{\mcC^\alpha}+\|v\|_{\mcC^\alpha}\big)^{4}\, \|u-v\|_{\mcC^\alpha}.
\end{equation}   }
\end{thm}

\medskip

\begin{Dem}
Let us give a detailed proof of the third order expansion, that claims that 
$$
(\star) := f(u) - {\sf P}_{f'(u)}(u) - \frac{1}{2}\,\left\{{\sf P}_{f^{(2)}(u)}(u^2) - 2{\sf P}_{f^{(2)}(u)u}(u) \right\} 
$$
is a $3\alpha$-H\"older function. We invite the reader to follow what comes next in the light of the proof given in Section \ref{SectionTaylor} in the time-independent, flat, model setting of the torus.

\ssk

As, by definition, the paraproduct operator ${\sf P}_{g}(\cdot)$ is a finite sum of different terms, each of them of the form
$$
\mcA^1_g(\cdot) := \int_0^1  \mcQ^{1\bullet}_t\Big( \mcQ^2_t(\cdot) \, \mcP^1_t(g) \Big) \,\frac{dt}{t},
$$
with $\mcQ^{1},\mcQ^2$ at least to ${\sf StGC}^3$, it is sufficient to prove that the following function
\begin{align*}
(\star) := f(u) - & \int_0^1  \left[\mcQ^{1\bullet}_t\Big( \mcQ^2_t(u) \, \mcP^1_t\big(f'(u)\big) \Big)  + \frac{1}{2} \mcQ^{1\bullet}_t\Big( \mcQ^2_t(u^2) \, \mcP^1_t\big(f^{(2)}(u)\big) \Big) \right. \\
& \qquad \qquad \left. - \mcQ^{1\bullet}_t\Big( \mcQ^2_t(u) \, \mcP^1_t\big(f^{(2)}(u)u\big) \Big)  \right] \,\frac{dt}{t}
\end{align*}
is an element of $\mcC^{3\alpha}$. Using Calder\'on's reproducing formula together with the normalization of the paraproduct, we have
$$ 
f(u) \simeq \int_0^1 \mcQ^{1\bullet}_t \mcQ^2_t\big(f(u) \mcP^1_{t}(1)\big) \, \frac{dt}{t}
$$
up to a remainder quantity corresponding to the low frequency part that is as smooth as we want. So one can write $(\star)$ under the form
\begin{equation} 
\label{eq:star}
(\star)=\int_0^1 \mcQ^{1\bullet}_t(\epsilon_t) \,\frac{dt}{t},
\end{equation}
with 
\begin{align*} 
\epsilon_t:= &  \mcQ^2_t\Big(f(u)\Big)\mcP^1_{t}(1) - \mcQ^2_t(u) \, \mcP^1_t\Big(f'(u)\Big)   \\
& - \frac{1}{2}  \mcQ^2_t(u^2) \, \mcP^1_t\Big(f^{(2)}(u)\Big)  +  \mcQ^2_t(u) \, \mcP^1_t\Big(f^{(2)}(u)u\Big).
\end{align*}
Due to the orthogonality/cancellation property of the operators $\mcQ^{1\bullet}_t$, it suffices for us to get an $L^\infty$ control of $\epsilon_t$.
Using the kernel representation of the different operators, we have for every $e\in \mcM$
\begin{align*}
\epsilon_t(e) = \iint_{\mcM^2} & K_{\mcQ^2_t}(e,e') K_{\mcP^1_t}(e,e'') \Big\{  f\big(u(e')\big) -u(e') f'\big(u(e'')\big)   \\
& - \frac{1}{2}u^2(e')f^{(2)}\big(u(e'')\big) + u(e')f^{(2)}\big(u(e'')\big) u(e'')\Big\} \, \nu(de')\nu(de'')
\end{align*}
Note also that we have from the usual Tayor formula for $f$
\begin{align*}
& f\big(u(e')\big) -u(e') f'\big(u(e'')\big) - \frac{1}{2}u^2(e')f^{(2)}\big(u(e'')\big) + u(e')f^{(2)}\big(u(e'')\big) u(e'')   \\
& = \iiint_{[0,1]^3} f^{(3)}\Big(u(e'')+s_3s_2s_1 \big(u(e')-u(e'')\big)\Big) s_2s_1\big(u(e')-u(e'')\big)^3\, ds_3ds_2ds_1   \\
& \ + f(u(e''))+u(e'')f'\big(u(e'')\big)+\frac{1}{2}u^2(e'')f^{(2)}\big(u(e'')\big).
\end{align*}
When we integrate against $K_{\mcQ^2_t}(e,e') K_{\mcP^1_t}(e,e'')$ a quantity depending only in $e''$ has no contribution, since the latter kernel satisfies a cancellation property along the $e'$-variable. So we have exactly
\begin{align*}
& \epsilon_t(e) = \iint_{\mcM^2} K_{\mcQ^2_t}(e,e') K_{\mcP^1_t}(e,e'')   \\
& \left(\ \iiint_{[0,1]^3} f^{(3)}\Big(u(e'')+s_3s_2s_1\big(u(e')-u(e'')\big)\Big) s_2s_1\big(u(e')-u(e'')\big)^3 \,  ds_3ds_2ds_1\right)\nu(de')\nu(de'').
\end{align*}
Since $K_{\mcQ^2_t}$ and $K_{\mcP^1_t}$ are both pointwisely dominated by the Gaussian kernel $\mcG_{t}$, and using the fact that $f^{(3)}$ is bounded on the range of $u$, we obtain the uniform control
\begin{align*}
& \big|\epsilon_t(e)\big| \lesssim \iint_{\mcM^2} \mcG_t(e,e') \mcG_t(e,e'') \big(u(e')-u(e'')\big)^3 \, \nu(de')\nu(de'') \\
& \lesssim \|u\|_{\mcC^\alpha}^3 \, t^{3\alpha/2},
\end{align*}
from which the fact that $(\star)$ belongs to $\mcC^{3\alpha}$ follows from \eqref{eq:star}. We used for that purpose the identity
$$
u(e')-u(e'') = \big(u(e')-u(e)\big) + \big(u(e)-u(e'')\big),
$$
together with Proposition \ref{prop:equivalence} on the characterization of parabolic regularity in terms of increments, to see that
$$ 
\left|u(e')-u(e'')\right| \lesssim \big(d(e',e)+d(e'',e)\big)^\alpha \|f\|_{\mcC^\alpha}.
$$
The fourth order expansion of the statement is proved by a very similar reasoning left to the reader. One proceeds exactly as in the proof of Theorem \ref{ThmTaylorExpansion} to prove the stability estimate \eqref{EqStabilityTaylor}.
\end{Dem}

\bigskip

Observe the fact that one can give a paracontrolled expansion formula with the $\widetilde{\sf P}$ operator in place of the ${\sf P}$ operator, as a consequence of Proposition \ref{PropStructureTildePi} on the structure of the modified paraproduct $\widetilde{\sf P}$, and the proof of Theorem \ref{ThmTaylorExpansion-bis}.

\bigskip

\section[\hspace{0.7cm} Continuity results]{Continuity results}
\label{SectionAppendixContinuity}

Recall the definitions of the corrector
$$
{\sf C}(f,g,h) := {\sf \Pi}\left(\widetilde {\sf P}_{f}(g),h\right) - f\,{\sf \Pi}(g,h),
$$
the (modified) commutators
\begin{equation*}
\begin{split}
&{\sf D}(f,g,h) = {\sf \Pi}\left(\widetilde {\sf P}_{f}(g),h\right) - {\sf P}_f\Big({\sf \Pi}(g,h)\Big),   \\
&{\sf R}(f,g,h) = {\sf P}_f\Big(\widetilde{\sf P}_gh\Big) - {\sf P}_{fg}h,   \\
&{\sf S}(f,g,h) = {\sf P}_f\left( \widetilde {\sf P}_gh\right) - {\sf P}_g\left({\sf P}_fh\right),
\end{split}
\end{equation*}
and their iterates, introduced in Section \ref{SectionToolKit}. We use also their $^\circ$ variants, built with ${\sf P}$ in the role of $\widetilde{\sf P}$; we have for instance 
$$
{\sf R}^\circ(f,g,h) = {\sf P}_f\Big({\sf P}_gh\Big) - {\sf P}_{fg}h.
$$
All these operators are initially defined on the space of smooth functions. We prove in this section the continuity results on these operators stated in Section \ref{SectionToolKit}.

\medskip

Let $\chi$ stand for a smooth non-negative function on $[0,\infty)$, equal to $1$ in a neighbourhood of $0$, with $\chi(r)=0$ for $r\geq r_M$. Given $e=(t,x)\in\mcM$, set 
$$
\delta_\ell(e,e') := \chi\big(d(x',x)\big)\,\big\langle V_\ell(x), \pi_{xx'}\big\rangle_{T_xM},
$$
for $e'=(t',x')\in\mcM$ and $\exp_x(\pi_{xx'})=x'$, and $1\leq \ell\leq\ell_0$, and 
$$
{\sf \Pi}_{(1)}^\ell(g,h) := {\sf \Pi}\big({\sf P}_{\delta_\ell(e,\cdot)}g,h\big)(e).
$$

\medskip

\begin{defn}   \label{DefnRefinedCorrector}
The \textbf{\textsf{refined paraproduct}} is defined by the formula
$$
{\sf C}_{(1)}(f,g,h)(e) := {\sf C}(f,g,h)(e) - \sum_{\ell=1}^{\ell_0} \gamma_\ell (V_\ell f)\,{\sf \Pi}_{(1)}^\ell(g,h)(e).
$$
\end{defn} 

\medskip

Recall from Section 3.1.1 the simple definition of the refined corrector in the model setting of the flat torus.

\bigskip

\subsection[\hspace{-0.3cm} Boundedness of commutators/correctors]{Boundedness of commutators/correctors}
\label{Subsectionapp1}

We start by looking at the case of the swap and merging operators {\sf S} and {\sf R}. We do not emphasize in the next statement the choice of parameter $\kappa$ in the time weight. This has no consequence on the use of these continuity results for the study of singular PDEs, as we only use Schauder estimates in weighted spaces to deal with the terms from the enhancement of the noise in the study of the linear (PAM) equation, not for all the well-defined terms builts from the corrector, commutator and their iterates.

\medskip

\begin{prop} \label{prop:R1App}
{\it
\begin{itemize}
   \item[\textcolor{gray}{$\bullet$}] Let $\alpha,\beta,\gamma$ be H\"older regularity exponents with $\gamma\in(-3,3), \beta\in(0,1)$ and $\alpha\in (-\infty,0)$. Then if 
$$
\beta+\gamma<3, \qquad \textrm{ and } \qquad \delta:=\alpha+\beta+\gamma\in(-3,3),
$$ 
we have
\begin{equation}
\label{eq:R1}
\big\| {\sf S}(f,g,h) \big\|_{\mcC^{\delta}} +\big\| {\sf S}^\circ (f,g,h) \big\|_{\mcC^{\delta}} \lesssim \|f\|_{\mcC^\alpha} \, \|g\|_{\mcC^\beta} \, \|h\|_{\mcC^\gamma},
\end{equation}
so the modified commutator on paraproducts extends naturally into a trilinear continuous map from $\mcC^\alpha\times\mcC^\beta\times\mcC^\gamma$ to $\mcC^\delta$.   \vspace{0.15cm}
   
   \item[\textcolor{gray}{$\bullet$}] If $\alpha=0$ then the product $fg$ has a sense for $f\in L^\infty(\mcM)$ and $g\in \mcC^\beta$, with $0<\beta<1$, and we have 
\begin{equation}
\label{eq:R1bis}
\big\| {\sf R}^\circ(f,g,h) \big\|_{\mcC^{\beta+\gamma}} + \big\| {\sf R}(f,g,h) \big\|_{\mcC^{\beta+\gamma}} \lesssim \\|f\|_{L^\infty}\,|g\|_{\mcC^\beta} \, \|h\|_{\mcC^\gamma}.
\end{equation} 

   \item[\textcolor{gray}{$\bullet$}] If $\alpha,\beta\in(0,1)$ and $\gamma\in(-3,3)$, then we have 
\begin{equation}
\label{eq:R1ter}
\big\| {\sf R}^\circ(f,g,h) \big\|_{\mcC^{\alpha+\beta+\gamma}} \lesssim \|f\|_{\mcC^\alpha} \, \|g\|_{\mcC^\beta} \, \|h\|_{\mcC^\nu}.
\end{equation} 
\end{itemize}   }
\end{prop}

\medskip

\begin{Dem}
Recall the operators ${\sf P}^{(b)}_g(\cdot)$, respectively $ \widetilde{\sf P}^{(b)}_g(\cdot)$, are given by a finite sum of operators of the form
$$
\mcA^1_f(\cdot) := \int_0^1  \mcQ^{1\bullet}_t\Big( \mcQ^2_t(\cdot) \, \mcP^1_t(f) \Big) \,\frac{dt}{t},
$$
respectively
$$ 
\widetilde \mcA^1_f(\cdot) := \int_0^1  \widetilde \mcQ^{1\bullet}_t\Big( \widetilde \mcQ^2_t(\cdot) \, \mcP^1_t(f) \Big) \,\frac{dt}{t},
$$
where $\mcQ^{1},\mcQ^2,\widetilde \mcQ^2$ belong at least to ${\sf StGC}^3$ and $\widetilde \mcQ^1$ is an element of ${\sf GC}^3$. We describe similarly the operator ${\sf P}^{(b)}_f(\cdot)$ as a finite sum of operators of the form
$$ 
\mcA^2_f(\cdot) := \int_0^1  \mcQ^{3\bullet}_t\Big( \mcQ^4_t(\cdot)  \mcP^2_t(f) \Big) \,\frac{dt}{t}.
$$
Thus, we need to study a generic modified commutator
$$ 
\mcA^2_f\left(\widetilde \mcA^1_g(h)\right) - \mcA^1_g\left(\mcA^2_f(h)\right),
$$
and introduce for that purpose the intermediate quantity
$$ 
\mcE(f,g,h) := \int_0^1 \mcQ^{3\bullet}_s\Big( \mcQ^4_s(h) \cdot \mcP^1_s(g) \cdot \mcP^2_s(f) \Big) \,\frac{ds}{s}.  
$$
(We proceed similarly for the study of ${\sf S}^\circ$.) Note here that due to the normalization ${\sf P}_{\bf 1} \simeq \textrm{Id}$, up to some strongly regularizing operator, there is no loss of generality in assuming that 
\begin{equation} 
\label{eq:norma}
\int_0^1 \widetilde \mcQ^{1\bullet}_t  \widetilde \mcQ^2_t \, \frac{dt}{t} = \int_0^1\mcQ^{1\bullet}_t \mcQ^2_t \, \frac{dt}{t} = \int_0^1\mcQ^{3\bullet}_t \mcQ^4_t \, \frac{dt}{t} = \textrm{Id}. 
\end{equation}

\medskip

\noindent \textsf{\textbf{Step 1. Study of $\mcA^2_f\left(\widetilde \mcA^1_g(h)\right) - \mcE(f,g,h)$.}} We shall use a family $\mcQ$ in ${\sf StGC}^a$, for some $a>|\delta|$, to control the H\"older norm of that quantity. By definition, and using the normalization \eqref{eq:norma}, for every $r\in(0,1)$, the quantity $\mcQ_r \Big(\mcA^2_f\left(\widetilde \mcA^1_g(h)\right) - \mcE(f,g,h) \Big)$ is equal to
{\small \begin{align*}
&\int_0^1\int_0^1  \mcQ_r \mcQ^{3\bullet}_s\Big\{ \mcQ^4_s \widetilde \mcQ^{1\bullet}_t\Big( \widetilde \mcQ^2_t(h)  \mcP^1_t(g) \Big) \cdot \mcP^2_s(f)  \Big\} \,\frac{ds\,dt}{st} - \int_0^1 \mcQ_r \mcQ^{3\bullet}_s\Big( \mcQ^4_s(h) \cdot \mcP^1_s(g) \cdot \mcP^2_s(f) \Big) \,\frac{ds}{s} \\
&= \int_0^1\int_0^1  \mcQ_r \mcQ^{3\bullet}_s\Big\{ \mcQ^4_s \widetilde \mcQ^{1\bullet}_t\Big( \widetilde \mcQ^2_t(h)  \big(\mcP^1_t(g)-\mcP^1_s(g)\big) \Big) \cdot  \mcP^2_s(f)  \Big\} \,\frac{dsdt}{st},
\end{align*}   }
where in the last line the variable of $\mcP^1_s(g)$ is that of $\mcQ^{3\bullet}_s$, and so it is frozen through the action of $\widetilde \mcQ^4_s\mcQ^{1\bullet}_t$. Then using that $g\in \mcC^\beta$ with $\beta \in (0,1)$, we  know by Proposition \ref{prop:hol} that we have, for $\tau\geq \sigma$,
$$ 
\Big|\big(\mcP^1_s g\big)(x,\tau) - \big(\mcP^1_t g\big)(y,\sigma)\Big| \lesssim \left(s+t+ \rho\big((x,\tau),(y,\sigma)\big)^2\right)^\frac{\alpha}{2} \|g\|_{\mcC^\beta}.
$$
Note that it follows from equation \eqref{EqIteratedG} that the kernel of $\mcQ^4_s \widetilde \mcQ^{*1}_t$ is pointwisely bounded by $\mcG_{t+s}$, and allowing different constants in the definition of $\mcG$, we have
\begin{equation}
\label{eq:gg}
\mcG_{t+s}\big((x,\tau),(y,\sigma)\big) \left(s+t+ d(x,y)^2\right)^\frac{\beta}{2} \lesssim (s+t)^\frac{\beta}{2}\, \mcG_{t+s}\big((x,\tau),(y,\sigma)\big).
\end{equation}
So using the cancellation property of the operators $\mcQ$, resp. $\mcQ^i$ and $\widetilde \mcQ^i$, at an order no less than $a$, resp. $3$, we deduce that
\begin{align*}
& \left\| \mcQ_r \Big(\mcA^2_f\left(\widetilde \mcA^1_g(h)\right) - \mcE(f,g,h)\Big) \right\|_\infty  \\
& \qquad \lesssim \|f\|_{\mcC^\alpha}\|g\|_{\mcC^\beta} \|h\|_{\mcC^\gamma} \int_0^1\int_0^1 \left(\frac{sr}{(s+r)^2}\right)^\frac{a}{2} \left(\frac{st}{(s+t)^2}\right)^\frac{3}{2} t^\frac{\gamma}{2} (s+t)^\frac{\beta}{2} s^\frac{\alpha}{2}   \,\frac{ds\,dt}{st},
\end{align*}
where we used that $\alpha$ is negative to control $\mcP^2_s(f) $. The integral over $t\in(0,1)$ can be computed since $\gamma>-3$ and $\beta+\gamma<3$, and we have
\begin{align*}
& \left\| \mcQ_r \Big(\mcA^2_f\left(\widetilde \mcA^1_g(h)\right) - \mcE(f,g,h)\Big) \right\|_\infty  \\
& \qquad \lesssim \|f\|_{\mcC^\alpha} \|g\|_{\mcC^\beta} \|h\|_{\mcC^\gamma} \int_0^1\int_0^1 \left(\frac{sr}{(s+r)^2}\right)^\frac{a}{2}  s^\frac{\delta}{2} \,\frac{ds}{s} \\
& \qquad \lesssim \|f\|_{\mcC^\alpha} \|g\|_{\mcC^\beta} \|h\|_{\mcC^\gamma} r^\frac{\delta}{2},
\end{align*}
uniformly in $r\in(0,1)$ because $|a|>\delta$. That concludes the estimate for the high frequency part. We repeat the same reasoning for the low-frequency part by replacing $\mcQ_r$ with $\mcQ_1$ and conclude that
$$ 
\Big\| \mcA^2_f\left(\widetilde \mcA^1_g(h)\right) - \mcE(f,g,h) \Big\|_{\mcC^{\delta}} \lesssim \|f\|_{\mcC^\alpha} \|g\|_{\mcC^\beta} \|h\|_{\mcC^\gamma}.
$$

\medskip

\textsf{\textbf{Step 2. Study of $\mcA^1_g\left(\mcA^2_f(h)\right) - \mcE(f,g,h)$.}} This term is almost the same as that of Step 1 and can be treated in exactly the same way. Note that $\mcQ_r \Big(\mcA^1_g\left(\mcA^2_f(h)\right) - \mcE(f,g,h) \Big)$ is equal, for every $r\in (0,1)$, to 
{\small \begin{align*}
&\int_0^1\int_0^1  \mcQ_r \mcQ^{1\bullet}_t\Big( \mcQ^2_t\mcQ^{3\bullet}_s\big( \mcQ^4_s(h)  \mcP^2_s(f)  \big) \cdot  \mcP^1_t(g)  \Big) \,\frac{ds\,dt}{st} - \int_0^1 \mcQ_r\mcQ^{3\bullet}_s\Big(\mcQ^4_s(h) \cdot \mcP^1_s(g) \cdot \mcP^2_s(f)  \Big) \,\frac{ds}{s} \\
&= \int_0^1\int_0^1 \mcQ_r\mcQ^{1\bullet}_t\Big\{ \mcQ^2_t\mcQ^{3\bullet}_s\Big( \mcQ^4_s(h) \big(\mcP^1_t(g)-\mcP^1_s(g)\big)  \cdot \mcP^2_s(f)  \Big) \Big\} \,\frac{ds\,dt}{st},
\end{align*}   }
where in the last line the variable of $\mcP^1_t(g)$ is that of $\mcQ^{1\bullet}_t$, so it is frozen through the action of $\mcQ^{3\bullet}_s$. The same proof as in Step 1 can be repeated, which gives the first statement of the theorem.

\medskip

\textsf{\textbf{Step 3. Proof of the second statement.}} For the second statement, Step 1 still holds. So it only remains to compare $\mcE(f,g,h)$ with $\mcA^2_{fg}(h)$. This amounts to compare $\mcP_t^2(fg)$ with $\mcP_t^1(fg \mcP_t^2(f)$. Using the regularity of $g\in \mcC^\beta$ and the uniform boundedness of $f\in L^\infty$, we get
$$ 
\left\| \mcP_t^2(fg) - \mcP_t^1(g) \mcP_t^2(uf) \right\|_{L^\infty} \lesssim t^{\beta/2}
$$
which allows us to conclude.

\medskip

\textsf{\textbf{Step 4. Proof of the third statement.}} The key observation is that ${\sf R}^\circ$ satisfies more cancellation that ${\sf R}$, namely
$$ 
{\sf R}^\circ({\bf 1},g,h)= {\sf R}^\circ(f,{\bf 1},h) =0.
$$
Taking advantage of that fact, for any $\mcQ \in {\sf StGC}^3$, we have for $s\in(0,1)$ and $e\in \mcM$, the identity
$$ 
\mcQ_s\big({\sf R}^\circ(f,g,h)\big)(e) = \mcQ_s\Big({\sf R}^\circ\big(f-f(e),g-g(e),h\big)\Big)(e).
$$
The difference structure has been taken into account in this cancellation/identity and we can now estimate each piece in the definition of ${\sf R}^\circ$ separately. Using previous arguments, one has
\begin{align*}
& \Big|\mcQ_s\big({\sf P}_{f-f(e)} {\sf P}_{g-g(e)}(h)\big)(e)\Big| \\
& \lesssim \int_{[0,1]^3} \iint_{\mcM^2} \left(\frac{st_1}{(s+t_1)^2}\right)^3 \left(\frac{t_1t_2}{(t_1+t_2)^2}\right)^3 \left(\frac{t_2 t_3}{(t_2+t_3)^2}\right)^3 \mcG_{s+t_1}(e,e') \mcG_{t_1+t_2}(e,e'') \\
& \qquad \qquad \qquad \big|f(e')-f(e)\big| \big|g(e)-g(e'')\big|  \nu(de') \nu(de'') t_3^{\gamma/2} \|h\|_{\mcC^\gamma} \, \frac{dt_1}{t_1} \frac{dt_2}{t_2} \frac{dt_3}{t_3} \\
& \lesssim_{f,g,h} \int_{[0,1]^3} \iint_{\mcM^2} \left(\frac{st_1}{(s+t_1)^2}\right)^3 \left(\frac{t_1t_2}{(t_1+t_2)^2}\right)^3 \left(\frac{t_2 t_3}{(t_2+t_3)^2}\right)^3 \mcG_{s+t_1}(e,e') \mcG_{t_1+t_2}(e,e'') \\
& \qquad \qquad \qquad \rho(e,e')^{\alpha} \rho(e,e'')^\beta \nu(de') \nu(de'') t_3^{\gamma/2}  \, \frac{dt_1}{t_1} \frac{dt_2}{t_2} \frac{dt_3}{t_3}   \\
& \lesssim_{f,g,h} \int_{[0,1]^3} \left(\frac{st_1}{(s+t_1)^2}\right)^3 \left(\frac{t_1t_2}{(t_1+t_2)^2}\right)^3 \left(\frac{t_2 t_3}{(t_2+t_3)^2}\right)^3 (s+t_1+t_2)^{(\alpha+\beta)/2} t_3^{\gamma/2}  \, \frac{dt_1}{t_1} \frac{dt_2}{t_2} \frac{dt_3}{t_3}    \\
& \lesssim s^{(\alpha+\beta+\gamma)/2} \|f\|_{\mcC^\alpha} \|g\|_{\mcC^\beta} \|h\|_{\mcC^\gamma},
\end{align*}
uniformly in $e\in \mcM$ and $s\in (0,1)$, which concludes the proof of the desired estimate for this first term. (The intermediate implicit constants in the upper bounds are constant multiples of $\|f\|_{\mcC^\alpha} \|g\|_{\mcC^\beta} \|h\|_{\mcC^\gamma}$.) The second estimate is obtained similarly, observing that
$$
\Big|\big(f(e')-f(e)\big)\big(g(e')-g(e)\big)\Big| \lesssim \rho(e',e)^{\alpha+\beta} \|f\|_{\mcC^\alpha} \|g\|_{\mcC^\beta},
$$
and gives in the end
$$ 
\big\|\mcQ_s\big({\sf R}^\circ(f,g,h)\big)\big\|_{L^\infty} \lesssim s^{(\alpha+\beta+\gamma)/2} \,\|f\|_{\mcC^\alpha} \|g\|_{\mcC^\beta} \|h\|_{\mcC^\gamma}.
$$
\end{Dem}

\bigskip

\begin{rem} \label{rem:R1} {\it 
The above proof actually shows the following property of the operator
\begin{equation} 
{\sf S}_{f,h} := g \mapsto {\sf S}(f,g,h)
\label{eq:opT}
\end{equation}
where $f\in \mcC^\alpha$ and $h\in\mcC^\gamma$ are fixed. For all families $\mcQ^1, \mcQ^2 \in {\sf GC}^a$ for some $a>0$, the linear operator $\mcQ_t^1  {\sf S}_{f,h} \mcQ_s^{2\bullet}$ has a kernel pointwisely bounded by
$$ 
(t+s)^\frac{\beta+\gamma}{2} \, \left(\frac{st}{(s+t)^2}\right)^\frac{a}{2} \mcG_{t+s}\big(e,e'\big)\, \|f\|_{\mcC^\alpha}\,\|h\|_{\mcC^\gamma}.
$$   }
\end{rem}

\medskip

\begin{prop}  \label{prop:C1}   {\it
\begin{itemize}
   \item[\textcolor{gray}{$\bullet$}] Let $\alpha,\beta,\gamma$ be H\"older regularity exponents with $\alpha\in(0,1), \beta\in(-3,3)$ and $\gamma\in (-\infty,3]$. Set 
$$
\delta := (\alpha+\beta)\wedge 3 +\gamma.
$$ 
If
$$
0 < \alpha + \beta + \gamma < 1 \qquad \textrm{ and }\qquad \beta + \gamma < 0
$$
then the corrector ${\sf C}$ extends continuously into a trilinear map from $\mcC^\alpha \times\mcC^\beta \times\mcC^\gamma$ to $\mcC^\delta$.   \vspace{0.15cm}
   
   \item[\textcolor{gray}{$\bullet$}] If $\alpha,\beta,\gamma$ are positive then the commutator ${\sf D}$ is a continuous trilinear map from $\mcC^\alpha \times\mcC^\beta \times\mcC^\gamma$ to $\mcC^\delta$.
\end{itemize}   }
\end{prop}

\medskip

\begin{Dem} 
The result on ${\sf C}$ was already proved in \cite[Proposition 3.6]{BB15} in a more general setting. We only focus here on proving the boundedness of ${\sf D}$. As already done above, we represent the operator ${\sf P}^{(b)}_f(\cdot)$ under the form
$$
\mcA_f(\cdot) := \int_0^1  \mcQ^{1\bullet}_t\Big( \mcQ^2_t(\cdot) \, \mcP^1_t(f) \Big) \,\frac{dt}{t},
$$
and the resonant term $\Pi^{(b)}(g,h)$ as
$$ 
\mcB(g,h) := \int_0^1  \mcP^{2\bullet}_t\Big( \mcQ^3_t(g)  \mcQ^4_t(h) \Big) \,\frac{dt}{t}.
$$
Thus, we need to study a generic modified commutator
\begin{align*}
& (\star):=\mcB\big(\mcA_f(g),h\big) - \mcA_{f}\big(\mcB(g,h)\big)   \\
& = \int_0^1\int_0^1   \mcP^{2\bullet}_t\Big( \mcQ^3_t \mcQ^{1\bullet}_s\big( \mcQ^2_s(g) \, \mcP^1_t(f) \big) \,  \mcQ^4_t(h) \Big) \,\frac{ds}{s}\frac{dt}{t}   \\
&-\int_0^1\int_0^1 \mcQ^{1\bullet}_s\Big( \mcQ^2_s\mcP^{2\bullet}_t\big( \mcQ^3_t(g)  \mcQ^4_t(h) \big) \, \mcP^1_s(f) \Big) \,\frac{ds}{s}\frac{dt}{t},
\end{align*}
and introduce for that purpose the intermediate quantity
$$ 
\mcE(f,g,h) := \int_0^1 \mcP^{2\bullet}_t\Big( \mcP_t^1(f) \mcQ^3_t(g) \mcQ_t^4(h)\Big) \,\frac{dt}{t}.  
$$
Then we compare the two quantities with $\mcE(f,g,h)$, such as done previously. Each of these two comparisons makes appear an exact commutation on the function $f$, due to our choice of normalization for our paraproducts. Using the $\mcC^\alpha$ regularity on $f$ together with the cancellation property of the $\mcQ$ operators, we get 
\begin{align*}
\|\mcQ_r(\star)\|_{L^\infty} & \lesssim \int_0^1\int_0^1 \left(\frac{r}{r+t}\right)^{3} \left(\frac{st}{(s+t)^2}\right)^3 s^{\beta/2} t^{\gamma/2} (s+t)^{\alpha/2}\, \frac{dt}{t}\frac{ds}{s} \\
& + \int_0^1 \int_0^1 \left(\frac{rs}{(r+s)^2}\right)^{3} \left(\frac{s}{s+t}\right)^3 t^{\beta/2} t^{\gamma/2} (s+t)^{\alpha/2} \, \frac{dt}{t} \frac{ds}{s}\\
& \lesssim \int_0^1 \left(\frac{r}{r+t}\right)^{3} t^{(\alpha+\beta+\gamma)/2} \, \frac{dt}{t} + \int_0^1 \left(\frac{r}{r+t}\right)^3 t^{\beta/2} t^{\gamma/2} (r+t)^{\alpha/2} \, \frac{dt}{t} \\
& \lesssim r^{\delta/2},
\end{align*}
which shows that $(\star)$ belongs to $\mcC^{\delta}$.
\end{Dem}

\bigskip

\subsection[\hspace{-0.3cm} Boundedness of iterated commutators/correctors]{Boundedness of iterated commutators/correctors}
\label{Subsectionapp2}

We now turn to the study of the continuity properties of the iterated versions of commutators/correctors, and start with the (modified) iterated commutator on paraproducts.

\medskip

\begin{prop} \label{prop:commuR}   {\it
 If $\gamma\in(0,1)$, $\alpha,\nu\in(0,1/2)$ and $\beta\in(-3,3)$ then we have 
\begin{equation}
\label{eq:RPi-high}
\big\| {\sf R}^\circ(u ,f,{\sf P}_a g) - {\sf P}_a {\sf R}^\circ (u,f,g) \big\|_{\mcC^{\alpha+\beta+\gamma+\nu}} \lesssim \|f\|_{\mcC^\alpha} \, \|g\|_{\mcC^\beta} \, \|u\|_{\mcC^\nu} \|a\|_{\mcC^\gamma}.
\end{equation}   }
\end{prop}

\medskip

This statement is a combination of both \eqref{eq:R1} and \eqref{eq:R1ter} of Proposition \ref{prop:R1App}; we let the reader write the proof. We now state the continuity result for the $4$-linear iterated swap operator defined in \eqref{Eq4LinearSwap}; a similar continuity result holds for the $5$-linear iterated operator defined in \eqref{Eq5LinearSwap}; its proof is left to the reader.

\medskip

\begin{prop} \label{prop:R2} \label{PropR2}   {\it 
Let $\alpha,\beta_1, \beta_2,\gamma$ be H\"older regularity exponents with $\gamma\in(-3,3)$, $\beta_1,\beta_2 \in(0,1)$ and $\alpha\in (-\infty,0)$. Then if 
$$
\alpha+\beta_2+\gamma<3, \qquad \textrm{ and } \qquad \delta:=\alpha+\beta+\gamma+\nu\in(-3,3),
$$ 
we have
\begin{equation}
\label{eq:R2}
\big\| {\sf S}\big(f,(g_1,g_2),h\big) \big\|_{\mcC^{\delta}} \lesssim \|f\|_{\mcC^\alpha} \, \|g_1\|_{\mcC^\beta_1} \, \|g_2\|_{\mcC^\beta_2} \|h\|_{\mcC^\gamma},
\end{equation}
so the commutator defines a quadrilinear continuous map from $\mcC^\alpha\times\mcC^{\beta_1}\times\mcC^{\beta_2}\times \mcC^\gamma$ to $\mcC^\delta$.   }
\end{prop}

\medskip

\begin{Dem}
Fix some functions $f \in \mcC^\alpha$ and $g_2\in C^{\beta_2}$; we have
$$
{\sf S}\big(f,(g_1,g_2),h\big) := {\sf S}\Big(f,\widetilde {\sf P}_{g_1} g_2, h\Big) - {\sf P}_{g_1}{\sf S}(f,g_2,h).
$$
With the same notations as in the proof of Proposition \ref{prop:R1App}, for which we have relations \eqref{eq:norma}, we write 
\begin{align*}
{\sf P}_{g_1} \big[{\sf S}(f,g_2,h)\big] &= \int_0^1 \mcQ^{1\bullet}_t\Big(\mcQ_t^2 \big[{\sf S}(f,g_2,h)\big] \cdot \mcP^1_t g_1\Big) \,\frac{dt}{t}   \\
&= \int_0^1 \int_0^1 \mcQ^{1\bullet}_t \Big(\mcQ_t^2 \big[{\sf S}(f,\widetilde \mcQ_s^{1\bullet} \widetilde \mcQ_s^2 g_2,h)\big]  \cdot \mcP^1_t g_1\Big) \,\frac{ds}{s} \frac{dt}{t}.
\end{align*}
Expanding ${\sf S}\big(f,\widetilde {\sf P}_{g_1} g_2, h\big)$ correspondingly, we get with ${\sf S}_{f,h}$ defined in \eqref{eq:opT},
\begin{align}
\label{eq:decompR2}
{\sf S}\big(f,(g_1,g_2),h\big) = \int_0^1 \int_0^1 \mcQ_t^{1\bullet}\Big\{ \mcQ_t^2 {\sf S}_{f,h}\Big(\widetilde \mcQ_s^{1\bullet} \big( \widetilde \mcQ_s^2 g_2 \cdot \big(\mcP_t^1 g_1 - \mcP_s^1 g_1\big) \big)\Big)
\Big\} \,\frac{ds}{s}\frac{dt}{t}, 
\end{align}
where the variable of $\mcP^1_tg_1h$ is that of $\mcQ^{1\bullet}_t$. Since $g_1$ belongs to $\mcC^{\beta_1}$, with $\beta_1 \in (0,1)$, we know from Proposition \ref{prop:hol} that 
$$
\Big|\big(\mcP^1_t g_1\big)(e) - \big(\mcP^1_s g_1\big)(e')\Big| \lesssim \left(t+s+\rho(e,e')^2\right)^\frac{\beta_1}{2} \|g_1\|_{\mcC^{\beta_1}},
$$
for all $e,e' \in \mcM$. As above, fix a collection $\mcQ$ of ${\sf StGC}^a$, for some $a>3$, to control H\"older norms. We need to estimate
$$ 
\Big\| \mcQ_r {\sf S}\big(f,(g_1,g_2),h\big) \Big\|_{L^\infty(\mcM)}.
$$
Using decomposition \eqref{eq:decompR2} , we have
\begin{equation} 
\label{qrr2}
\Big\|  \mcQ_r {\sf S}\big(f,(g_1,g_2),h\big) \Big\|_{L^\infty(\mcM)} \lesssim \int_0^1 \int_0^1 \left(\frac{rt}{(r+t)^2}\right)^\frac{a}{2} I_{s,t} \,\frac{ds}{s}\frac{dt}{t},
\end{equation}
where
$$ 
I_{s,t} := \sup_{e\in \mcM} \, \mcQ_t^2 {\sf S}_{f,h}\Big(\widetilde \mcQ_s^{1\bullet} \big( \widetilde \mcQ_s^2 g_2 \cdot \big(\mcP_t^1 g_1(e) - \mcP_s^1 g_1\big) \big)\Big)(e).
$$
Due to Remark \ref{rem:R1}, we have a pointwise estimate of the kernel of the operator $\mcQ_t^2{\sf S}_{f,h}\big( \mcQ_s^{1\bullet}(\cdot)\big)$, so with the pointwise regularity estimate on $h$ and \eqref{eq:gg}, we deduce that
\begin{align*}
 I_{s,t} & \lesssim (s+t)^\frac{\alpha+\beta_1+\gamma}{2} \,\big\|\widetilde \mcQ_s^2 g_2\big\|_{L^\infty} \,\|f\|_{\mcC^\alpha}\, \|g_1\|_{\mcC^{\beta_1}} \, \|h\|_{\mcC^\gamma} \\
 & \lesssim (s+t)^\frac{\delta}{2} \,\|f\|_{\mcC^\alpha} \, \|g_1\|_{\mcC^{\beta_1}} \, \|g_2\|_{\mcC^{\beta_2}} \, \|h\|_{\mcC^\gamma}.
\end{align*}
It follows from that estimate and the fact that $|\sigma|<a$, that
$$ 
\Big\|  \mcQ_r {\sf S}\big(f,(g_1,g_2),h\big) \Big\|_{L^\infty(\mcM)} \lesssim r^\frac{\delta}{2} \, \|f\|_{\mcC^\alpha} \, \|g_1\|_{\mcC^{\beta_1}} \, \|g_2\|_{\mcC^{\beta_2}} \, \|h\|_{\mcC^\gamma},
$$
uniformly in $r\in(0,1)$. A similar analysis of the low frequency of ${\sf S}\big(f,(g_1,g_2),h\big)$ can be done, which completes the proof of the H\"older estimate.
\end{Dem}

\medskip

\begin{prop} \label{prop:R}
{\it Pick $\beta_1,\beta_2 \in(0,1)$ with $\beta_1+\beta_2<1$, and $\gamma\in(-3,3)$. We have
\begin{align}
&\big\| {\sf R}(f, \widetilde{\sf P}_{g_1}g_2,h) - {\sf P}_{g_1}{\sf R}(f,g_2,h) \big\|_{\mcC^{\beta_1+\beta_2+\gamma}} + \big\| {\sf R}^\circ(f , {\sf P}_{g_1}g_2,h) - {\sf P}_{g_1} {\sf R}^\circ(f,g_2,h) \big\|_{\mcC^{\beta_1+\beta_2+\gamma}} \nonumber \\
& \qquad \qquad \lesssim \|f\|_{L^\infty} \, \|g_1\|_{\mcC^{\beta_1}} \, \|g_2\|_{\mcC^{\beta_2}} \|h\|_{\mcC^\gamma}. \label{eq:RPi2}
\end{align}   }
\end{prop}

\ssk

\begin{Dem} 
Fix a parameter $r\in(0,1)$ and look for a control (uniformly in $e\in\mcM$) of 
$$ 
(\star) := \mcQ_{r} \big[  {\sf R}(f , \widetilde {\sf P}_{g_1}g_2,h) - {\sf P}_{g_1} {\sf R}(f,g_2,h) \big](e).
$$
We follow the arguments of Proposition \ref{prop:R1App} -- more precisely of \eqref{eq:R1bis} and \eqref{eq:R1ter}, since
$$ 
(\star) = \mcQ_{r} \big[{\sf R}(f , \widetilde {\sf P}_{g_1}g_2,h) \big](e)  - \int_0^1 \mcQ_r\mcQ_s^{\bullet} \big[ \mcP_s g_1 \cdot \mcQ_s [{\sf R}(f,g_2,h)] \big](e),
$$
and we have seen there how tocontrol the composition of a $\mcQ$ operator with ${\sf R}$. So we repeat the exact same reasoning with replacing the function $f$ by $\widetilde {\sf P}_{g_1}f$ and now the commutator (with the extra paraproduct ${\sf P}_{g_1}$) brings terms of the form
$$ 
\big[{\sf P}_{g_1}g_2(e') - {\sf P}_{g_1}g_2(e'')\big] - \mcP_r(g_1)(e''')\big[ g_2(e')-g_2(e'')\big]
$$
for points $e',e'',e'''$ of the parabolic space $\mcM$, multiplied by Gaussian kernels localizing points at the scale $\rho(e,e')+\rho(e,e'')+\rho(e,e''')\lesssim r^{1/2}$.
Using Proposition \ref{prop:derive}, we deduce that since $g_1\in \mcC^{\beta_1}$ and $g_2\in \mcC^{\beta_2}$ then
$$ 
\left|\big[{\sf P}_{g_1}g_2(e') - {\sf P}_{g_1}g_2(e'')\big] - \mcP_r(g_1)(e''')\big[ g_2(e')-g_2(e'')\big]\right| \lesssim r^{(\beta_1+\beta_2)/2}.
$$
The result follows by repeating the computations done in proving \eqref{eq:R1ter}.
\end{Dem}

\medskip

We now look at the iterated corrector. The proof of continuity for the lower and upper iterates are almost the same and the reader can see clearly on the model case of iterated integrals what the difference is.

\medskip

\begin{prop} \label{prop:C2}
{\sf 
Let $\alpha,\beta_1\in(0,1), \beta_2\in(-3,3)$ and $\gamma\in (-\infty,3]$. Assume that
$\alpha+\beta_1+\beta_2<3$ with
$$
\delta:=\alpha + \beta_1 + \beta_2 + \gamma \in (0,1),\qquad \alpha+\beta_2 + \gamma < 0 \quad \textrm{and} \quad \beta_1+\beta_2+\gamma<0.
$$
Then the $4$-linear upper iterated corrector ${\sf C}$ is a continuous map from $\mcC^\alpha \times\mcC^{\beta_1}\times\mcC^{\beta_2} \times\mcC^{\gamma}$ to $\mcC^\delta$.}
\end{prop}

\ssk

\begin{Dem}
Fix $f \in \mcC^\alpha$ and $h \in \mcC^{\gamma}$ and set
$$
\overline {\sf C}(\cdot) := {\sf C}(f,\cdot,h),
$$
so 
$$ 
{\sf C}\big(f,(g_1,g_2),h\big) = \overline{\sf C}\Big(\widetilde{\sf P}_{g_1}g_2\Big) - g_1\,\overline{\sf C}(g_2).
$$
Using the same notation as previously, and omitting for convenience the indices on the different collections $\mcQ$ and $\mcP$, we write 
\begin{align*}
\overline {\sf C}\Big(\widetilde {\sf P}_{g_1}g_2\Big) &= \int_0^1 \overline {\sf C} \widetilde \mcQ^{\bullet}_s\Big(\widetilde\mcQ_s g_2 \cdot \mcP_s g_1\Big) \,\frac{ds}{s},   \\
g_1\,\overline {\sf C}(g_2)  &= g_1 \,\overline {\sf C}\Big(\widetilde {\sf P}_{\bf 1} (g_2)\Big) =  g_1 \,\int_0^1 \overline {\sf C} \widetilde \mcQ^{\bullet}_s \Big( \widetilde  \mcQ_s g_2 \cdot \mcP_s {\bf 1}\Big) \,\frac{ds}{s}.
\end{align*}
Note that due to the conservation property of the heat semigroup associated with $L$, the quantity $\mcP_s {\bf 1}$ is either constant equal to $1$ or to $0$, depending on whether $\mcP_s$ encodes some cancellation or not. Thus, given $e=(x,\tau)\in\mcM$, and setting
$$
F_{s,e}:= \widetilde \mcQ_s g_2 \cdot \big(\mcP_s g_1- \mcP_s({\bf 1}) \cdot g_2(e)\big),
$$ 
we have
\begin{align*}
 {\sf C}\big(f,(g_1,g_2),h\big)(e)= \overline {\sf C}\Big( \widetilde {\sf P}_{g_1}g_2\Big)(e) - _1(e)\,\overline {\sf C}(g_2)(e) = \int_0^1 \overline {\sf C} \Big(\widetilde \mcQ_s^\bullet F_{s,e}\Big) (e) \,\frac{ds}{s}.
\end{align*}
As before, we can use that $g_1\in \mcC^{\beta_1}$, with $\beta_1\in (0,1)$. We have for $e,e' \in \mcM$
$$
\big|g_1(e) - g_1(e')\big| \lesssim \rho(e,e')^{\beta} \,\|g_1\|_{\mcC^{\beta_1}},
$$
and therefore, using the ``Gaussian bounds'' for $\mcP_s$, 
\begin{align*}	
\big|\big(\mcP_s g_1\big)(e') - \big(\mcP_s{\bf 1}\big)(e')\, g_1(e)\big| \lesssim \big(s+\rho(e,e')^2 \big)^\frac{\beta}{2} \|g_1\|_{\mcC^{\beta_1}}. 
\end{align*}
As done in the proof of Proposition \ref{prop:C1}, we introduce an intermediate quantity of the form
\begin{align*}
S\big(a,b,c\big) := \int_0^1 \mcP_t\Big(\mcQ_tb \cdot \mcQ_t c \cdot \mcP_t a\Big) \,\frac{dt}{t},
\end{align*}
and write
\begin{align} \label{eq:I1-I2} \nonumber
\overline {\sf C}\Big(\widetilde \mcQ_s^\bullet F_{s,e}\Big)(e) 
&= {\sf \Pi}\big(\widetilde {\sf P}_f(\widetilde \mcQ_s^\bullet F_{s,e}),h\Big) (e) - S\Big(f,\widetilde \mcQ_s^\bullet F_{s,e},h\Big)(e)   \\ \nonumber 
& \qquad + S\Big(f,\widetilde \mcQ_s^\bullet F_{s,e},h\Big)(e) - f(e) \cdot {\sf \Pi}\big(\widetilde\mcQ_s^\bullet F_{s,e},h\Big)(e)   \\
& =: I_1(s) + I_2(s). 
\end{align}

\ssk

\textcolor{gray}{$\bullet$} We start with the estimate for $I_2$. One can then write with generic notations for the resonant term $\Pi$
\begin{align*}
\Big(S\big(f,g_2,h\big) - f\cdot {\sf \Pi}(g_2,h)\Big) (e) = \int_0^1 \mcP_t\Big( \mcQ_tg_2 \cdot \mcQ_t h \cdot \big(\mcP_t f -f(e) \big) \Big) (e) \,\frac{dt}{t},
\end{align*}
and it is known that the integrand is pointwisely bounded by $t^\frac{\alpha+\nu_1+\nu_2}{2}$. Since this argument only uses pointwise estimates, we can replace $b$ by $\mcQ_s^\bullet F_{s,e}$. Therefore, by writing 
\begin{align*}
\int_0^1 I_2(s) \,\frac{ds}{s} = \int_0^1 \int_0^1 \mcP_t \Big(\mcQ_t \widetilde \mcQ_s^\bullet F_{s,e} \cdot \mcQ_t h \cdot \big(\mcP_t f -f(e) \big) \Big) (e) \,\frac{dt}{t} \frac{ds}{s}
\end{align*}
and using 
\begin{equation} 
\label{eq:QsQt}
\big\| \mcQ_t \widetilde\mcQ_s^\bullet \phi \big\|_{L^\infty(\mcM)} \lesssim \left(\frac{st}{(s+t)^2}\right)^{3/2} \big\|\phi \big\|_{L^\infty(\mcM)},
\end{equation}
with $\phi=F_{s,e}$, we obtain 
\begin{align*}
&\left\| \int_0^1 I_2(s) \,\frac{ds}{s} \right\|_{L^\infty(\mcM)}   \\
& \leq \int_0^1 \int_0^1 \left\|e \mapsto \mcP_t\big(\mcQ_t\widetilde \mcQ_s^\bullet F_{s,e} \cdot \mcQ_t h \cdot (f(e)-\mcP_t f) \big) (e) \right\|_{L^\infty} \,\frac{dt}{t} \frac{ds}{s}   \\
& \lesssim \|g_1\|_{\mcC^{\beta_1}}\|g_2\|_{\mcC^{\beta_2}} \|f\|_{\mcC^\alpha} \|h\|_{\mcC^{\gamma}}   \\ 
&\quad \times \int_0^1 \int_0^1 \left(\frac{st}{(s+t)^2}\right)^\frac{3}{2} \mcG_{t+s}(e,e') \Big(s+\rho( e,e')^2 \Big)^\frac{\beta_1}{2}  s^{\beta_2/2} t^\frac{\alpha+\gamma}{2}  \,\frac{ds}{s}\frac{dt}{t}   \\
& \lesssim  \|f\|_{\mcC^\alpha} \|g_1\|_{\mcC^{\beta_1}}\|g_2\|_{\mcC^{\beta_2}}\|h\|_{\mcC^{\gamma}}
\int_0^1 \int_0^1 \left(\frac{st}{(s+t)^2}\right)^\frac{3}{2} s^{\beta_2/2} (s+t)^{\beta_1/2} t^\frac{\alpha+\gamma}{2}  \,\frac{ds}{s}\frac{dt}{t}   \\
&\lesssim  \|f\|_{\mcC^\alpha} \|g_1\|_{\mcC^{\beta_1}}\|g_2\|_{\mcC^{\beta_2}}\|h\|_{\mcC^{\gamma}},
\end{align*}
since $\alpha+\beta_1+\beta_2+\gamma>0$. 

\ssk

\textcolor{gray}{$\bullet$}  Let us now estimate the regularity of $I_2(s)$. Let $e,e' \in \mcM$ with $\rho(e,e')\leq 1$. We split the integral in $t$ into two parts, corresponding to $t<\rho(e,e')^2$ or $t>\rho(e,e')^2$. In the first case, note that
\begin{align*}	
\int_0^{\rho(e,e')^2} t^{(\alpha+\beta_1+\beta_2+\gamma)/2} \,\frac{dt}{t} \lesssim \rho(e,e')^{\alpha+\beta_1+\beta_2+\gamma},
\end{align*}
so that by repeating the arguments above, we get the desired estimate. In the case $t>\rho^2$ with $\rho:=\rho(e,e')$, write for $s \in (0,1)$
\begin{align} \label{eq:C4-split1} \nonumber 
	&\int_{\rho^2}^1 \left\{ \mcP_t\Big(\mcQ_t \widetilde \mcQ_s^\bullet F_{s,e} \cdot \mcQ_t h \cdot \big(f(e)-\mcP_tf\big) \Big)(e) \right.\\ \nonumber 
	&\qquad\qquad\qquad\qquad- \left.\mcP_t\Big(\mcQ_t \widetilde\mcQ_s^\bullet F_{s,e'} \cdot \mcQ_t h \cdot\big(f(e')-\mcP_tf\big) \Big)(e') \right\} \,\frac{dt}{t} \\ \nonumber 
& = \int_{\rho^2}^1 \left\{ \mcP_t\Big(\mcQ_t \widetilde \mcQ_s^\bullet F_{s,e} \cdot \mcQ_t h \cdot \big(f(e)-\mcP_tf\big) \Big)(e) \right. \\ \nonumber 
&\qquad\qquad\qquad\qquad- \left. \mcP_t\Big(\mcQ_t \widetilde \mcQ_s^\bullet F_{s,e} \cdot \mcQ_t h \cdot \big(f(e)-\mcP_tf\big) \Big)(e')\right\} \,\frac{dt}{t}\\ \nonumber 
&\qquad + \big(g_1(e)-g_1(e')\big)\int_{\rho^2}^1 \mcP_t\Big(\mcQ_t \widetilde \mcQ_s^\bullet \widetilde\mcQ_s b \cdot Q_th \cdot \big(f(e')-\mcP_t f \big) \Big)(e') \,\frac{dt}{t}\\
&\qquad - \big(f(e)-f(e')\big)\int_{\rho^2}^1\mcP_t\Big(\mcQ_t \widetilde \mcQ_s^\bullet F_{s,e} \cdot \mcQ_t h\Big)(e')\,\frac{dt}{t}.
\end{align}
For the second and third term, we can assume $s \simeq t$ by \eqref{eq:QsQt}. One obtains
\begin{align*}	
	&\big|g_1(e) - g_1(e')\big| \int_{\rho^2}^1 \Big|\mcP_t\Big(\mcQ_t \widetilde \mcQ_s^\bullet \widetilde \mcQ_s g_2 \cdot \mcQ_th \cdot (f(e')-\mcP_t f)\Big)(e')\Big|\,\frac{dt}{t}\\
	& \qquad  \lesssim \|f\|_{\mcC^\alpha} \|g_1\|_{\mcC^{\beta_1}}\|g_2\|_{\mcC^{\beta_2}}\|h\|_{\mcC^{\gamma}}\, \rho^{\beta_1} 
	\int_{\rho^2}^1 t^\frac{\alpha+\beta_1+\gamma}{2} \,\frac{dt}{t} \\
	& \qquad  \lesssim \|f\|_{\mcC^\alpha} \|g_1\|_{\mcC^{\beta_1}}\|g_2\|_{\mcC^{\beta_2}}\|h\|_{\mcC^{\gamma}}\, \rho^{\alpha+\beta_1+\beta_2+\gamma},
\end{align*}
since $\alpha+\beta_2+\gamma$ is negative, and
\begin{align*}
	& \big| f(e)-f(e') \big|\int_{\rho^2}^1\Big|\mcP_t\Big(\mcQ_t \widetilde \mcQ_s^\bullet F_{s,e} \cdot \mcQ_t h \Big)(e')\Big|\,\frac{dt}{t}\\
	 & \qquad \lesssim \|f\|_{\mcC^\alpha} \|g_1\|_{\mcC^{\beta_1}}\|g_2\|_{\mcC^{\beta_2}}\|h\|_{\mcC^{\gamma}}\, \rho^{\alpha}  \int_{\rho^2}^1 t^\frac{\beta_1+\beta_2+\gamma}{2} \,\frac{dt}{t}\\
	& \qquad \lesssim \|f\|_{\mcC^\alpha} \|g_1\|_{\mcC^{\beta_1}}\|g_2\|_{\mcC^{\beta_2}}\|h\|_{\mcC^{\gamma}}\, \rho^{\alpha+\beta_1+\beta_2+\gamma},
\end{align*}
since $\beta_1+\beta_2+\gamma$ is also negative. For the first term in \eqref{eq:C4-split1}, we now repeat the arguments of the proof of Proposition \ref{prop:C1}, which rely on the Lipschitz regularity of the heat kernel as well as the fact that $(\alpha+\beta_1+\beta_2+\gamma) \in (0,1)$. Summarising, we have shown that for $e,e' \in \mcM$ with $\rho(e,e') \leq 1$
\begin{align*}
&\left|\int_0^1 \Big( I_2(s)(e) - I_2(s)(e')\Big) \,\frac{ds}{s} \right|   \\
&\qquad\quad  \lesssim  \rho(e,e')^{\alpha+\beta_1+\beta_2+\gamma}\, \|f\|_{\mcC^\alpha} \|g_1\|_{\mcC^{\beta_1}}\|g_2\|_{\mcC^{\beta_2}}\|h\|_{\mcC^{\gamma}}.
\end{align*}
Let us now come to $I_1(s)$ as defined in \eqref{eq:I1-I2}. Set $\phi:=\widetilde \mcQ_s^\bullet F_{s,e}$, and write
\begin{align*}
\Big| {\sf \Pi}\big(\widetilde {\sf P}_f(\phi),h\big) - S\big(f,g_2,h\big) \Big| \leq \int_0^1 \big|\mcP_t (A_t(\phi,f) \cdot \mcQ_t h) \big| \,\frac{dt}{t}
\end{align*}
with 
\begin{align*}
A_t(\phi,f) := \mcQ_t\left( \int_0^1 \mcP_t \widetilde \mcQ_r^\bullet \big( \widetilde\mcQ_r \phi \cdot \mcP_r f\big) \,\frac{dr}{r} - \mcP_t f \mcP_t \phi\right).
\end{align*}
Following the proof of Proposition \ref{prop:C1}, and using  \eqref{eq:QsQt}, one obtains 
\begin{align*}
	& \Big\| A_t\big(\widetilde \mcQ_s^\bullet F_{s,e},u\big)\Big\|_{L^\infty(\mcM)} \\
	& \qquad \lesssim \int_0^1 \left(\frac{r t}{(r+t)^2}\right)^\frac{3}{2} \left(\frac{sr}{(s+r)^2}\right)^\frac{3}{2}  s^\frac{\beta_2}{2}  (r+t)^\frac{\alpha+\beta_1}{2}  \,\frac{dr}{r}\|f\|_{\mcC^\alpha} \|g_1\|_{\mcC^{\beta_1}} \|g_2\|_{\mcC^{\beta_2}},
\end{align*}
hence
\begin{align*}
	& \left\| \int_0^1 I_1(s) \,\frac{ds}{s} \right\|_{L^\infty(\mcM)} 
	\lesssim \|f\|_{\mcC^\alpha} \|g_1\|_{\mcC^{\beta_1}}\|g_2\|_{\mcC^{\beta_2}}\|h\|_{\mcC^{\gamma}}   \\
	& \qquad \times  \int_0^1 \int_0^1 \int_0^1  \left(\frac{r t}{(r+t)^2}\right)^\frac{3}{2} \left(\frac{sr}{(s+r)^2}\right)^{3/2} s^\frac{\beta_2}{2}  (r+t)^\frac{\alpha+\beta_1}{2} t^\frac{\gamma}{2} \, \frac{dr}{r} \frac{ds}{s}\frac{dt}{t},
\end{align*}
and the triple integral is finite since  $(\alpha+\beta_1+\beta_2+\gamma)$ is positive.  

\ssk

\textcolor{gray}{$\bullet$}  For the regularity estimate of $I_1(s)$, consider 
\begin{align*}
	\int_0^1 \Big\{ \mcP_t \Big(A_t(\widetilde \mcQ_s^\bullet F_{s,e},f) \cdot \mcQ_t h\Big)(e) - \mcP_t \Big(A_t(\widetilde \mcQ_s^\bullet F_{s,e'},f) \cdot \mcQ_t h\Big)(e') \Big\} \,\frac{dt}{t}.
\end{align*}
The estimate of this expression is similar, though simpler, compared to the one of $I_2(s)$, as here $e$ is frozen only in one spot. As before, one deals with this terms using the heat kernel regularity of $\mcP_t$ and the regularity estimate for $a$.
\end{Dem}

\bigskip
\bigskip

\noindent \textcolor{gray}{$\bullet$} {\sf I. Bailleul} - {\small Univ Rennes, CNRS, IRMAR - UMR 6625, F-35000 Rennes, France.} {\it ismael.bailleul@univ-rennes1.fr}   \vspace{0.3cm}

\noindent \textcolor{gray}{$\bullet$} {\sf F. Bernicot} - {\small Laboratoire de Math\'ematiques Jean Leray, CNRS - Université de Nantes, 2 Rue de la Houssini\`ere 44322 Nantes Cedex 03, France.}  {\it frederic.bernicot@univ-nantes.fr}

\end{document}